\documentclass[a4paper,11pt]{article}

\usepackage{macros}
\usepackage{vmargin}

\setmarginsrb{2cm}{1.5cm}{2cm}{1cm}{0mm}{1cm}{0mm}{0.5cm}

\long\def\symbolfootnote[#1]#2{\begingroup%
\def\thefootnote{\fnsymbol{footnote}}\footnote[#1]{#2}\endgroup}

\graphicspath{{figures/}}

\begin{document}
\pagestyle{empty}

\begin{center}{\LARGE\bf Mixed timestepping schemes for nonsmooth mechanics with high frequency damping\symbolfootnote[2]{This is a preprint of a paper submitted to Multibody System Dynamics.}}\end{center} 

\begin{center}
\parbox{3in}{ \centering \textbf{Shahed Rezaei}\\
  Institute of Applied Mechanics\\
	RWTH Aachen\\
	Mies-van-der-Rohe-Str. 1\\
	52074 Aachen, Germany\\
  {\tt\small shahed.rezaei@rwth-aachen.de}
}
\parbox{3in}{ \centering \textbf{Thorsten Schindler}\\
  Institute of Applied Mechanics\\
	Technische Universit\"at M\"unchen\\
	Boltzmannstra\ss{}e 15\\
	85748 Garching, Germany\\
  {\tt\small thorsten.schindler@mytum.de}
}
\end{center}
\section*{Abstract}
This work deals with the integration of nonsmooth flexible multibody systems with impacts and dry friction. We develop a framework which improves a non-impulsive trajectory of state variables by impulsive correction after each time-step if necessary. This correction is automatic and is evaluated on the same kinematic level as the piecewise non-impulsive trajectory. The resulting overall mixed timestepping scheme is consistent with respect to impacts and friction as well as benefits from advantages of the base integration schemes used to calculate the approximation inside the time-step. Therefore, we compare the generalized-$\alpha$ method, the Bathe method and the ED-$\alpha$ method.%
\section*{Keywords}
timestepping scheme $\cdot$ nonsmooth mechanics $\cdot$ flexible multibody system $\cdot$ high frequency damping $\cdot$ friction $\cdot$ index reduction 
\section{Introduction}
Nonsmooth mechanical systems are characterized by jumps in the system's velocities or accelerations. Typical examples with impacts or dry friction can e.g. be found in automotive, railway or robotics applications~\cite{Bro99,Glo01,Aca08,Lei08,Pfe08,Joh08,Ste11}. Hence, we deal with impulsive and non-impulsive periods. There are two types of integration methods to handle nonsmooth motion: event-driven schemes and timestepping schemes. On the one hand, event-driven schemes, i.e., with event detection, are highly recommendable in non-impulsive periods because of their high integration order. The drawback of event-driven schemes is that they cannot represent an infinite accumulation of impacts. On the other hand, timestepping schemes are applicable in impulsive periods but they are of integration order one in both impulsive and non-impulsive periods~\cite{Aca08}. Integration schemes known from computational mechanics~\cite{Lau02,Doy11,Sim92,Bau10,Deu08,Woh11} usually suffer from oscillations in the relative contact velocities at least if impulsive periods occur~\cite{Kra12}.\par
In~\cite{Sch12,Sch14a,Sch14b}, classic timestepping schemes~\cite{Mor99,Jea99,Pao02,Pao02a} are improved by splitting the non-impulsive and impulsive forces, such that the advantages of event-driven and timestepping schemes are conserved and the disadvantages are avoided. The approach is theoretically based on time-discontinuous Galerkin methods. Hence, we assume jumps for the velocity across discretization intervals, such that non-impulsive forces benefit from higher order trial functions and impulsive reactions yield local integration order one automatically. For the piecewise linear velocity trial functions discussed, an effective algorithm and framework are presented. In the case of finite element applications, it is shown however that the base integration scheme may get unstable without additional damping because of its half-explicit nature. That is why, we derive different timestepping methods based on the general framework of impulsive corrections and non-impulsive base integration schemes but leave the time-discontinuous Galerkin setting. We use base integration schemes like the generalized-$\alpha$ method~\cite{Chu93}, the Bathe-method~\cite{Bat07} and the ED-$\alpha$ method~\cite{Bot97} within mixed timestepping schemes~\cite{Ese13,Aca11a} for high-frequency damping instead. A comparison between acceleration and velocity level~\cite{Arn07,Fueh91,Eic98,Sch12d} for multi-contact cases is studied in~\cite{Sch14a,Sch14b}. We evaluate both the non-impulsive forces and the impulsive reactions on the advanced velocity level using the projection formulation for convex sets together with semi-smooth Newton methods~\cite{Ala91,Chr98,Qi99,Sch11a}.\par
Timestepping methods for structural dynamics are evaluated among other criteria based on two well understood and appreciated characteristics: unconditional stability and high frequency dissipation~\cite{New59,Hil77,Chu93,Sim13,Bau10}. The first criterion ensures that the stability of the method does not depend on the time-step size. The other requirement concerns the algorithmic dissipation and damping properties for non-physical high-frequency modes. In structural dynamics, large scale movements and lower modes are frequently of particular interest whereas high frequencies are not considered or only get into the system as spurious oscillations due to the discretization. However, the high-frequency modes affect the numerical and convergence properties of the system especially in the nonlinear regime. Hence, it is desirable to reduce these frequencies by numerical damping. The intensity of the damping should be controllable and should affect the lower modes as less as possible. In order to meet the above mentioned criteria, a variety of integration methods has been developed and optimized with regard to the desired properties. The HHT-formalism~\cite{Hil77} and the generalized-$\alpha$-method~\cite{Chu93} are two well known methods based on the Newmark-scheme~\cite{New59}. Having a controlling parameter for high frequency damping and second-order accuracy for the generalized coordinates and the Lagrange multipliers mark the generalized-$\alpha$ method as a good choice for numerical integration~\cite{Arn07}. An extended version of the generalized-$\alpha$ method for nonsmooth problems is presented in~\cite{Che13,Brue14}. However, the global accuracy remains order one in the presence of constraint forces. An undesirable property is that the damping parameters have to be selected and acceptable values depend on the characteristics of the problem being solved. The damping parameter in Newmark-based integrators plays a very important role as the high frequency oscillations may effect the results and even convergence of the finite element solver~\cite{Bot04,Sch14a,Sch14b}. Recently, new attempts have been made to overcome such problems, among those, we focus on two schemes: the Bathe-method and the ED-$\alpha$ method, which both use extra information in each time interval. The Bathe-method combines the use of the trapezoidal rule and the Euler backward method~\cite{Bat07}. It has been shown that the additional midpoint information can result in a more robust simulation in the nonlinear analysis and can avoid nonphysical large contact forces. The Bathe-method is also effective in the linear analysis~\cite{Bat12}. The ED-$\alpha$ (Energy Decaying) scheme is proved to be one of the most effective methods in solving stiff nonlinear finite element problems. The ED-$\alpha$ method is based on a variational interpretation of Runge-Kutta methods and time-discontinuous Galerkin (TDG) methods~\cite{Bot97}. It has been shown that the ED-$\alpha$ scheme benefits from unconditional stability in the nonlinear regime and also damps out the unresolved and spurious high frequencies using tunable parameters~\cite{Bot04}.\par
Among different examples in the area of nonlinear structural dynamics, flexible multibody systems~\cite{Sha05,Bau10} are of interest for the present paper. Sudden impacts and persistent contacts in addition to the nonlinear behavior of the system due to large displacements and deformation increase the difficulty to solve the system of equations properly and increase the level of sensitivity regarding stability~\cite{Sch14a,Sch14b}. As an example, we consider a flexible slider-crank mechanism with clearance and dry friction. The numerical damping capabilities are explored by means of an elastic connecting rod, while the other parts of the system are assumed to be rigid bodies. The derivation and evaluation of the system matrices and vectors is oriented towards~\cite{Sha05} by applying the floating frame approach. For the elastic connecting rod, we consider a beam type element with three degrees of freedom at each node. Hence in Sect.~\ref{sec:equations}, we explain the equations of motion of nonsmooth mechanical systems and mixed timestepping schemes, which benefit from high-frequency damping. Their damping behavior is compared for a simple linear system with large stiffness and bilateral constraint. In Sect.~\ref{sec:application}, the slider-crank mechanism is introduced. Sect.~\ref{sec:examples} carries out the validation of the algorithm by comparing the results with the rigid body studies in~\cite{Flo10}. The convergence order is studied by comparing the results with a Simpack\footnote{http://www.simpack.com/} simulation. The damping properties of different time integration schemes are discussed in detail based on the amplitudes of the frequency spectrum and percentage of period elongations. In addition, a modal analysis is undertaken. Sect.~\ref{sec:summary} summarizes mixed timestepping schemes with high-frequency damping based on powerful and stable implicit time integration algorithms from computational mechanics such as the generalized-$\alpha$ method, the Bathe-method, and the ED-$\alpha$ method for complex nonlinear systems in structural dynamics. This paper is based on the student work~\cite{Rez13}. 

\section{Equations of motion and time-discretization\label{sec:equations}}
\begin{figure}[ht]
  \centering
  \includegraphics[width=0.85\columnwidth]{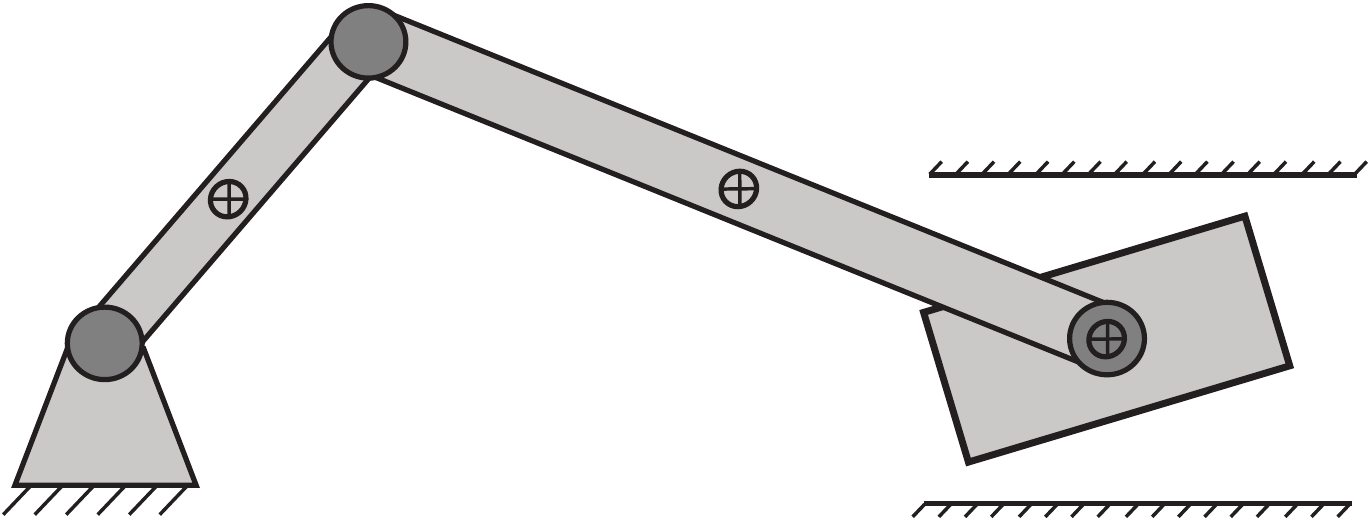}
  \caption{Slider crank mechanism.}    
  \label{fig:slider_crank}
\end{figure}

In this section, we first derive the general form of the dynamic equations of motion for flexible multibody systems subject to unilateral contacts and friction (Fig.~\ref{fig:slider_crank}). 
An example is the slider-crank mechanism for a predefined gap between slider and the bordering wall. We discuss how such mechanical system can be modeled and how the time evolution of such systems can be obtained by numerical integration. Therefore, we introduce and discuss the properties of the generalized-$\alpha$ method, the Bathe-method, and the ED-$\alpha$ method for solving the equations of motion implicitly in time. For each method, we solve for impulsive forces $\vLambda$ and non-impulsive forces $\vlambda$ based on the idea of separation of impact and contact forces~\cite{Sch12,Sch14a,Sch14b}.
\subsection{Nonsmooth approach}
\begin{figure}[ht]
  \centering
  \footnotesize
  \def\svgwidth{0.85\columnwidth}
  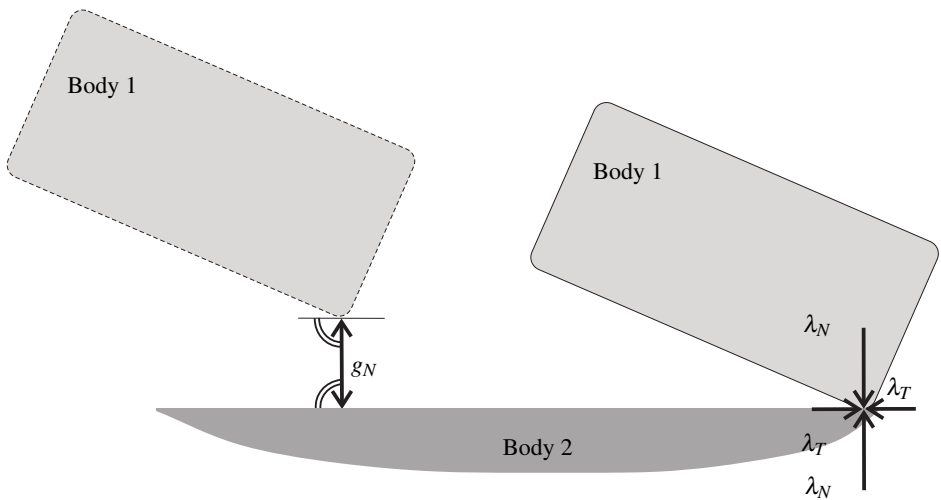
  \caption{Normal gap as well as normal and tangential contact force.}    
  \label{fig:gap_def}
\end{figure}

The nonsmooth approach is a method for modeling mechanical systems with unilateral contacts and friction using set-valued force laws. We consider the slider in Fig.~\ref{fig:slider_crank}, which slides, sticks or impacts on the cylinder wall. For each corner of the slider, normal $g_N$ and tangential $g_T$ gaps as well as normal $\lambda_N$ and tangential $\lambda_T$ reaction forces have to be evaluated according to Fig.~\ref{fig:gap_def}.

\begin{figure}[ht]
  \centering
  \footnotesize
  \def\svgwidth{0.95\columnwidth}
  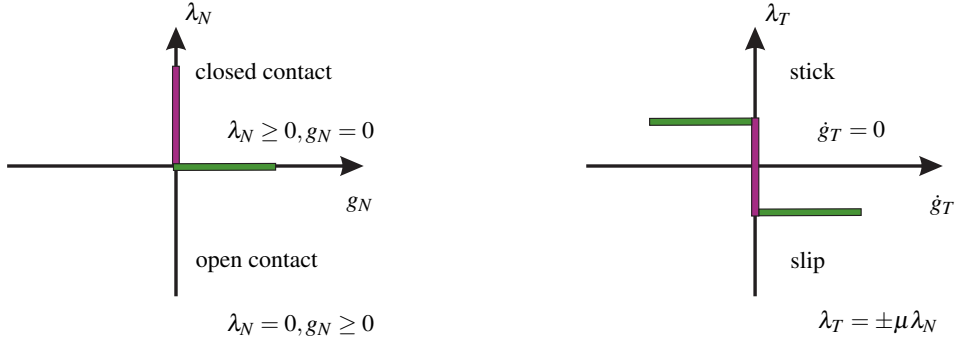
  \caption{Contact laws.}
  \label{fig:contact_law}
\end{figure}

The associated set-valued friction law of Coulomb-type is shown on the right-hand side of Fig.~\ref{fig:contact_law}. 
Regarding the sliding case, the friction force is given explicitly. Regarding the sticking case, the friction force is set-valued and determined according to an additional algebraic constraint. The nonsmooth approach changes the underlying mathematical structure if required and leads to a proper description of the mechanical systems (left-hand side of Fig.~\ref{fig:contact_law}). As a consequence of the changing mathematical structure, impacts can occur and the time evolution of the positions and the velocities cannot be assumed to be continuously differentiable anymore. Therefore, additional impact equations and impact laws have to be defined, which are discussed in the following sections.
\subsubsection{Timestepping integrators}
Timestepping schemes are based on a time-discretization of the system dynamics including the contact conditions in normal and tangential direction. The whole set of discretized equations and constraints is used to compute the next state of the motion. In contrast to event-driven schemes, these methods need no event-detection and are very robust in application. Moreover in the following section, a time discretization is used where the constraints are satisfied on velocity level. The accuracy of classic timestepping schemes is low. We propose higher order integration methods with high-frequency damping to process non-impulsive periods of the motion with larger time-step sizes and  increased integration order.
\subsubsection{Normal contact law}
The contact law for normal constraints is depicted on the left-hand side of Fig.~\ref{fig:contact_law}. The normal contact force $\vlambda_{N}$ vanishes if the bodies are separated ($\vg_{N} > 0$) and can only be positive if the bodies are in contact ($\vg_{N} = 0$). These constraints together form a complementarity condition:
\begin{align}
  0 \leq \vg_{N}\perp\vlambda_{N} \geq 0 \;. \label{eq:mechanical_system_normal_contact}
\end{align}
\subsubsection{Coulomb's friction law}
Coulomb's friction law is used to model friction forces. It states that the sliding friction force is proportional to the normal force of a contact. The amount of the static friction force is less than or equal to the maximum static friction force, which is also proportional to the normal force of a contact. The sliding friction force has the opposite direction of the relative velocity of the frictional contact. For closed contacts, the overall contact law in tangential direction is
\begin{align}
	\begin{aligned}
  &\norm{\vlambda_{T}} \leq \mu\vlambda_{N}\quad\text{for}\quad\dot{\vg}_{T}=\vnull\wedge\vg_{N}\leq\vnull \;,\\
  &\vlambda_{T}=-\frac{\dot{\vg}_{T}}{\norm{\dot{\vg}_{T}}}\mu\vlambda_{N}\quad\text{for }\quad\dot{\vg}_{T}\neq\vnull\wedge\vg_{N}\leq\vnull \;.
  \end{aligned} \label{eq:mechanical_system_tangential_contact}
\end{align}
The coefficient of friction $\mu$, in general, is a function of the relative velocity. We assume it to be a constant value.
\subsubsection{Projection formulation}
The projection formulation transforms the conditions shown in Fig.~\ref{fig:contact_law} in equivalent equations using convex analysis. The proximal point to $x$ in a convex set C is defined by the projection~\cite{Roc97}
\begin{align}
  \text{proj}_C\;:\;\MR\rightarrow\MR\;,\;x\mapsto\text{proj}_C(x)=\text{arg}\mathop{\text{min}}\limits_{x^* \in C} \left\|x-x^*\right\|\;.
\end{align}
With this definition, we express the normal contact law on velocity level as
\begin{align}
	\vlambda_{N_i}-\text{proj}_{\MR^+_0}\left(\vlambda_{N_i}-r\dot{\vg}_{N_i}\right)&=\vnull\label{eq:lambda_velocity}	
\end{align}
for all contacts belonging to the index set of closed constraints
\begin{align}
	\mathcal{I}_1&=\left\{k\in\mathcal{I}_0\;:\;\vg_{N_k}\leq 0\right\}\;,
\end{align}
where $\mathcal{I}_0$ contains all constraints. The arbitrary auxiliary parameter $r>0$ represents the slope of the regularizing function; it may be used for stabilizing the solution process~\cite{Sch11a}. In the same manner, we formulate the Coulomb friction law~\eqref{eq:mechanical_system_tangential_contact} as
\begin{align}
	\vlambda_{T_i}-\mathrm{proj}_{C_T(\vlambda_{N_i})}\left(\vlambda_{T_i}-r\dot{\vg}_{T_i}\right)&=\vnull\label{eq:lambda_tangential_velocity}	
\end{align}
for all contacts in $\mathcal{I}_1$, where the corresponding convex set is given by
\begin{align}
  C_T\;:\;\MR\rightarrow\mathcal{P}\left(\MR^2\right)\;:\;y\mapsto C_T(y)=\left\{ x \in\MR^2 \; | \;\; \|x\| \leq \mu|y| \right\}\;.
\end{align}
\subsubsection{Newton's law of impact}
For countable time instances $t_j$, the evolution of the slider-crank mechanism might get impulsive. For some component $k^*$ of the gap function $g_{N_{k^*}}\left(\vq\left(t_j\right)\right)=0$ but $g_{N_{k^*}}\left(\vq\left(t\right)\right)>0$ for $t_j-\delta \leq t<t_j$ and an appropriate $\delta>0$. This possibly leads to jumps in the velocity variables. Their derivatives do not exist anymore in the classical sense~\cite{Sch12,Sch14a,Sch14b}. We have to define the left-hand and right-hand limits:
\begin{align}
  \dot{\vg}_{N_j}^-:=\lim_{t\uparrow t_j}\dot{\vg}_N\left(t\right)&\;,\quad\dot{\vg}_{T_j}^-:=\lim_{t\uparrow t_j}\dot{\vg}_T\left(t\right)\;,\\
  \dot{\vg}_{N_j}^+:=\lim_{t\downarrow t_j}\dot{\vg}_N\left(t\right)&\;,\quad\dot{\vg}_{T_j}^+:=\lim_{t\downarrow t_j}\dot{\vg}_T\left(t\right)\;.
\end{align}
Then, the Lagrange multipliers describe the finite impulsive interaction in the sense of distributions:
\begin{align}
  \vLambda_{N_j}=\lim_{\delta\downarrow 0}\int_{t_j-\delta}^{t_j}\vlambda_N\text{d} t\;,\quad\vLambda_{T_j}=\lim_{\delta\downarrow 0}\int_{t_j-\delta}^{t_j}\vlambda_T\text{d} t \;.\label{eq:Lambda_int}
\end{align}
Newton's impact law describes the elasticity of the collision by considering the local velocities before ($ \dot{\vg}_{j}^-$) and after ($\dot{\vg}_{j}^+$) the impact:
\begin{align}
  &0 \leq \dot{\vg}_{N_j}^+ +\vepsilon_N \dot{\vg}_{N_j}^-\ \bot \ \vLambda_{N_j} \geq 0\;, \label{eq:mechanical_system_impact} \\
  &\norm{\vLambda_{T_j}} \leq \mu\vLambda_{N_j} \quad\text{for}\quad \dot{\vg}_{T_j}^++\vepsilon_T\dot{\vg}_{T_j}^-=\vnull \;, \\
  &\vLambda_{T_j}=-\frac{\dot{\vg}_{T_j}^+}{\norm{\dot{\vg}_{T_j}^+}}\mu\vLambda_{N_j} \quad\text{for}\quad \dot{\vg}_{T_j}^++\vepsilon_T\dot{\vg}_{T_j}^-\neq\vnull \;, \label{eq:mechanical_system_tangential_impact}
\end{align}
where $\vepsilon_N$ and $\vepsilon_T$ are the coefficients of restitution in normal and tangential direction, respectively. Therefore for the normal impact equations, we get
\begin{align}
  \vLambda_{N_j}-\mathrm{proj}_{\MR^+_0}\left(\vLambda_{N_j}-r(\dot{\vg}^+_{N_j}+\vepsilon_N\dot{\vg}^-_{N_j})\right)&=\vnull\;.\label{eq:Lambda_i_velocity}
\end{align}
For the tangential impact equations, we get
\begin{align}
	\vLambda_{T_j}-\mathrm{proj}_{C_T(\vLambda_{N_j})}\left(\vLambda_{T_j}-r(\dot{\vg}^+_{T_j}+\vepsilon_T\dot{\vg}^-_{T_j})\right)&=\vnull\;.\label{eq:Lambda_i_tangential_velocity}
\end{align}
\subsection{Equation of motion for a multibody system}
By determining the kinetic energy of the deformable bodies, the virtual work of the internal and external forces and the kinematic constraints that describe mechanical joints as well as specified bilateral or unilateral constraints, one can use Lagrange's equation to write the system equations of motion for multibody systems:
\begin{align}
  \vM\ddot{\vq} + \vC\dot{\vq} + \vK\vq  = \vh + \vW_{N}\vlambda_{N} + \vW_{T}\vlambda_{T} \;, \label{eq:mechanical_system_position}
\end{align}
where $\vM$ is the mass matrix, $\vC$ is the damping matrix, $\vK$ is the stiffness matrix, $\vW_N$ and $\vW_T$ are constraint matrices in normal and tangential direction, respectively, $\vlambda_N$ and $\vlambda_T$ are normal and tangential contact forces, respectively, and $\vh$ combines the effect of external forces and the quadratic velocity vector. If an impact occurs, the impact equations
\begin{align}
  &\vM_j \left[\vv_j^+-\vv_j^-\right] = \vW_{N_j}\vLambda_{N_j}+\vW_{T_j}\vLambda_{T_j} \;, \label{eq:mechanical_system_velocity_jump}
\end{align}
Newton's impact law~\eqref{eq:mechanical_system_impact} with restitution coefficient~$\vepsilon_N\in\left[0,1\right]$ and Newton's impact law~\eqref{eq:mechanical_system_tangential_impact} with $\vepsilon_T\in\left[0,1\right]$ have to be solved.

\subsection{Summary of computational algorithm}
\begin{figure}[ht]
  \centering
  \footnotesize
  \def\svgwidth{0.95\columnwidth}
  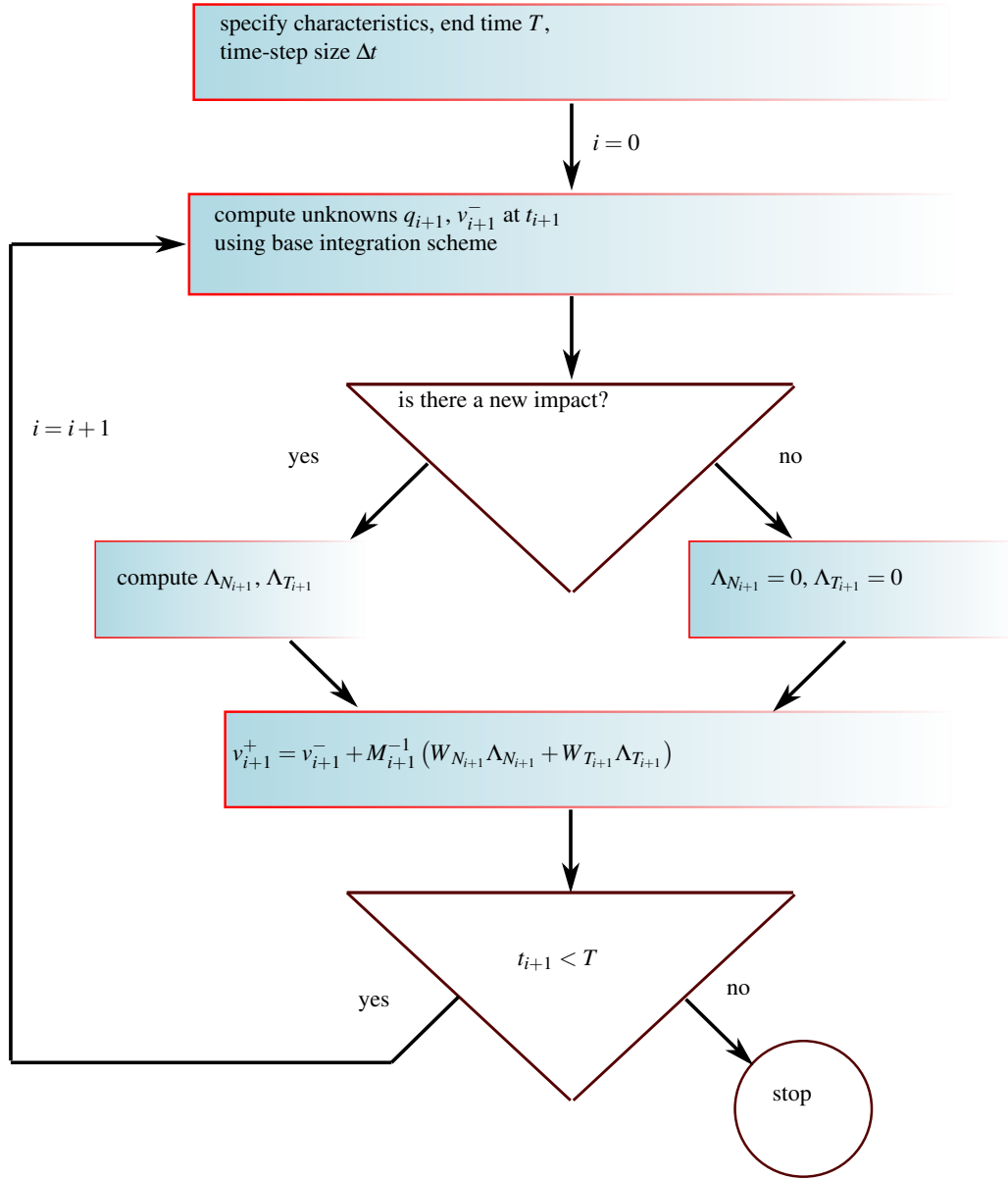
  \caption{Flowchart of computational algorithm~\cite{Sch14a,Sch14b}.}
  \label{fig:flow}
\end{figure}

The idea is to combine both the non-impulsive motion and the impulsive motion within one consistent integration scheme. First, this integration scheme has to model impacts and velocity jumps automatically if necessary. Second, it has to switch to effective higher order integration with all kinds of nice benefits of sophisticated integration schemes for differential algebraic equations. We introduce a framework, which is derived from a time-discontinuous Galerkin setting as proposed in~\cite{Sch14a,Sch14b}. It is more abstract and includes integration schemes for differential algebraic equations as base integration schemes on velocity level for one time-step from $t_i$ to $t_{i+1}$. With this propagation at the end of the time-step, we can check by the same activity rules, i.e., index set calculations on velocity level with~$\mathcal{I}_1$, if new contacts have been closed:
\begin{align}
  \exists k^*\,:\,g_{N_{k^*}}\left(\vq_{i}\right)>0\wedge g_{N_{k^*}}\left(\vq_{i+1}\right)\leq 0 \;.
\end{align}
In this case, we just correct the solution for the velocity variables and calculate $\vv_{i+1}^+$ with the impulsive forces $\vLambda_{N_{i+1}}$ and $\vLambda_{T_{i+1}}$. If
\begin{align}
  \nexists k^*\,:\,g_{N_{k^*}}\left(\vq_{i}\right)>0\wedge g_{N_{k^*}}\left(\vq_{i+1}\right)\leq 0\;,
\end{align}
we just set $\vv_{i+1}^+=\vv_{i+1}^-$. The overall algorithm can be summarized as shown in Fig.~\ref{fig:flow}.\par
\begin{figure}[ht]
  \centering
  \footnotesize
  \def\svgwidth{0.95\columnwidth}
  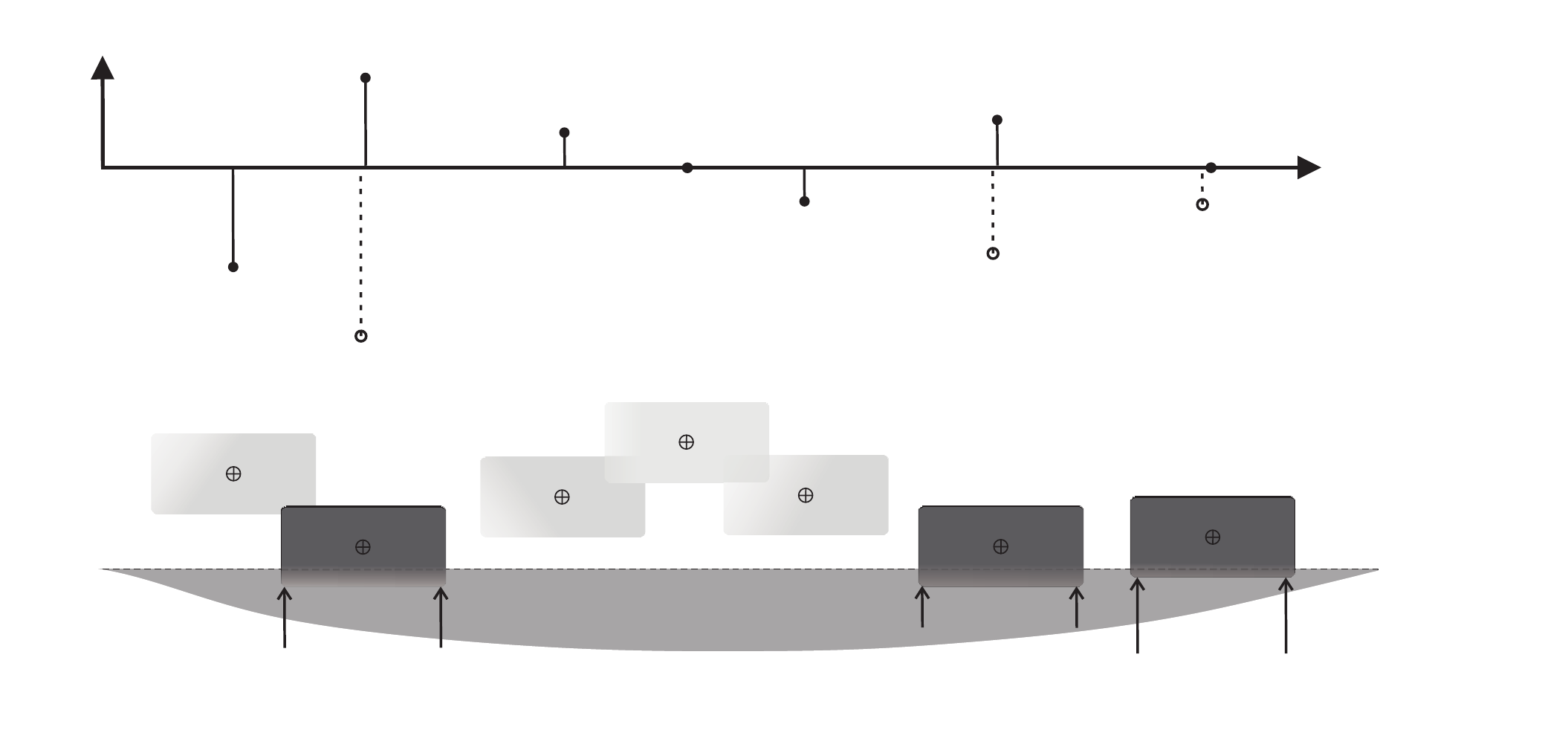
  \caption{Transition from impulsive to non-impulsive reactions.}
  \label{fig:impact_contact}
\end{figure}

As an example, Fig.~\ref{fig:impact_contact} shows the slider at different time-steps before and after an impact.
The dotted line represents the normal gap velocity before an impact. Going from $t_{i+5}$ to $t_{i+6}$, there is no \emph{new} active contact point. Hence, we just have to solve for non-impulsive forces on velocity level to avoid further penetration. We notice the usual drift-off effect (${g}_{N_{k}}\left(\vq_{i+6}\right) \leq 0$) after imposing the contact forces $\vlambda_N$ on velocity level which only constrains $\dot{\vg}_{N_{i+6}}\geq\vnull$.\par
We proceed with explaining the generalized-$\alpha$ method as a base integration scheme. Then, we focus on how to calculate impulsive corrections. The Bathe-method as well as the ED-$\alpha$ method are further base integration schemes, which can be used instead of the generalized-$\alpha$ method. We explain them and compare the properties of all three base integration schemes at the end of the following section. 
\subsection{Generalized-$\valpha$ method} 
According to~\cite{Arn07}, the generalized-$\alpha$ method for flexible multibody systems can be summarized as follows:
\begin{align}
  &\vM_{i+1} \va_{i+1} + \vC_{i+1} \vv^{-}_{i+1} + \vK_{i+1} \vq_{i+1} = \nonumber\\ 
  &\qquad\vh^{-}_{i+1} + \vW_{N_{i+1}}\vlambda_{N_{i+1}} + \vW_{T_{i+1}}\vlambda_{T_{i+1}}\;,\label{eq:mechanical_system} \\
  &\left(1-\alpha_m\right)\vA_{i+1}+\alpha_m\vA_i = \left(1-\alpha_f\right)\va_{i+1}+\alpha_f\va_i \;, \label{eq:general_alpha_recurrence}\\
  &\vq_{i+1} = \vq_i+\Delta t\vv^{+}_i+\Delta t^2\left[\left(0.5-\beta\right)\vA_i+\beta\vA_{i+1}\right] \;,\label{eq:general_alpha_position}\\
  &\vv^{-}_{i+1} = \vv^{+}_i+\Delta t\left[\left(1-\gamma\right)\vA_i+\gamma\vA_{i+1}\right] \;, \label{eq:general_alpha_velocity}\\
  &\vq_0 = \vq\left(0\right) \;, \quad \vv^{+}_0=\vv\left(0\right) \;, \label{eq:mechanical_system_initial_condition}\\
  &\va_0 = \vM^{-1}_0\left(\vh_0 - \vK_0\vq_0 - \vC_0\vv^{+}_0\right) \;, \quad \vA_0 = \va_0 \;,\label{eq:mechanical_system_initial_acceleration}
\end{align}\label{eq:general_alpha}
where $\vq_i$ is the vector of generalized coordinates, $\vv_i$ is the vector of generalized velocities, $\va_i$ is the vector of generalized accelerations, $\vA_i$ is the vector of acceleration-like auxiliary variables, defined by the recurrence relation~\eqref{eq:general_alpha_recurrence}, and $\Delta t$ is the time-step size.
\subsubsection{General characteristics}
\begin{figure}[ht]
  \centering
  \footnotesize
  \def\svgwidth{0.95\columnwidth}
  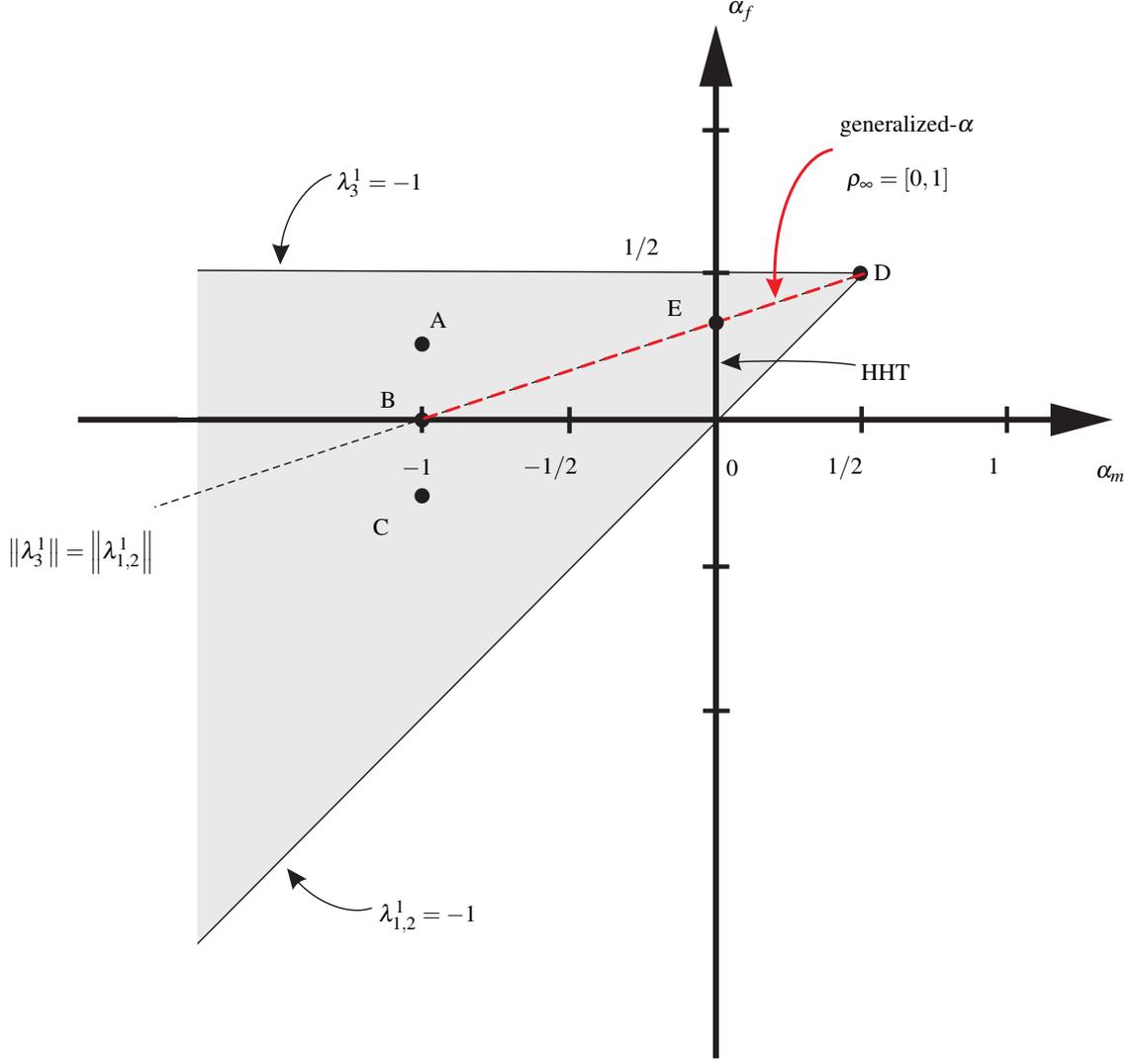
  \caption{Generalized-$\alpha$ stability region in $\alpha_m$-$\alpha_f$ space, different test case $A(-1,1/2)$, $B(-1,0)$, $C(-1,-1/2)$, $D(1/2,1/2)$, $E(0,1/3)$~\cite{Chu93}.}
  \label{fig:alpha_space}
\end{figure}

\begin{figure}[ht]
  \centering
  \includegraphics[width=\columnwidth]{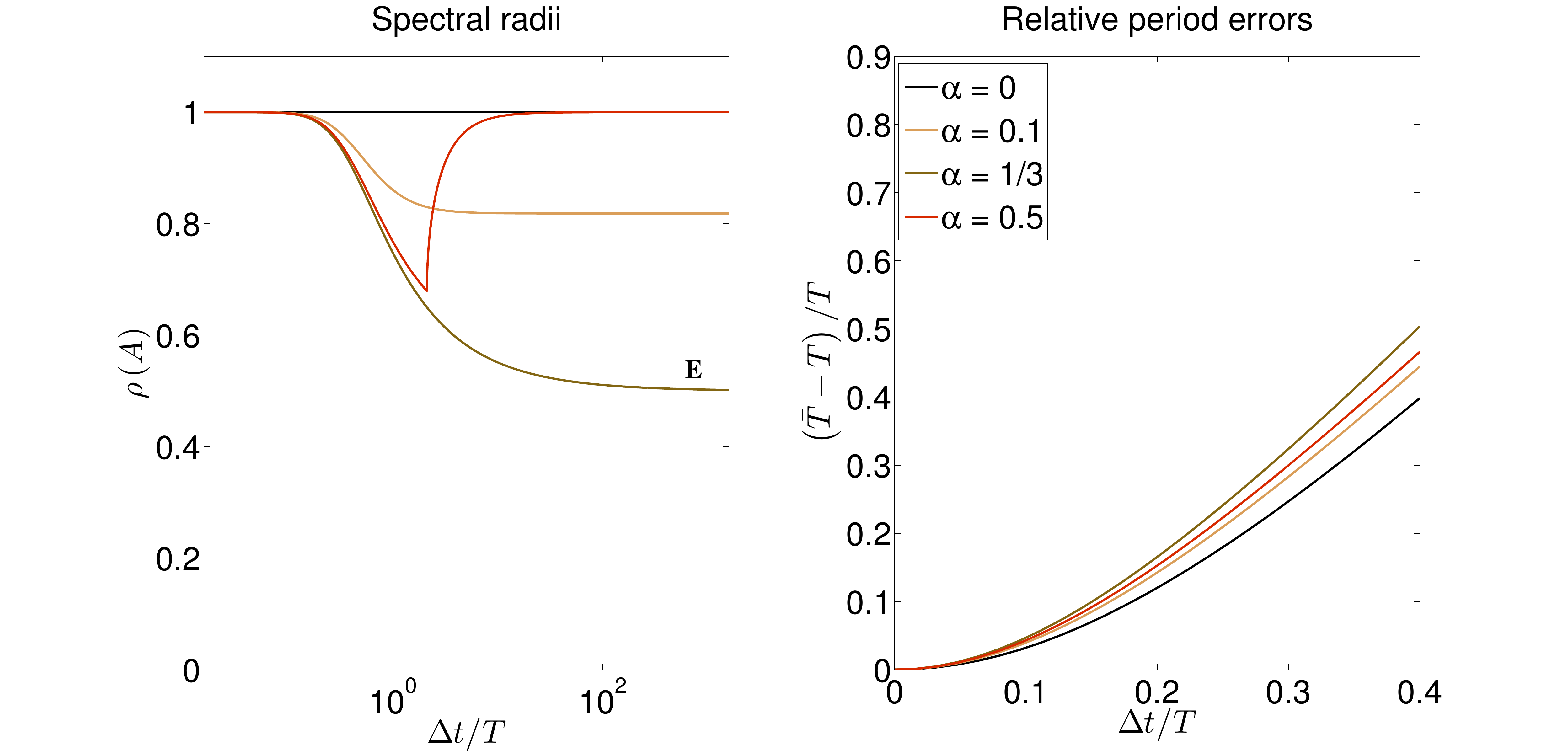}
  \caption{Properties of the HHT-method for different $\alpha=\alpha_f$, $\alpha_m=0$.}
  \label{fig:hht_sum}
\end{figure}

\begin{figure}[ht]
  \centering
  \includegraphics[width=\columnwidth]{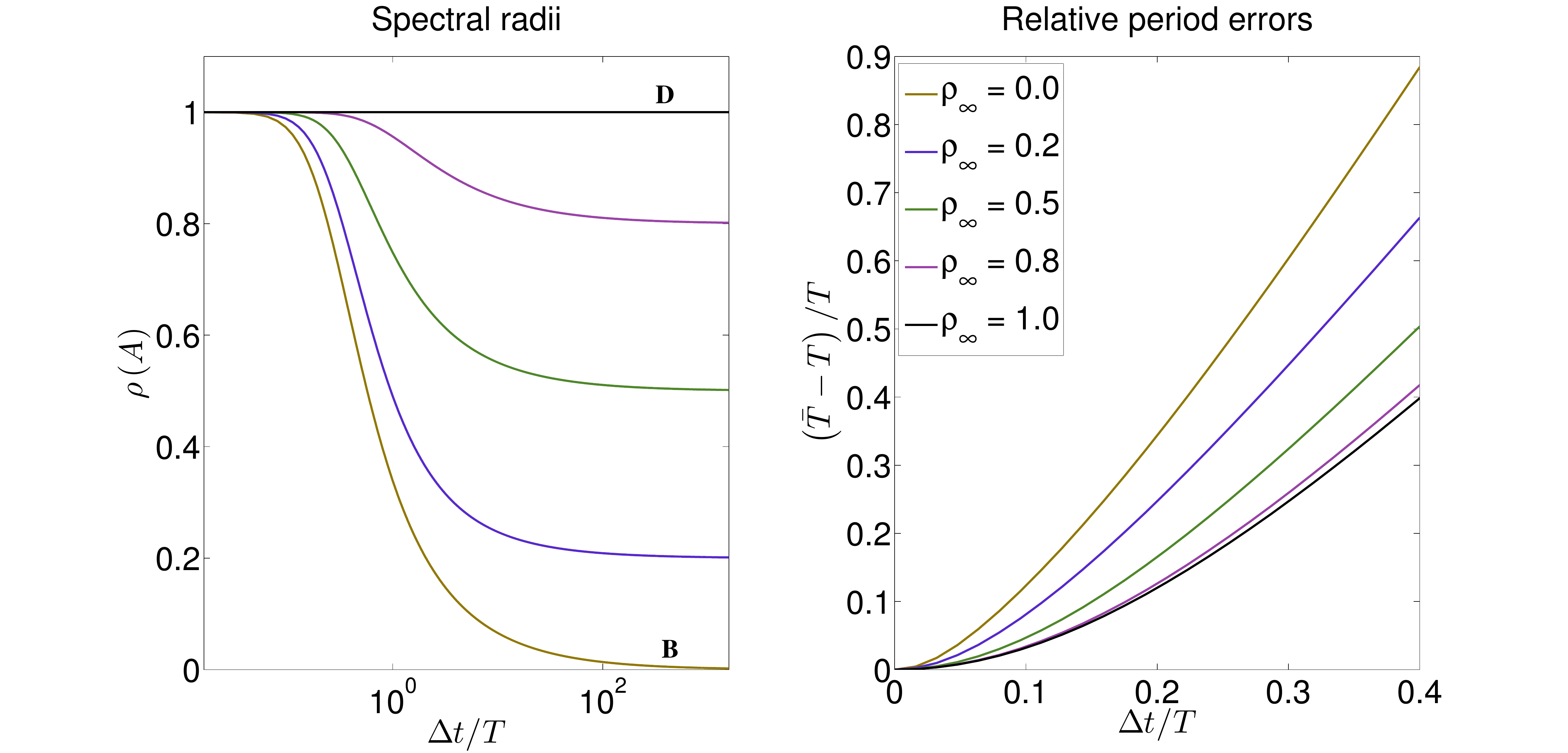}
  \caption{Properties of the generalized-$\alpha$ method for different $\rho_{\infty}$.}
  \label{fig:ga_sum}
\end{figure}

\begin{figure}[ht]
  \centering
  \includegraphics[width=\columnwidth]{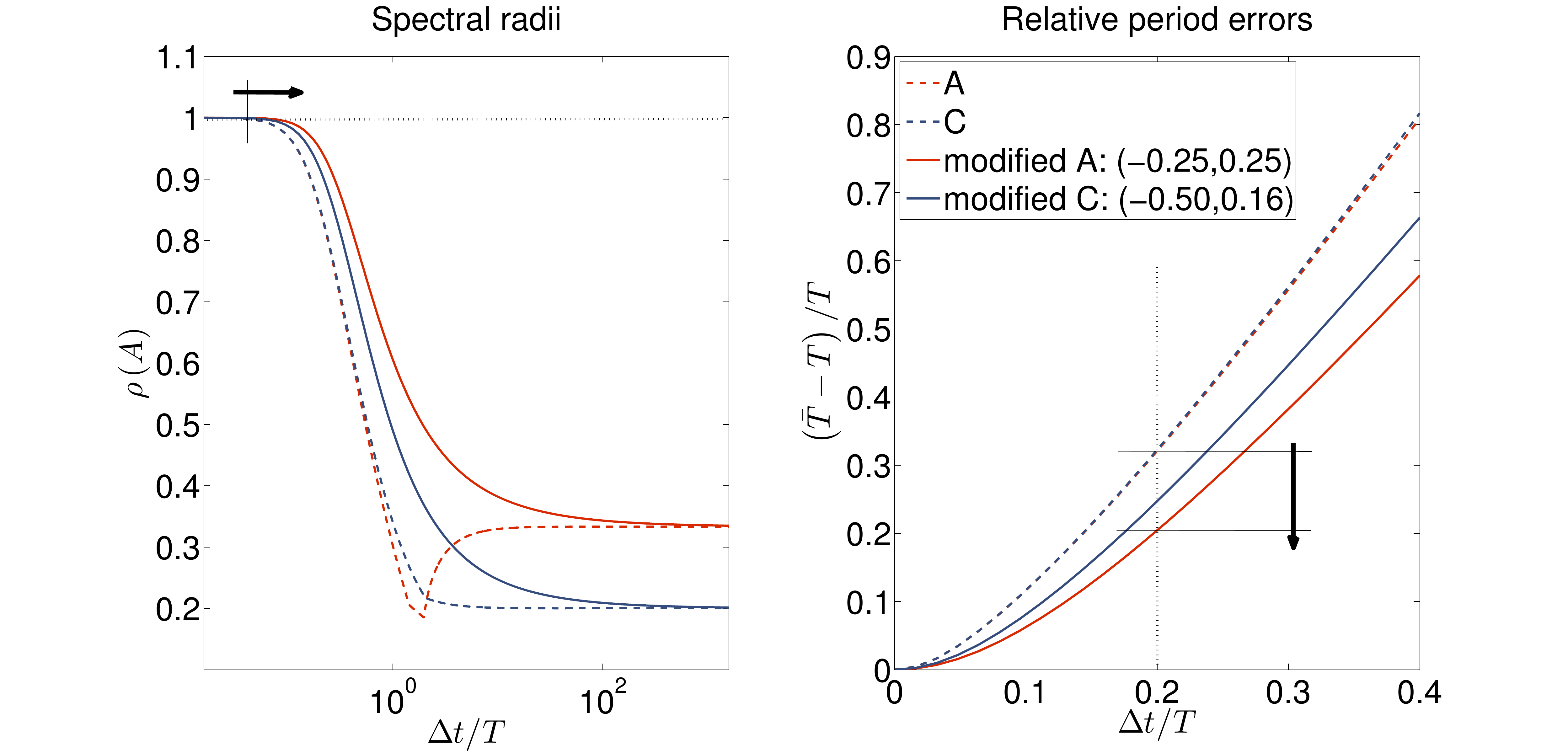}
  \caption{Smooth transition of cases $A$ and $C$ by modifying $\alpha_f$ and $\alpha_m$.}
  \label{fig:ABC_modified}
\end{figure}

The displacement and velocity update~\eqref{eq:general_alpha_position} and \eqref{eq:general_alpha_velocity} are identical to those of the Newmark algorithm. The structure of these update equations is obtained using Taylor series expansion about $t_i$. The crucial task is to determine the relationship between the algorithmic parameters, $\alpha_m$, $\alpha_f$, $\gamma$, and $\beta$. With appropriate expressions for $\gamma$ and $\beta$, and if $\alpha_m=0$, the algorithm reduces to the HHT-$\alpha$ method. The generalized-$\alpha$ method, with parametric values given in \eqref{eq:general_alpha_par}, is unconditionally stable for linear problems, second order accurate possessing an optimal combination of high-frequency and low-frequency dissipation:
\begin{align}
  \begin{aligned}
    \alpha_m &=\frac{2\rho_{\infty}-1}{\rho_{\infty}+1} \;, \\
    \alpha_f &=\frac{\rho_{\infty}}{\rho_{\infty}+1} \;, \\
    \gamma &= \frac{1}{2}-\alpha_m+\alpha_f \;, \\
    \beta &= \frac{1}{4}\left(1-\alpha_m+\alpha_f\right)^2 \;,
\end{aligned} \label{eq:general_alpha_par}
\end{align} 
where $\rho_{\infty}\in\left[0,1\right]$ is the spectral radius of the amplification matrix at the high frequency limit. The stability region is indicated by the shaded area in Fig.~\ref{fig:alpha_space} (\ref{sec:appA}).
In order to have a better insight on the effect of the different parameters on the eigenvalues of the amplification matrix, we plot the spectral radius and the relative period error with respect to $\Delta t/T$, i.e., the time-step divided by the period $T=2\pi/\omega$. First, we consider a special case of the generalized-$\alpha$ method by setting $\alpha_m=0$ (HHT method). From Fig.~\ref{fig:alpha_space}, we are allowed to choose $\alpha_f$ between $0$ and $0.5$ (Fig.~\ref{fig:hht_sum}). We see a \emph{cusp} after we pass the point $E$.
Using optimum values for the generalized-$\alpha$ method as introduced in \eqref{eq:general_alpha_par}, we see the behavior of the method for different $\rho_{\infty}$ in Fig.~\ref{fig:ga_sum} from the no-dissipation case ($\rho_{\infty}$=1) to the so-called asymptotic annihilation case ($\rho_{\infty}$=0), or moving along the red dotted line from point $D$ to point $B$.
In order to observe the properties of different regions in Fig.~\ref{fig:alpha_space}, we plot the spectral radius for point $A$ and point $C$ (Fig.~\ref{fig:ABC_modified}). If we select the $\alpha_m$ and $\alpha_f$ values away from the dotted line ($\norm{\lambda^{\infty}_{3}}=\norm{\lambda^{\infty}_{1,2}}$), we may expect the cusp. Figure~\ref{fig:ABC_modified} shows how we can modify the $\alpha_m$ and $\alpha_f$ values for point $A$ and point $C$ in order to have the same $\rho_{\infty}$ but with smooth transition instead.
Further properties can be found in~\ref{sec:appA}.
\subsubsection{Calculation of contact forces on velocity level}\label{sec:lambda}
Using \eqref{eq:mechanical_system_normal_contact} and \eqref{eq:mechanical_system_tangential_contact} together with \eqref{eq:mechanical_system_position}, we calculate contact forces $\vlambda_{N_{i+1}}$ and $\vlambda_{T_{i+1}}$. We interpret \eqref{eq:mechanical_system_normal_contact} and \eqref{eq:mechanical_system_tangential_contact} on velocity level~\cite{Sch14a,Sch14b}:
\begin{align}
	\vlambda_{N_{i+1}}&=\begin{cases}
		0&\mathrm{if~} \vg_{N_i}>0\\
\mathrm{proj}_{\MR^+_0}\left[\vlambda_{N_{i+1}}-r~\dot{\vg}^{-}_{N_{i+1}}\right]&\mathrm{else}
	\end{cases}\;,\label{eq:lambda_component_wise_normal}\\
	\vlambda_{T_{i+1}}&=\begin{cases}
		0&\mathrm{if~} \vg_{N_i}>0\\
\mathrm{proj}_{C_T(\vlambda_{N_{i+1}})}\left[\vlambda_{T_{i+1}}-r~\dot{\vg}^{-}_{T_{i+1}}\right]&\mathrm{else}
	\end{cases}\;,\label{eq:lambda_component_wise_tangential}
\end{align}
using the projection-formulation row-by-row and a predictor for closed contacts $\vg_{N_i}\leq 0$. Thereby, the local velocities satisfy
\begin{align}
  \dot{\vg}^{-}_{N_{i+1}}=\vW^T_{N_{i+1}}\vv^{-}_{i+1}\;,\quad \dot{\vg}^{-}_{T_{i+1}}=\vW^T_{T_{i+1}}\vv^{-}_{i+1}\;,\label{eq:implicit_velocity_gap}
\end{align}
and depend on $\vlambda_{N_{i+1}}$ and $\vlambda_{T_{i+1}}$ because of the \emph{implicit representation}. Hence, we calculate equivalent forces $\vlambda_{N_{i+1}}$ and $\vlambda_{T_{i+1}}$ such that the given local velocities are projected into their respective admissible space. Using \eqref{eq:general_alpha_velocity}, \eqref{eq:general_alpha_recurrence} and \eqref{eq:mechanical_system}, we substitute:
\begin{align}
  &\dot{\vg}^{-}_{N_{i+1}}=\vF_{N_{i+1}}+\Delta t\gamma\frac{1-\alpha_f}{1-\alpha_m} \left\{\vG_{N_{i+1}}\vlambda_{N_{i+1}}+\vG_{NT_{i+1}}\vlambda_{T_{i+1}}\right\}\;,\label{eq:velocity_gap_normal}\\
  &\dot{\vg}^{-}_{T_{i+1}}=\vF_{T_{i+1}}+\Delta t\gamma\frac{1-\alpha_f}{1-\alpha_m} \left\{\vG_{TN_{i+1}}\vlambda_{N_{i+1}}+\vG_{T_{i+1}}\vlambda_{T_{i+1}}\right\}\;,\label{eq:velocity_gap_tangential}
\end{align}
where
\begin{align}
  &\vF_{N_{i+1}}=\vW_{N_{i+1}}^T\vv^{+}_i+\Delta t\left(1-\gamma\right)\vW_{N_{i+1}}^T\vA_i+\Delta t\gamma\vW_{N_{i+1}}^T\left(\frac{1-\alpha_f}{1-\alpha_m}\widehat{\vM}^{-1}_{i+1}\widehat{\vR}_{i+1}+\frac{\alpha_f}{1-\alpha_m}\va_i-\frac{\alpha_m}{1-\alpha_m}\vA_i\right)\;,\\
  &\vF_{T_{i+1}}=\vW_{T_{i+1}}^T\vv^{+}_i+\Delta t\left(1-\gamma\right)\vW_{T_{i+1}}^T\vA_i+\Delta t\gamma\vW_{T_{i+1}}^T\left(\frac{1-\alpha_f}{1-\alpha_m}\widehat{\vM}^{-1}_{i+1}\widehat{\vR}_{i+1}+\frac{\alpha_f}{1-\alpha_m}\va_i-\frac{\alpha_m}{1-\alpha_m}\vA_i\right)\;, 
\end{align}
and
\begin{align}
  \vG_{N_{i+1}}&=\vW^T_{N_{i+1}}\widehat{\vM}^{-1}_{i+1}\vW_{N_{i+1}}\;, \quad \vG_{NT_{i+1}}=\vW^T_{N_{i+1}}\widehat{\vM}^{-1}_{i+1}\vW_{T_{i+1}}\;,\\
  \vG_{T_{i+1}}&=\vW^T_{T_{i+1}}\widehat{\vM}^{-1}_{i+1}\vW_{T_{i+1}}\;, \quad \vG_{TN_{i+1}}=\vW^T_{T_{i+1}}\widehat{\vM}^{-1}_{i+1}\vW_{N_{i+1}}
\end{align}
are called Delassus matrices with the effective mass matrix
\begin{align}
  \widehat{\vM}_{i+1}=\vM_{i+1}+\Delta t_i\gamma\frac{1-\alpha_f}{1-\alpha_m}\vC_{i+1}+\Delta t^2_i\beta\frac{1-\alpha_f}{1-\alpha_m}\vK_{i+1}
\end{align}
and the effective right-hand side
\begin{align}
  \widehat{\vR}_{i+1}&=\vh_{i+1}^--\vC_{i+1}\vv_i^+-\left(\Delta t_i\left(1-\gamma\right)-\Delta t_i\gamma\frac{\alpha_m}{1-\alpha_m}\right)\vC_{i+1}\vA_i-\Delta t_i\gamma\frac{\alpha_f}{1-\alpha_m}\vC_{i+1}\va_i-\vK_{i+1}\vq_i-\Delta t_i\vK_{i+1}\vv_i^+\\
	&\quad-\left(\Delta t^2_i\left(0.5-\beta\right)-\Delta t^2_i\beta\frac{\alpha_m}{1-\alpha_m}\right)\vK_{i+1}\vA_i-\Delta t^2_i\beta\frac{\alpha_f}{1-\alpha_m}\vK_{i+1}\va_i \;.
\end{align}
We focus on active contacts and transform \eqref{eq:lambda_component_wise_normal} and \eqref{eq:lambda_component_wise_tangential} formally using row-by-row interpretation:
\begin{align}
	\vlambda_{N_{i+1},\mathcal{I}_1^i}&=\mathrm{proj}_{\MR^+_0}\left[\vlambda_{N_{i+1},\mathcal{I}_1^i}-r\dot{\vg}^{-}_{N_{i+1},\mathcal{I}_1^i}\right]\;,\label{eq:active_contacts_normal_formally}\\
	\vlambda_{T_{i+1},\mathcal{I}_1^i}&=\mathrm{proj}_{C_T(\vlambda_{N_{i+1},\mathcal{I}_1^i})}\left[\vlambda_{T_{i+1},\mathcal{I}_1^i}-r\dot{\vg}^{-}_{T_{i+1},\mathcal{I}_1^i}\right]\;.\label{eq:active_contacts_tangential_formally}
\end{align}
Thereby, the index set of closed constraints at $t_i$ is given by
\begin{align}
	\mathcal{I}_1^i&=\left\{k\in\mathcal{I}_0\;:\;\vg_{N_k}(\vq_i)\leq 0\right\}\;.
\end{align}
In the multi-contact case, active contacts might be depending. Hence, we cannot use a nonsmooth Newton method but we solve \eqref{eq:active_contacts_normal_formally}, \eqref{eq:active_contacts_tangential_formally} by a nonsmooth variant of the Gauss-Newton method, i.e., we choose an approximate root $\left(\bar{\vlambda}_{N,\mathcal{I}_1},\bar{\vlambda}_{T,\mathcal{I}_1}\right)$ of the function
\begin{align}
	&f\;:\;\MR^{|\mathcal{I}_1|}\times\MR^{2|\mathcal{I}_1|}\rightarrow\MR^{|\mathcal{I}_1|}\times\MR^{2|\mathcal{I}_1|}\;,\nonumber\\
 &(\vlambda_{N,\mathcal{I}_1},\vlambda_{T,\mathcal{I}_1})\mapsto f(\vlambda_{N,\mathcal{I}_1},\vlambda_{T,\mathcal{I}_1})\nonumber\\
 &=\begin{pmatrix} \vlambda_{N,\mathcal{I}_1}-\mathrm{proj}_{\MR^+_0}\left[\vlambda_{N,\mathcal{I}_1}-r\dot{\vg}^{-}_{N,\mathcal{I}_1}(\vlambda_{N,\mathcal{I}_1},\vlambda_{T,\mathcal{I}_1})\right]\\
  \vlambda_{T,\mathcal{I}_1}-\mathrm{proj}_{C_T(\vlambda_{N,\mathcal{I}_1})}\left[\vlambda_{T,\mathcal{I}_1}-r\dot{\vg}^{-}_{T,\mathcal{I}_1}(\vlambda_{N,\mathcal{I}_1},\vlambda_{T,\mathcal{I}_1})\right]\end{pmatrix}\;,
	\label{eq:contact_f}
\end{align} 
with the Moore-Penrose pseudoinverse operator pinv. Using the Gauss-Newton algorithm, we calculate unknown contact forces with the following iterative algorithm:
\begin{align}
  \begin{aligned}
    &\left(\bar{\vlambda}_{N,\mathcal{I}_1},\bar{\vlambda}_{T,\mathcal{I}_1}\right)=\left(\vnull,\vnull\right)\;,
    ~\bar{f}=f\left(\bar{\vlambda}_{N,\mathcal{I}_1},\bar{\vlambda}_{T,\mathcal{I}_1}\right)\\
		&\mathrm{while~}\left\|\bar{f}\right\|>\mathrm{tol}\\
  &\hspace{0.1cm}\nabla f\left(\bar{\vlambda}_{N,\mathcal{I}_1},\bar{\vlambda}_{T,\mathcal{I}_1}\right)=\begin{pmatrix}\vI & \vnull\\\vnull & \vI\end{pmatrix}-\begin{pmatrix}\vTheta_N(\vI-r\Delta t\gamma\frac{1-\alpha_f}{1-\alpha_m}\vG_{N,\mathcal{I}_1}) & \vTheta_N(-r\Delta t\gamma\frac{1-\alpha_f}{1-\alpha_m}\vG_{NT,\mathcal{I}_1})\\\vTheta_T(-r\Delta t\gamma\frac{1-\alpha_f}{1-\alpha_m}\vG_{TN,\mathcal{I}_1}) & \vTheta_T(\vI-r\Delta t\gamma\frac{1-\alpha_f}{1-\alpha_m}\vG_{T,\mathcal{I}_1})\end{pmatrix}\\
  &\hspace{0.1cm}\left(\bar{\vlambda}_{N_{\text{new}},\mathcal{I}_1},\bar{\vlambda}_{T_{\text{new}},\mathcal{I}_1}\right)=\left(\bar{\vlambda}_{N,\mathcal{I}_1},\bar{\vlambda}_{T,\mathcal{I}_1}\right)-\mathrm{pinv}(\nabla f\left(\bar{\vlambda}_{N,\mathcal{I}_1},\bar{\vlambda}_{T,\mathcal{I}_1}\right))\bar{f}\\
  &\hspace{0.1cm}\left(\bar{\vlambda}_{N,\mathcal{I}_1},\bar{\vlambda}_{T,\mathcal{I}_1}\right)=\left(\bar{\vlambda}_{N_{\text{new}},\mathcal{I}_1},\bar{\vlambda}_{T_{\text{new}},\mathcal{I}_1}\right)\\
  &\hspace{0.1cm}\bar{f}=f\left(\bar{\vlambda}_{N,\mathcal{I}_1},\bar{\vlambda}_{T,\mathcal{I}_1}\right)\\
 		&\mathrm{end}
  \end{aligned}\label{eq:Newton}
\end{align}
The occurring Heaviside functions
\begin{align}
  &\vTheta_N\;:\;\MR^{|\mathcal{I}_1|}\times\MR^{2|\mathcal{I}_1|}\rightarrow\text{diag}^{{|\mathcal{I}_1|},{|\mathcal{I}_1|}}\;,\quad\vTheta_{N_{kk}}\left(\vlambda_{N,\mathcal{I}_1},\vlambda_{T,\mathcal{I}_1}\right)=\begin{cases}
		0&\mathrm{if~}{\lambda}_{N_k}-r\dot{\vg}^{-}_{N_{k}}<0\\
 		1&\mathrm{else}
	\end{cases}\;,\\
  &\vTheta_T\;:\;\MR^{|\mathcal{I}_1|}\times\MR^{2|\mathcal{I}_1|}\rightarrow\text{diag}^{{|\mathcal{I}_1|},{|\mathcal{I}_1|}}\;,\quad\vTheta_{T_{kk}}\left(\vlambda_{N,\mathcal{I}_1},\vlambda_{T,\mathcal{I}_1}\right)=\begin{cases}
    0&\mathrm{if~}\norm{{\lambda}_{T_k}-r\dot{\vg}^{-}_{T_{k}}}>\mu\abs{{\lambda}_{N_k}}\\
    1&\mathrm{else}
  \end{cases}\;,
\end{align}
are interpreted row-by-row and we use $r=0.1$ without adaptation to improve convergence of the numerical scheme~\cite{Sch11a}. 
\subsection{Assignment of impulses} \label{sec:Lambda}
If on the other hand, there is at least one inactive gap function ($g_{N_{k^*}}\left(\vq_{i}\right)>0$) at $t_{i}$ that becomes active ($g_{N_{k^*}}\left(\vq_{i+1}\right)\leq0$), an impact occurs in a rigidly connected component of its source. Hence, we update the velocity at $t_{i+1}$ using \eqref{eq:mechanical_system_velocity_jump}:
\begin{align}
  \vv_{i+1}^+= \vv_{i+1}^- + \vM^{-1}_{i+1}\vW_{N_{i+1}}\vLambda_{N_{i+1}}+\vM^{-1}_{i+1}\vW_{T_{i+1}}\vLambda_{T_{i+1}}\;,\label{eq:vi-nonsmooth}
\end{align}
where $\left(~\right)^-$ and $\left(~\right)^+$ denote the parameter at $t_{i+1}$ before and after the impact, respectively. For the computation of $\vv_{i+1}^+$, $\vLambda_{N_i}$ and $\vLambda_{T_i}$, we write \eqref{eq:mechanical_system_impact} and \eqref{eq:mechanical_system_tangential_impact} on velocity level row-by-row:
\begin{align}
	&\vLambda_{N_{i+1}}=\begin{cases}
    0&\mathrm{if~}\vg_{N_{i+1}}>0\\
\mathrm{proj}_{\MR^+_0}\left[\vLambda_{N_{i+1}}-r\underbrace{(\dot{\vg}_{N_{i+1}}^++\vepsilon_N\dot{\vg}_{N_{i+1}}^-)}_{\dot{\bar{\vg}}_{N_{i+1}}}\right]&\mathrm{else~}
  \end{cases}\label{eq:Lambda_component_wise_normal}\;,\\
	&\vLambda_{T_{i+1}}=\begin{cases}
    0&\mathrm{if~}\vg_{N_{i+1}}>0\\
\mathrm{proj}_{C_T(\vLambda_{N_{i+1}})}\left[\vLambda_{T_{i+1}}-r\underbrace{(\dot{\vg}_{T_{i+1}}^++\vepsilon_T\dot{\vg}_{T_{i+1}}^-)}_{\dot{\bar{\vg}}_{T_{i+1}}}\right]&\mathrm{else~}
  \end{cases}\;.\label{eq:Lambda_component_wise_tangential}
\end{align}
Equation~\eqref{eq:mechanical_system_velocity_jump} can be used to eliminate $\dot{\vg}_{N_{i+1}}^{+}$ and $\dot{\vg}_{T_{i+1}}^+$, row-by-row resulting in
\begin{align}
  \vLambda_{N_{i+1},\mathcal{I}_1^{i+1}}&=\mathrm{proj}_{\MR^+_0}\left[\vLambda_{N_{i+1},\mathcal{I}_1^{i+1}}-r\dot{\bar{\vg}}_{N_{i+1},\mathcal{I}_1^{i+1}}\right]\;,\\
  \vLambda_{T_{i+1},\mathcal{I}_1^{i+1}}&=\mathrm{proj}_{C_T(\vLambda_{N_{i+1},\mathcal{I}_1^{i+1}})}\left[\vLambda_{T_{i+1},\mathcal{I}_1^{i+1}}-r\dot{\bar{\vg}}_{T_{i+1},\mathcal{I}_1^{i+1}}\right]\;.
\end{align}
Again we look for the roots of
\begin{align}
	\begin{aligned} 		
    &f\;:\;\MR^{|\mathcal{I}_1|}\times\MR^{2|\mathcal{I}_1|}\rightarrow\MR^{|\mathcal{I}_1|}\times\MR^{2|\mathcal{I}_1|}\;,\\
    &\left(\vLambda_{N,\mathcal{I}_1},\vLambda_{T,\mathcal{I}_1}\right)\mapsto f(\vLambda_{N,\mathcal{I}_1},\vLambda_{T,\mathcal{I}_1})= \begin{pmatrix}\vLambda_{N,\mathcal{I}_1}-\mathrm{proj}_{\MR^+_0}\left[\vLambda_{N,\mathcal{I}_1}-r\dot{\bar{\vg}}_{N,\mathcal{I}_1}(\vLambda_{N,\mathcal{I}_1},\vLambda_{T,\mathcal{I}_1})\right]\\
    \vLambda_{T,\mathcal{I}_1}-\mathrm{proj}_{C_T(\vLambda_{N,\mathcal{I}_1})}\left[\vLambda_{T,\mathcal{I}_1}-r\dot{\bar{\vg}}_{T,\mathcal{I}_1}(\vLambda_{N,\mathcal{I}_1},\vLambda_{T,\mathcal{I}_1})\right]
     \end{pmatrix}\;,
	\end{aligned} \label{eq:impact_f}
\end{align} 
with a nonsmooth Gauss-Newton method as in \eqref{eq:Newton}. The derivative of $f$ at $\left(\vLambda_{N,\mathcal{I}_1},\vLambda_{T,\mathcal{I}_1}\right)$ is given by
\begin{align}
  \begin{aligned}
    &\nabla f\left(\vLambda_{N,\mathcal{I}_1},\vLambda_{T,\mathcal{I}_1}\right)=\begin{pmatrix}\vI-\vTheta_{N}(\vI-r\vG_{N,\mathcal{I}_1}) & \vTheta_{N}(-r\vG_{NT,\mathcal{I}_1})\\ 
    \vTheta_{T}(-r\vG_{TN,\mathcal{I}_1}) & \vI-\vTheta_{T}(\vI-r\vG_{T,\mathcal{I}_1})\end{pmatrix}\;.
  \end{aligned}
\end{align}
The occurring Heaviside functions
\begin{align}
  &\vTheta_{N}\;:\;\MR^{|\mathcal{I}_1|}\times\MR^{2|\mathcal{I}_1|}\rightarrow\text{diag}^{|\mathcal{I}_1|,|\mathcal{I}_1|}\;,\vTheta_{N_{kk}}\left(\vLambda_{N,\mathcal{I}_1},\vLambda_{T,\mathcal{I}_1}\right)=\begin{cases}
		0&\mathrm{if~} \Lambda_{N_k}-r\dot{\bar{\vg}}_{N_{k}}<0\\
 		1&\mathrm{else}
	\end{cases}\;,\\
  &\vTheta_{T}\;:\;\MR^{|\mathcal{I}_1|}\times\MR^{2|\mathcal{I}_1|}\rightarrow\text{diag}^{|\mathcal{I}_1|,|\mathcal{I}_1|}\;\vTheta_{T_{kk}}\left(\vLambda_{N,\mathcal{I}_1},\vLambda_{T,\mathcal{I}_1}\right)=\begin{cases}
    0&\mathrm{if~}\norm{\Lambda_{T_k}-r\dot{\bar{\vg}}_{T_{k}}}>\mu\abs{\Lambda_{N_k}}\\
    1&\mathrm{else}
  \end{cases}\;,
\end{align}
are interpreted row-by-row. We use $r=0.1 \Delta t$ without adaptation to improve convergence of the numerical scheme. In nonlinear dynamics, we have to apply an iterative algorithm to calculate unknowns at $t_{i+1}$ based on updating mass, stiffness, damping, gap functions and external forces. In case of no impact applying the fixed point iteration method, we obtain the unknowns at $t_{i+1}$ using the generalized-$\alpha$ algorithm for each time step:
\begin{align}
	\begin{aligned}
  \toprule 
	&\text{Generalized-$\alpha$ time integration scheme}\\
  \midrule 
  &\text{Set~} k=0\\
  &\vq^0_{i+1}=\vq_i,\;   \vv^0_{i+1}=\vv_i,\;   \va^0_{i+1}=\va_i,\; \vA^0_{i+1}=\vA_i\\
  &\vM^0_{i+1}=\vM_i,\; \vh^0_{i+1}=\vh_i,\; \vW^0_{N_{i+1}}=\vW_{N_i},\; \vW^0_{T_{i+1}}=\vW^0_{T_i},\;\vK^0_{i+1}=\vK_i,\; \vC^0_{i+1}=\vC_i \\		
  &\mathbf{while~} \mathbf{true}\\
  &\hspace{0.5cm}\text{calculate~}\vlambda^{k}_{N_{i+1}}, \vlambda^{k}_{T_{i+1}}\text{~with Gauss-Newton algorithm}\\
 &\hspace{0.5cm}\va^{k+1}_{i+1}=\left[\widehat{\vM}^{k}_{i+1}\right]^{-1} \left(\widehat{\vR}^{k}_{i+1}+\vW^{k}_{N_{i+1}}\vlambda^{k}_{N_{i+1}}+\vW^{k}_{T_{i+1}}\vlambda^{k}_{T_{i+1}}\right)\\
  &\hspace{0.5cm}\vA^{k+1}_{i+1}=\frac{1-\alpha_f}{1-\alpha_m}\va^{k+1}_{i+1}+\frac{\alpha_f}{1-\alpha_m}\va_i-\frac{\alpha_m}{1-\alpha_m}\vA_i\\
  &\hspace{0.5cm}\vv^{k+1}_{i+1}=\vv_i+\Delta t_i\left(1-\gamma\right)\vA_i+\Delta t_i\gamma\vA^{k+1}_{i+1}\\
  &\hspace{0.5cm}\vq^{k+1}_{i+1}=\vq_i+\Delta t_i\vv_i+\Delta t^2_i\left(0.5-\beta\right)\vA_i+\Delta t^2_i\beta\vA^{k+1}_{i+1}\\
  &\hspace{0.5cm}\mathbf{if~} \left\|\vv^{k+1}_{i+1}-\vv^{k}_{i+1}\right\|<\mathbf{tol~} \vv^{-}_{i+1}=\vv^{k+1}_{i+1}~ \mathbf{break~}\\
  &\hspace{0.5cm}\vM^{k+1}_{i+1}=\vM\left( \vq^{k+1}_{i+1}\right), \; \vh^{k+1}_{i+1}=\vh\left( \vq^{k+1}_{i+1},\vv^{k+1}_{i+1}\right) \\
   &\hspace{0.5cm}\vC^{k+1}_{i+1}=\vC\left( \vq^{k+1}_{i+1},\vv^{k+1}_{i+1}\right), \; \vK^{k+1}_{i+1}=\vK\left( \vq^{k+1}_{i+1},\vv^{k+1}_{i+1}\right) \\
  &\hspace{0.5cm}\vW^{k+1}_{N_{i+1}}=\vW_N\left( \vq^{k+1}_{i+1}\right), \; \vW^{k+1}_{T_{i+1}}=\vW_T\left( \vq^{k+1}_{i+1}\right)\\
  &\hspace{0.5cm} k = k+1\\
		&\mathbf{end} \\
  \bottomrule
	\end{aligned}\label{eq:fixed_point_general_alpha}
\end{align}
\subsection{Bathe-method}
The Bathe-method \cite{Bat12} is an effective implicit time integration scheme, which has been proposed for the finite element solution of nonlinear problems in structural dynamics. Various important attributes have been demonstrated. In particular, it has been shown that the scheme remains stable without the use of adjustable parameters. For this method, the complete time step $\Delta t$ is subdivided into two equal sub-steps. For the first sub-step the trapezoidal rule is used and for the second sub-step the 3-point Euler backward method is employed, as it is described in \eqref{eq:bathe_half_velocity} to \eqref{eq:bathe_acceleration}.
\subsubsection{General characteristics}
\begin{figure}[ht]
  \centering
  \includegraphics[width=\columnwidth]{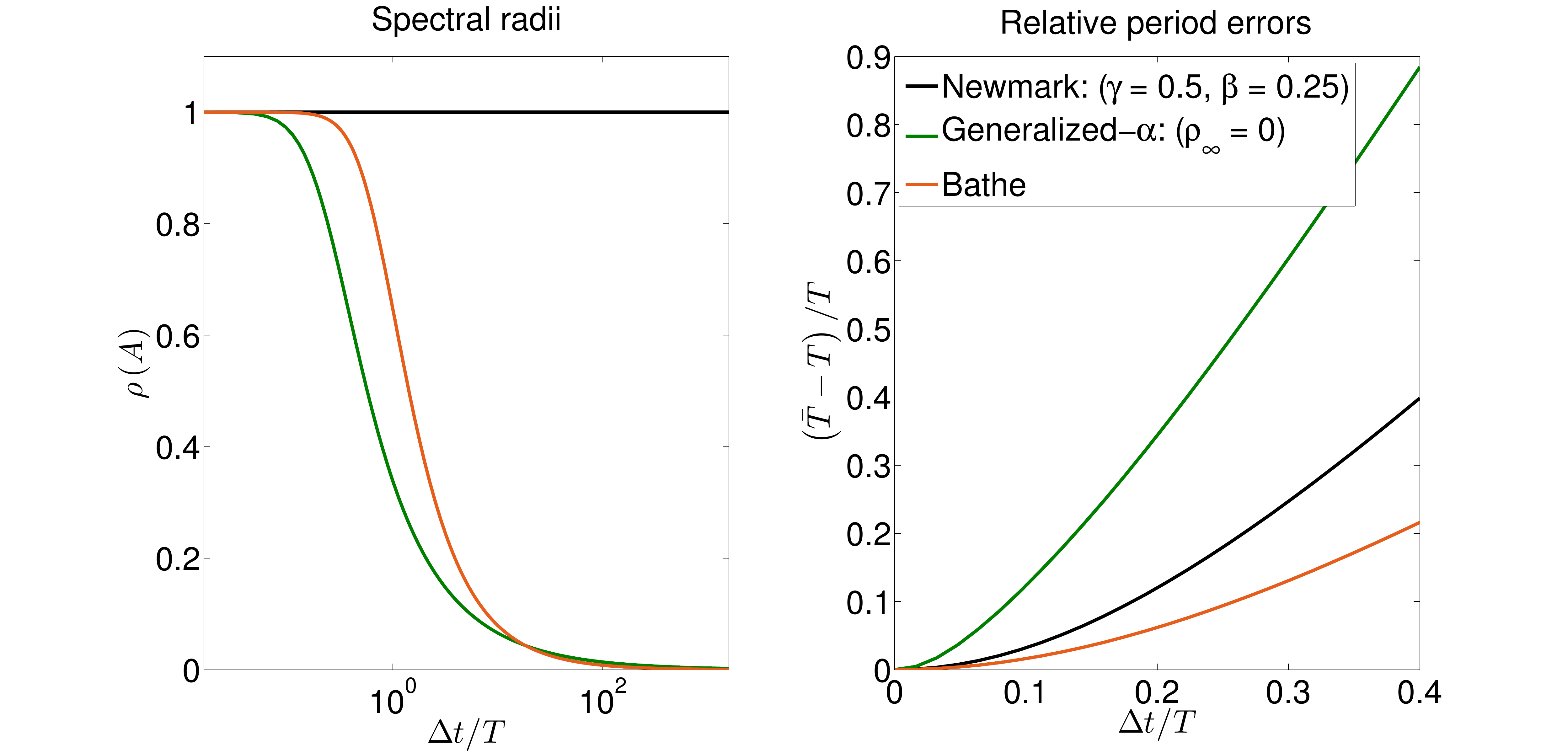}
  \caption{Comparison of second-order methods: Newmark method, generalized-$\alpha$ method and Bathe-method.}
  \label{fig:bathe_newmark_compare}
\end{figure}

According to \cite{Bat12}, we have to consider the dynamic equilibrium for time $t+\Delta t$ and $t+\Delta t/2$ which is indicated using index $i+1$ and index $i+1/2$:
\begin{align}
  &\vM_{i+\frac{1}{2}} \va_{i+\frac{1}{2}} + \vC_{i+\frac{1}{2}} \vv_{i+\frac{1}{2}} + \vK_{i+\frac{1}{2}} \vq_{i+\frac{1}{2}} = \vh_{i+\frac{1}{2}} + \vW_{N_{i+\frac{1}{2}}}\vlambda_{N_{i+\frac{1}{2}}} + \vW_{T_{i+\frac{1}{2}}}\vlambda_{T_{i+\frac{1}{2}}}\;,\label{eq:mechanical_system_i_half}\\
  &\vM_{i+1} \va_{i+1} + \vC_{i+1} \vv^{-}_{i+1} + \vK_{i+1} \vq_{i+1} = \vh^{-}_{i+1}+ \vW_{N_{i+1}}\vlambda_{N_{i+1}} + \vW_{T_{i+1}}\vlambda_{T_{i+1}}\;,\label{eq:mechanical_system_i_1_bathe}\\
  &\vv_{i+\frac{1}{2}} = \vv^{+}_i+\frac{\Delta t_i}{4}\left( \va_i + \va_{i+\frac{1}{2}} \right) \;, \label{eq:bathe_half_velocity}\\
  &\vq_{i+\frac{1}{2}} = \vq_i+\frac{\Delta t_i}{4}\left( \vv^{+}_i + \vv_{i+\frac{1}{2}} \right) \;, \label{eq:bathe_half_position}\\
  &\vv^{-}_{i+1} = \frac{1}{\Delta t_i}\vq_i - \frac{4}{\Delta t_i}\vq_{i+\frac{1}{2}} + \frac{3}{\Delta t_i}\vq_{i+1} \;, \label{eq:bathe_velocity}\\
  &\va_{i+1} = \frac{1}{\Delta t_i}\vv^{+}_i - \frac{4}{\Delta t_i}\vv_{i+\frac{1}{2}} + \frac{3}{\Delta t_i}\vv^{-}_{i+1} \;, \label{eq:bathe_acceleration}\\
  &\va_0 = \vM^{-1}_0\left(\vh_0 - \vK_0\vq_0 - \vC_0\vv_0\right) \;. \label{eq:bathe_mechanical_system_initial_acceleration}
\end{align}\label{eq:bathe}
In Fig.~\ref{fig:bathe_newmark_compare} for $\rho_{\infty}=0$, we compare the Bathe-method and the generalized-$\alpha$ method. We notice that the Bathe-method preserves the low-frequency oscillation better than the similar generalized-$\alpha$ method, and of course both avoid high-frequency vibrations. Total annihilation is important to avoid high frequency noises generated by the time integration algorithm. The Bathe-method shows a better behavior for the period error in comparison to the generalized-$\alpha$ method and the classic Newmark method.
Further properties can be found in~\ref{sec:appB}.
\subsubsection{Calculation of contact forces on velocity level} \label{sec:lambda_bathe}
Using the Bathe-method \eqref{eq:bathe_half_velocity}-\eqref{eq:bathe_acceleration}, we have to modify the implicit representation of gap velocities in \eqref{eq:implicit_velocity_gap}. We have to solve for the unknowns first at $t+\Delta t / 2$ and then at $t+\Delta t$:
\begin{align}
  &\dot{\vg}_{N_{i+1/2}}=\vF_{N_{i+1/2}}+\frac{4}{\Delta t} \left\{\vG_{N_{i+1/2}}\vlambda_{N_{i+1/2}}+\vG_{NT_{i+1/2}}\vlambda_{T_{i+1/2}}\right\}\;,\\
  &\dot{\vg}_{T_{i+1/2}}=\vF_{T_{i+1/2}}+\frac{4}{\Delta t} \left\{\vG_{TN_{i+1/2}}\vlambda_{N_{i+1/2}}+\vG_{T_{i+1/2}}\vlambda_{T_{i+1/2}}\right\}\;,\end{align}
where
\begin{align}
  \vF_{N_{i+1/2}}&=\vW_{N_{i+1/2}}^T \left\{-\vv_i-\frac{4}{\Delta t}\vq_i+\frac{4}{\Delta t}\widehat{\vK}^{-1}_1 \widehat{\vR}_1 \right\}\;,\\
  \vF_{T_{i+1/2}}&=\vW_{T_{i+1/2}}^T \left\{-\vv_i-\frac{4}{\Delta t}\vq_i+\frac{4}{\Delta t}\widehat{\vK}^{-1}_1 \widehat{\vR}_1 \right\}\;,
\end{align}
and
\begin{align}
  \vG_{N_{i+1/2}}&=\vW^T_{N_{i+1/2}}\widehat{\vK}^{-1}_{1}\vW_{N_{i+1/2}}\;, \quad \vG_{NT_{i+1/2}}=\vW^T_{N_{i+1/2}}\widehat{\vK}^{-1}_{1}\vW_{T_{i+1/2}}\;,\\
  \vG_{T_{i+1/2}}&=\vW^T_{T_{i+1/2}}\widehat{\vK}^{-1}_{1}\vW_{T_{i+1/2}}\;, \quad \vG_{TN_{i+1/2}}=\vW^T_{T_{i+1/2}}\widehat{\vK}^{-1}_{1}\vW_{N_{i+1/2}}\;.
\end{align}
After calculation of the unknowns at $t_{i+1/2}$ for the second half of the interval, we have:
\begin{align}
  &\dot{\vg}^{-}_{N_{i+1}}=\vF_{N_{i+1}}+\frac{3}{\Delta t} \left\{\vG_{N_{i+1}}\vlambda_{N_{i+1}}+\vG_{NT_{i+1}}\vlambda_{T_{i+1}}\right\}\;,\\
  &\dot{\vg}^{-}_{T_{i+1}}=\vF_{T_{i+1}}+\frac{3}{\Delta t} \left\{\vG_{TN_{i+1}}\vlambda_{N_{i+1}}+\vG_{T_{i+1}}\vlambda_{T_{i+1}}\right\}\;,
\end{align}
where
\begin{align}
  \vF_{N_{i+1}}&=\vW_{N_{i+1}}^T \left\{\frac{1}{\Delta t}\vq_i-\frac{4}{\Delta t}\vq_{i+1/2}+\frac{3}{\Delta t}\widehat{\vK}^{-1}_2 \widehat{\vR}_2 \right\}\;,\\
  \vF_{T_{i+1}}&=\vW_{T_{i+1}}^T \left\{\frac{1}{\Delta t}\vq_i-\frac{4}{\Delta t}\vq_{i+1/2}+\frac{3}{\Delta t}\widehat{\vK}^{-1}_2 \widehat{\vR}_2 \right\}\;,
\end{align}
and
\begin{align}
\vG_{N_{i+1}}&=\vW^T_{N_{i+1}}\widehat{\vK}^{-1}_{2}\vW_{N_{i+1}}\;, \quad \vG_{NT_{i+1}}=\vW^T_{N_{i+1}}\widehat{\vK}^{-1}_{2}\vW_{T_{i+1}}\;,\\
\vG_{T_{i+1}}&=\vW^T_{T_{i+1}}\widehat{\vK}^{-1}_{2}\vW_{T_{i+1}}\;, \quad \vG_{TN_{i+1}}=\vW^T_{T_{i+1}}\widehat{\vK}^{-1}_{2}\vW_{N_{i+1}}\;.
\end{align}
The matrices $\widehat{\vK}_1$, $\widehat{\vK}_2$, $\widehat{\vR}_1$ and $\widehat{\vR}_2$ 
are defined in \eqref{eq:fixed_point_bathe_first}-\eqref{eq:fixed_point_bathe_second}. Once we have calculated the velocities $\vv^{-}_{i+1}$, the computation of the impulsive forces is given in Sect.~\ref{sec:Lambda}. The generalized-$\alpha$ method and the Bathe-method use the same procedure to update the calculated velocities after the impact. In case of no impacts, applying the fixed point iteration method for each time step of the Bathe-method, we get: 
\begin{align}
	\begin{aligned}
  \toprule
	&\text{Bathe time integration scheme: first half step}\\
  \midrule 
  &\text{Set~} k=0\\
  &\vq^0_{i+1/2}=\vq_i,\;  \vv^0_{i+1/2}=\vv_i,\;  \va^0_{i+1/2}=\va_i\;\\
  &\vM^0_{i+1/2}=\vM_i,\; \vh^0_{i+1/2}=\vh_i,\;\vW^0_{N_{i+1/2}}=\vW_{N_i},\; \vW^0_{T_{i+1/2}}=\vW^0_{T_i},\; \vK^0_{i+1/2}=\vK_i \\
  &\mathbf{while~} \mathbf{true}\\
  &\hspace{0.5cm}\text{calculate~}\vlambda^{k}_{N_{i+1/2}}, \vlambda^{k}_{T_{i+1/2}}\text{~with Gauss-Newton algorithm}\\
  &\hspace{0.5cm}\widehat{K}_1^k=\frac{16}{\Delta t^2}\vM^k_{i+1/2}+\frac{4}{\Delta t}\vC^k_{i+1/2}+\vK^k_{i+1/2}\\
  &\hspace{0.5cm}\widehat{\vR}_1^k=\vh^k_{i+1/2} + \vM^k_{i+1/2}\left( \frac{16}{\Delta t^2}\vq_i+\frac{8}{\Delta t}\vv_i+\va_i \right)+ \vC^k_{i+1/2}\left( \frac{4}{\Delta t}\vq_i+\vv_i \right)\\
  &\hspace{0.5cm}\vq^{k+1}_{i+1/2}=\left[\widehat{K}_1^k\right]^{-1} \left(\widehat{\vR}_1^k +\vW^k_{N_{i+1/2}}\vlambda^k_{N_{i+1/2}}+\vW^k_{T_{i+1/2}}\vlambda^k_{T_{i+1/2}}\right)\\
  &\hspace{0.5cm}\vv^{k+1}_{i+1/2}=-\vv_i+\frac{4}{\Delta t} \left(\vq^{k+1}_{i+1/2}-\vq_i \right)\\
  &\hspace{0.5cm}\va^{k+1}_{i+1/2}=-\va_i+\frac{4}{\Delta t} \left(\vv^{k+1}_{i+1/2}-\vv_i \right)\\
  &\hspace{0.5cm}\mathbf{if~} \left\|\vv^{k+1}_{i+1/2}-\vv^{k}_{i+1/2}\right\|<\mathbf{tol~}  \mathbf{break~}\\
	&\hspace{0.5cm} \text{update}~\vM^{k+1}_{i+1/2},\;\vh^{k+1}_{i+1/2},\;\vW^{k+1}_{N_{i+1/2}},\;\vW^{k+1}_{T_{i+1/2}}\\
  &\hspace{0.5cm} k = k+1\\
	&\mathbf{end}\\
  \bottomrule
	\end{aligned}\label{eq:fixed_point_bathe_first}
\end{align}
\begin{align}
	\begin{aligned}
  \toprule
	&\text{Bathe time integration scheme: second half step}\\
  \midrule 
  &\text{Set~} k=0\\
  &\vq^0_{i+1}=\vq_{i+1/2},\;   \vv^0_{i+1}=\vv_{i+1/2},\;   \va^0_{i+1}=\va_{i+1/2}\;\\
  &\vM^0_{i+1}=\vM_{i+1/2},\; \vh^0_{i+1}=\vh_{i+1/2},\;\vW^0_{N_{i+1}}=\vW_{N_{i+1/2}},\; \vW^0_{T_{i+1}}=\vW^0_{T_{i+1/2}},\;\vK^0_{i+1}=\vK_{i+1/2} \\
  &\mathbf{while~} \mathbf{true}\\
  &\hspace{0.5cm}\text{calculate~}\vlambda^{k}_{N_{i+1}}, \vlambda^{k}_{T_{i+1}}\text{~with Gauss-Newton algorithm}\\
  &\hspace{0.5cm}\widehat{K}_2^k=\frac{9}{\Delta t^2}\vM^k_{i+1}+\frac{3}{\Delta t}\vC^k_{i+1}+\vK^k_{i+1}\\
  &\hspace{0.5cm}\widehat{\vR}_2^k=\vh^k_{i+1}+ \vM^k_{i+1}\left( \frac{12}{\Delta t^2}\vq_{i+1/2}-\frac{3}{\Delta t^2}\vq_i+\frac{4}{\Delta t}\vv_{i+1/2}-\frac{1}{\Delta t}\vv_i \right)+\vC^k_{i+1}\left( \frac{4}{\Delta t}\vq_{i+1/2}-\frac{1}{\Delta t}\vq_i \right)\\
  &\hspace{0.5cm}\vq^{k+1}_{i+1}=\left[\widehat{K}_2^k\right]^{-1} \left(\widehat{\vR}_2^k +\vW^k_{N_{i+1}}\vlambda^k_{N_{i+1}}+\vW^k_{T_{i+1}}\vlambda^k_{T_{i+1}}\right)\\
  &\hspace{0.5cm}\vv^{k+1}_{i+1}=\frac{1}{\Delta t}\vq_i - \frac{4}{\Delta t}\vq_{i+1/2} + \frac{3}{\Delta t}\vq^{k+1}_{i+1} \\
  &\hspace{0.5cm}\va^{k+1}_{i+1}=\frac{1}{\Delta t}\vv_i - \frac{4}{\Delta t}\vv_{i+1/2} + \frac{3}{\Delta t}\vv^{k+1}_{i+1} \\
  &\hspace{0.5cm}\mathbf{if~} \left\|\vv^{k+1}_{i+1}-\vv^{k}_{i+1}\right\|<\mathbf{tol~} \vv^{-}_{i+1}=\vv^{k+1}_{i+1}~ \mathbf{break~}\\
	&\hspace{0.5cm} \text{update}~\vM^{k+1}_{i+1/2},\;\vh^{k+1}_{i+1/2},\;\vW^{k+1}_{N_{i+1/2}},\;\vW^{k+1}_{T_{i+1/2}}\\
  &\hspace{0.5cm} k = k+1\\
	&\mathbf{end} \\
  \bottomrule
	\end{aligned}\label{eq:fixed_point_bathe_second}
\end{align}
\subsection{ED-$\valpha$ method}
In energy decaying schemes, we develop robust algorithms for integrating stiff nonlinear finite element problems in time. The basic motivation behind these schemes is that classical algorithms that are unconditionally stable and high frequency dissipate in the linear regime, loose their properties in the nonlinear regime~\cite{Bot04}. We typically interpret the stages $\vq_j$, $\vv_j$, $\va_j$ as field variables associated with the time $t^{+}_{i}$. In this sense, the unknown fields are allowed to create a jump discontinuity at the beginning of the time step that is responsible for the high frequency damping behavior of the scheme~\cite{Bot97}.
\begin{align}
  &\vM_{j} \va_{j} + \vC_{j} \vv_{j} + \vK_{j} \vq_{j} = \vh_{j} + \vW_{N_{j}}\vlambda_{N_{j}} + \vW_{T_{j}}\vlambda_{T_{j}}\;,\label{eq:mechanical_system_j}\\
  &\vM_{i+1} \va_{i+1} + \vC_{i+1} \vv^{-}_{i+1} + \vK_{i+1} \vq_{i+1} = \vh^{-}_{i+1}+ \vW_{N_{i+1}}\vlambda_{N_{i+1}} + \vW_{T_{i+1}}\vlambda_{T_{i+1}}\;,\label{eq:mechanical_system_i_1_ed}\\
  &\vv_{j} = \vv^{+}_i+\Delta t_i\alpha_{AR}\left[\alpha\left(\va_j - \va_{i}\right) - \va_{i+1} + \va_{i} \right] \;, \label{eq:ed_half_velocity}\\
  &\vv^{-}_{i+1} = \vv^{+}_i+\frac{\Delta t_i}{2}\left( \va_j + \va_{i+1} \right) \;, \label{eq:ed_velocity}\\
  &\vq_{j} = \vq_i+\Delta t_i\alpha_{AR}\left[\alpha\left(\vv_j - \vv^{+}_{i}\right) - \vv_{i+1}^{-} + \vv^{+}_{i} \right] \;, \label{eq:ed_half_position}\\
  &\vq_{i+1} = \vq_i+\frac{\Delta t_i}{2}\left( \vv_j + \vv^{-}_{i+1} \right) \;, \label{eq:ed_position}\\
  &\va_0 = \vM^{-1}_0\left(\vh_0 - \vK_0\vq_0 - \vC\right) \;. \label{eq:ed_mechanical_system_initial_acceleration}
\end{align}
In this work for practical implementation of the scheme, the displacements $\vq_j$ and $\vq_{i+1}$ are eliminated, leaving a velocity-based iteration scheme in the $2 \times \text{ndof}$ unknowns $v_j$ and $v_{i+1}$. However, the overall procedure is more expensive than other one-stage schemes like the generalized-$\alpha$ method, since the matrices are twice as large.\par
\begin{figure}[ht]
  \centering
  \includegraphics[width=\columnwidth]{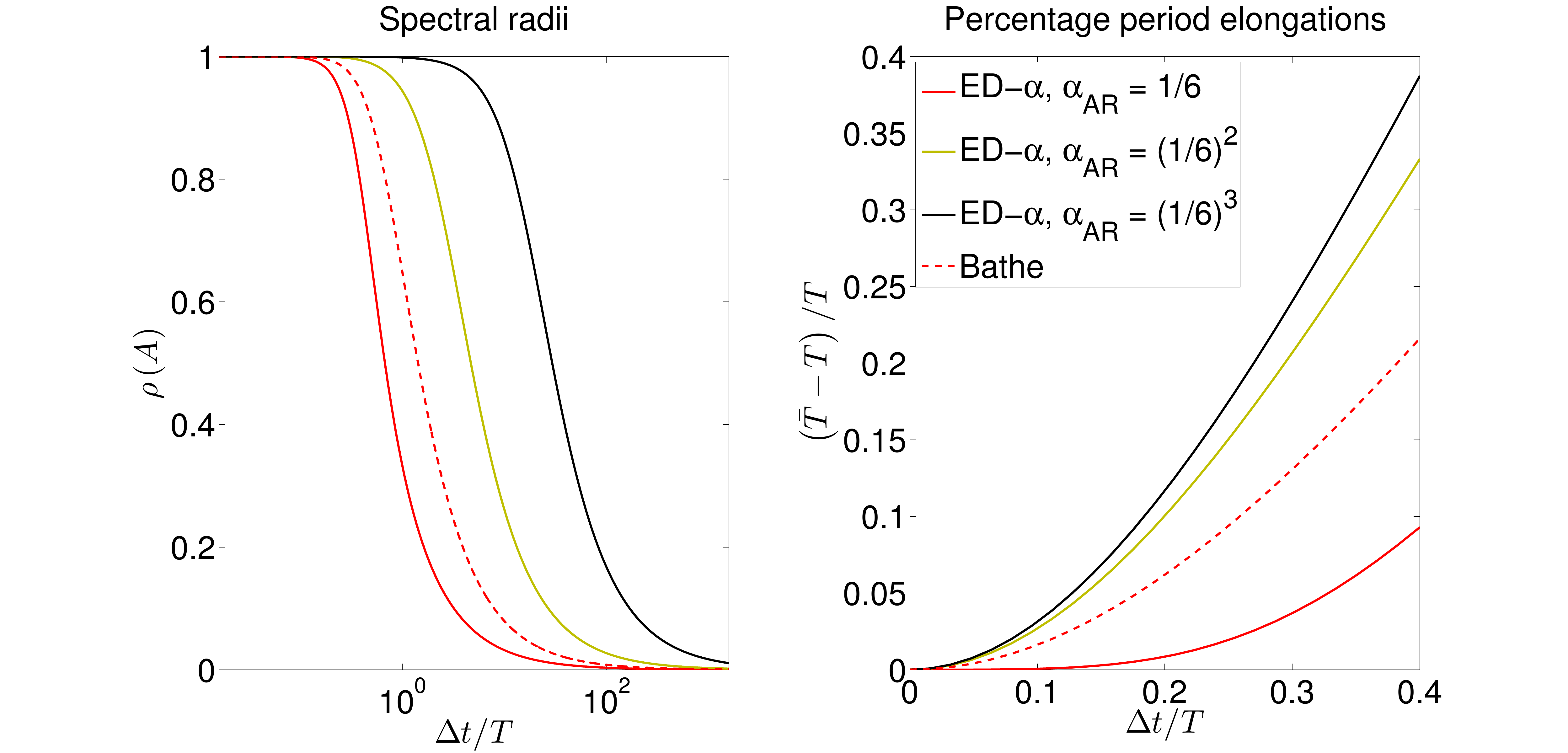}
  \caption{Comparison of methods: ED-$\alpha$ method with $\rho_{\infty}=0$ and Bathe-method.}
  \label{fig:ed_bathe_compare}
\end{figure}
\begin{figure}[ht]
  \centering
  \includegraphics[width=\columnwidth]{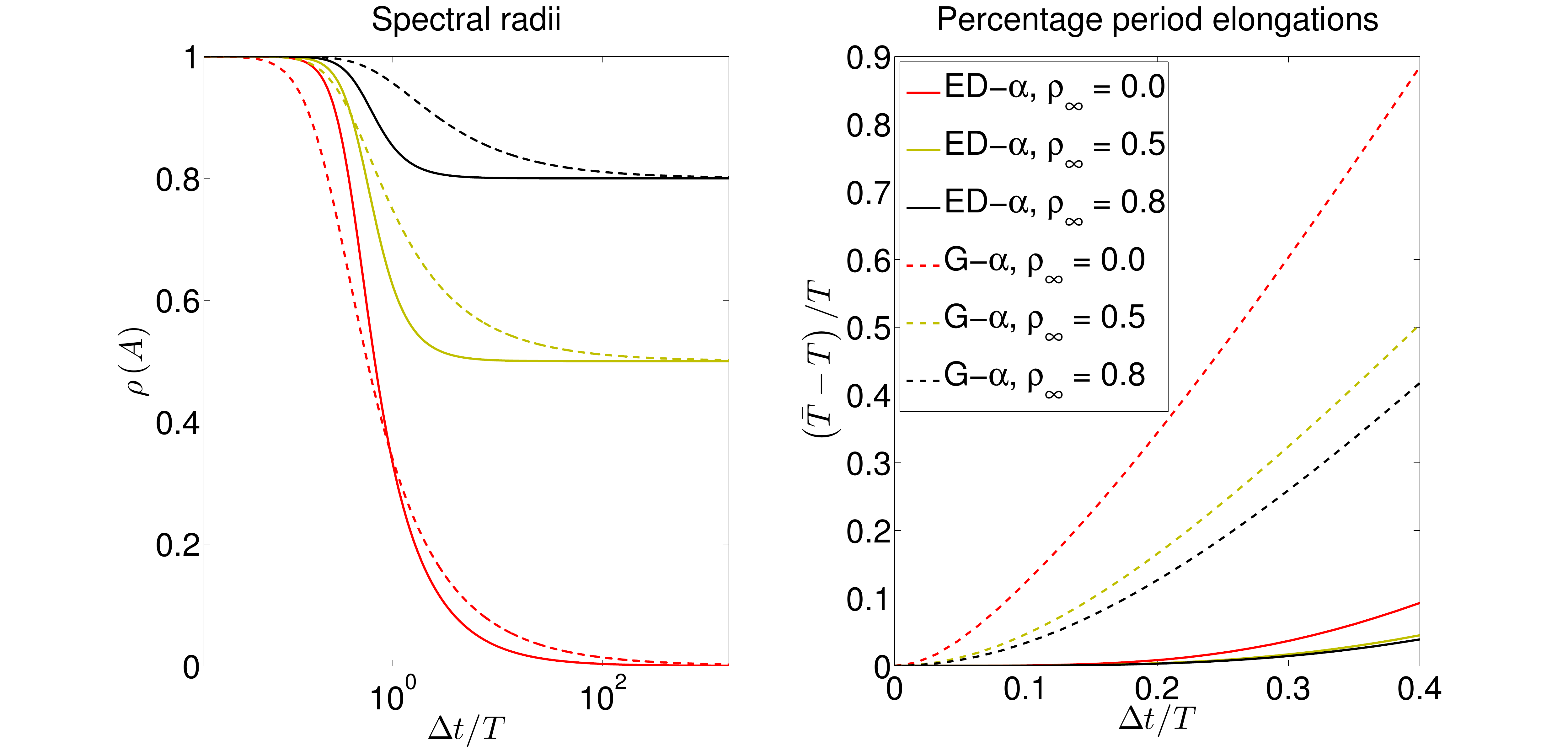}
  \caption{Comparison of methods: ED-$\alpha$ method with $\alpha_{AR}=1/6$ and generalized-$\alpha$ method}
  \label{fig:ed_ga_compare}
\end{figure}

Note that for $\alpha_{AR}=0$ or $\alpha=0$, we recover a conserving scheme. The parameter $\alpha_{AR}$ does not control the asymptotic value of the spectral radius but only controls the cut-off frequency of the scheme and so relative period errors (Fig.~\ref{fig:ed_bathe_compare}). The minimum period elongation is obtained for $\alpha_{AR}=1/6$. The method is second order accurate for arbitrary $\alpha_{AR}\geq 0$ and arbitrary ordinary differential equations; third order accuracy is obtained for the scalar linear model problem in Sect.~\ref{sec:model_problem} and the special value $\alpha_{AR}=1/6$. The parameter $\alpha$ is responsible for the asymptotic value of the spectral radius, in fact \eqref{eq:ed_alpha_par} is an optimal choice~\cite{Bot04} (Fig.~\ref{fig:ed_ga_compare}):
\begin{align}
  \begin{aligned}
    \alpha&=\frac{1-\rho_{\infty}}{1+\rho_{\infty}} \;, \\
    \alpha_{AR}&=\frac{1}{6} \;. 
\end{aligned}\label{eq:ed_alpha_par}
\end{align} 
In case of no impact using \eqref{eq:mechanical_system_j} to \eqref{eq:ed_position} we write:
\begin{figure}[ht]
  \centering
  \footnotesize
  \def\svgwidth{0.6\columnwidth}
  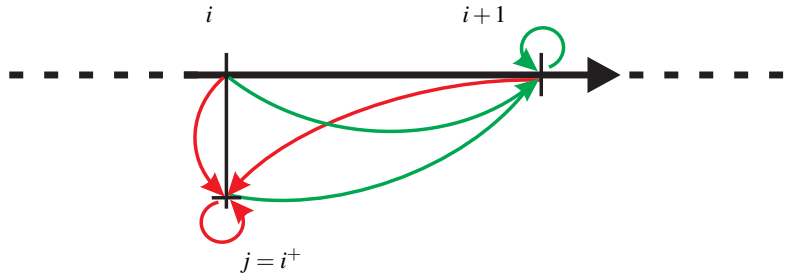
  \caption{ED-$\alpha$ method in time.}
  \label{fig:ed_vector}
\end{figure}
\begin{align}
	\begin{aligned}
  \toprule
	&\text{ED-$\valpha$ time integration scheme}\\
  \midrule 
  &\text{Set~} k=0\\
  &\vq^0_{j}=\vq_i\;,\; \vv^0_{j}=\vv_i\;,\; \va^0_{j}=\va_i\\
  &\vM^0_{j}=\vM_i\;,\; \vh^0_{j}=\vh_i\;,\; \vW^0_{N_{j}}=\vW_{N_i}\;,\; \vW^0_{T_{j}}=\vW^0_{T_i}\;,\; \vK^0_{j}=\vK_i \\ 
  &\mathbf{while~} \mathbf{true}\\
  &\hspace{0.5cm}\text{calculate~}\vlambda^{k}_{N_{j}},\;\vlambda^{k}_{T_{j}},\;\vlambda^{k}_{N_{i+1}},\;\vlambda^{k}_{T_{i+1}} \text{~with Gauss-Newton algorithm} \\		
  &\hspace{0.5cm}\vv^{k+1}_{c}=\left[\widehat{\vM}_c^k\right]^{-1} \left(\widehat{\vR}_c^k +\vW^k_{N_{c}}\vlambda^k_{N_{c}}+\vW^k_{T_{c}}\vlambda^k_{T_{c}}\right)\\
  &\hspace{0.5cm}\vq^{k+1}_{i+1}=c_1\vq_i+c_2\vv^{k+1}_{j}+c_3\vv^{k+1}_{i+1} \\
  &\hspace{0.5cm}\va^{k+1}_{i+1}=c_4\vv^{k+1}_{i+1}+c_5\vv_{i}+c_6\vv^{k+1}_{j}+c_7\va_{i} \\
	&\hspace{0.5cm}\vq^{k+1}_{j}=c_8\vv^{k+1}_{i+1}+c_9\vq_{i}+c_{10}\vv^{k+1}_{j}+c_{11}\vv_{i} \\
  &\hspace{0.5cm}\va^{k+1}_{j}=c_{12}\vv^{k+1}_{j}+c_{13}\vv_{i}+c_{14}\va_{i}+c_{15}\vv^{k+1}_{i+1} \\
  &\hspace{0.5cm}\mathbf{if~} \left\|\vv^{k+1}_{i+1}-\vv^{k}_{i+1}\right\|<\mathbf{tol~} \vv^{-}_{i+1}=\vv^{k+1}_{i+1}~\mathbf{break~}\\
	&\hspace{0.5cm} \text{update}~\vM^{k+1}_{i+1}\;,\;\vh^{k+1}_{i+1}\;,\;\vW^{k+1}_{N_{i+1}}\;,\;\vW^{k+1}_{T_{i+1}}\\
	&\hspace{0.5cm} \text{update}~\vM^{k+1}_{j}\;,\;\vh^{k+1}_{j}\;,\;\vW^{k+1}_{N_{j}}\;,\;\vW^{k+1}_{T_{j}}\\
  &\hspace{0.5cm} k = k+1\\
	&\mathbf{end}\\
  \bottomrule
	\end{aligned}\label{eq:fixed_point_ed}
\end{align}
The matrices $\widehat{\vM}_c$, $\widehat{\vR}_c$, $\vW_{N_{c}}$, $\vW_{T_{c}}$ and vectors $\vv_{c}$, $\vlambda_{c}$ are defined as
\begin{align}
  &\widehat{\vM}_c=\begin{pmatrix} c_4\vM_{i+1}+c_3\vK_{i+1} & c_6\vM_{i+1}+c_2\vK_{i+1} \\ c_{15}\vM_j+c_8\vK_j & c_{12}\vM_j+c_{10}\vK_j \end{pmatrix} \;, \\
  &\widehat{\vR}_c=\begin{pmatrix} \vh_{i+1}-c_5\vM_{i+1}\vv_{i}-c_7\vM_{i+1}\va_{i}-c_1\vK_{i+1}\vq_{i} \\ \vh_{j}-c_{13}\vM_{j}\vv_{i}-c_{14}\vM_{j}\va_{i}-c_9\vK_{j}\vq_{i}-c_{11}\vK_{j}\vv_{i} \end{pmatrix} \;, \\
  &\vW_{N_{c}}=\begin{pmatrix} \vW_{N_{i+1}} & \vnull \\ \vnull & \vW_{N_{j}} \end{pmatrix} \;, \quad \vW_{T_{c}}=\begin{pmatrix} \vW_{T_{i+1}} & \vnull \\ \vnull & \vW_{T_{j}} \end{pmatrix}\;, \\
  &\vv_{c}=\begin{pmatrix} \vv_{i+1} \\ \vv_{j} \end{pmatrix} \;, \quad \vlambda_{c}=\begin{pmatrix} \vlambda_{i+1} \\ \vlambda_{j} \end{pmatrix} \; .
\end{align}
The coefficients $c_1$ to $c_{15}$ are
\begin{align}
  \begin{aligned}
    c_1 &= 1\;, & &\\
    c_2 &= \frac{\Delta t}{2}\;, & c_3 &= \frac{\Delta t}{2}\;, \\
    c_4 &= \frac{2~\alpha_{AR}~\alpha}{\Delta t~\alpha_{AR}\left(\alpha+1\right)}\;, & c_5 &= \frac{2~\left(0.5-\alpha_{AR}~\alpha\right)}{\Delta t~\alpha_{AR}\left(\alpha+1\right)}\;,\\
    c_6 &= \frac{-1}{\Delta t~\alpha_{AR}\left(\alpha+1\right)}\;, & c_7 &= \frac{\alpha_{AR}\left(1-\alpha\right)}{\Delta t~\alpha_{AR}\left(\alpha+1\right)}\;, \\
    c_8 &= -\Delta t~\alpha_{AR}\;, & c_9 &= 1\;,\\
    c_{10} &= \Delta t~\alpha_{AR}~\alpha\;, & c_{11} &= \Delta t~\alpha_{AR}\left(1-\alpha\right)\;,\\
    c_{12} &= \frac{1}{\Delta t~\alpha_{AR}\left(\alpha+1\right)}\;, & c_{13} &= \frac{-2~\left(0.5+\alpha_{AR}\right)}{\Delta t~\alpha_{AR}\left(\alpha+1\right)}\;,\\
    c_{14} &= \frac{-\alpha_{AR}\left(1-\alpha\right)}{\Delta t~\alpha_{AR}\left(\alpha+1\right)}\;, & c_{15} &= \frac{2~\alpha_{AR}}{\Delta t~\alpha_{AR}\left(\alpha+1\right)}\;.
  \end{aligned}
\end{align}
\subsection{Comparison of base integration schemes}\label{sec:model_problem}
\begin{figure}[ht]
  \centering
  \footnotesize
  \def\svgwidth{0.6\columnwidth}
  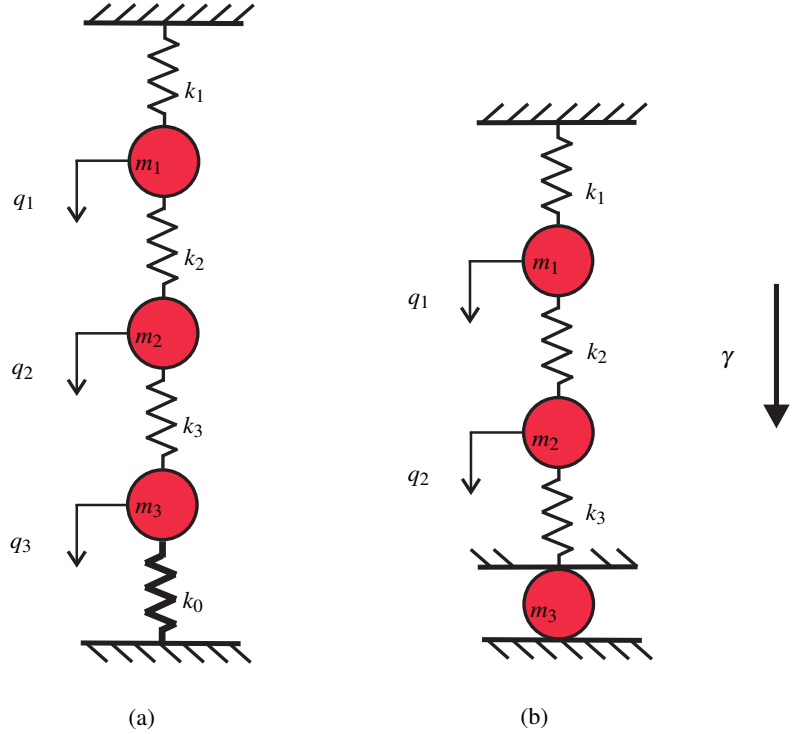
  \caption{Model problem of a three degree of freedom mass-spring system, (a): model with stiff spring, (b): model with bilateral constraint.}
  \label{fig:balls}
\end{figure}

Our objective in this section is to present the solution of a simple linear system as a \emph{model problem} to represent the stiff and flexible parts. We compare the presented base integration schemes and discuss their properties within the overall framework. We also compare the results for bilateral contact forces using the acceleration level and velocity level approach to satisfy the constraints. Let us consider the solution of the 3 degree of freedom mass-spring system shown in Fig.~\ref{fig:balls} for which the governing equations are
\begin{align}
&\begin{pmatrix}m_1 & 0 & 0 \\ 0 & m_2 & 0 \\ 0 & 0 & m_3\end{pmatrix}\begin{pmatrix}\ddot{q}_1 \\ \ddot{q}_2 \\ \ddot{q}_3\end{pmatrix}+\begin{pmatrix}k_1+k_2 & -k_2 & 0 \\ -k_2 & k_2+k_3 & -k_3 \\ 0 & -k_3 & k_3+k_0\end{pmatrix}\begin{pmatrix}q_1\\q_2\\q_3\end{pmatrix}=\begin{pmatrix}m_1\gamma\\m_2\gamma\\m_3\gamma\end{pmatrix}\;.\label{eq:mass_spring_eq}
\end{align}
We use the parameters $k_1=\unit[1]{N/m}$, $k_2=\unit[1]{N/m}$, $k_3=\unit[1]{N/m}$, $k_0=\unit[10^7]{N/m}$, $m_1= \unit[1]{kg}$, $m_2= \unit[1]{kg}$, $m_3=\unit[1]{kg}$, $\gamma=\unit[9.81]{m/s^2}$. In the left part of Fig.~\ref{fig:balls}, a stiff spring is used to represent, for example, an almost rigid connections~\cite{Bat12}, while the other springs represent the flexible parts of the structural model.\par
\begin{figure}[ht]
  \centering
  \includegraphics[width=\columnwidth]{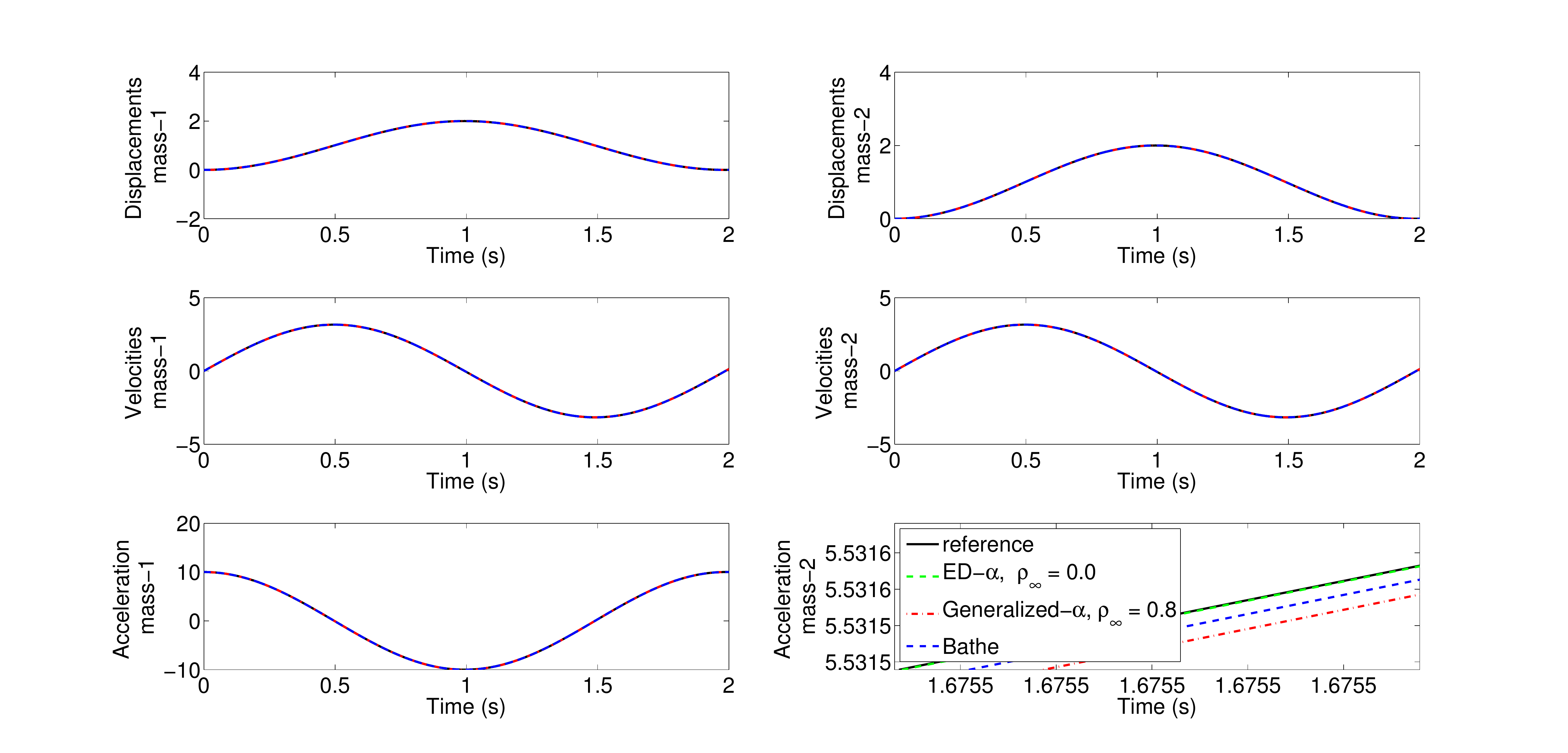}
  \caption{Displacement, velocity and acceleration of masses 1 and 2 for various integration methods.}
  \label{fig:balls_a_1}
\end{figure}

\begin{figure}[ht]
  \centering
  \includegraphics[width=1\columnwidth]{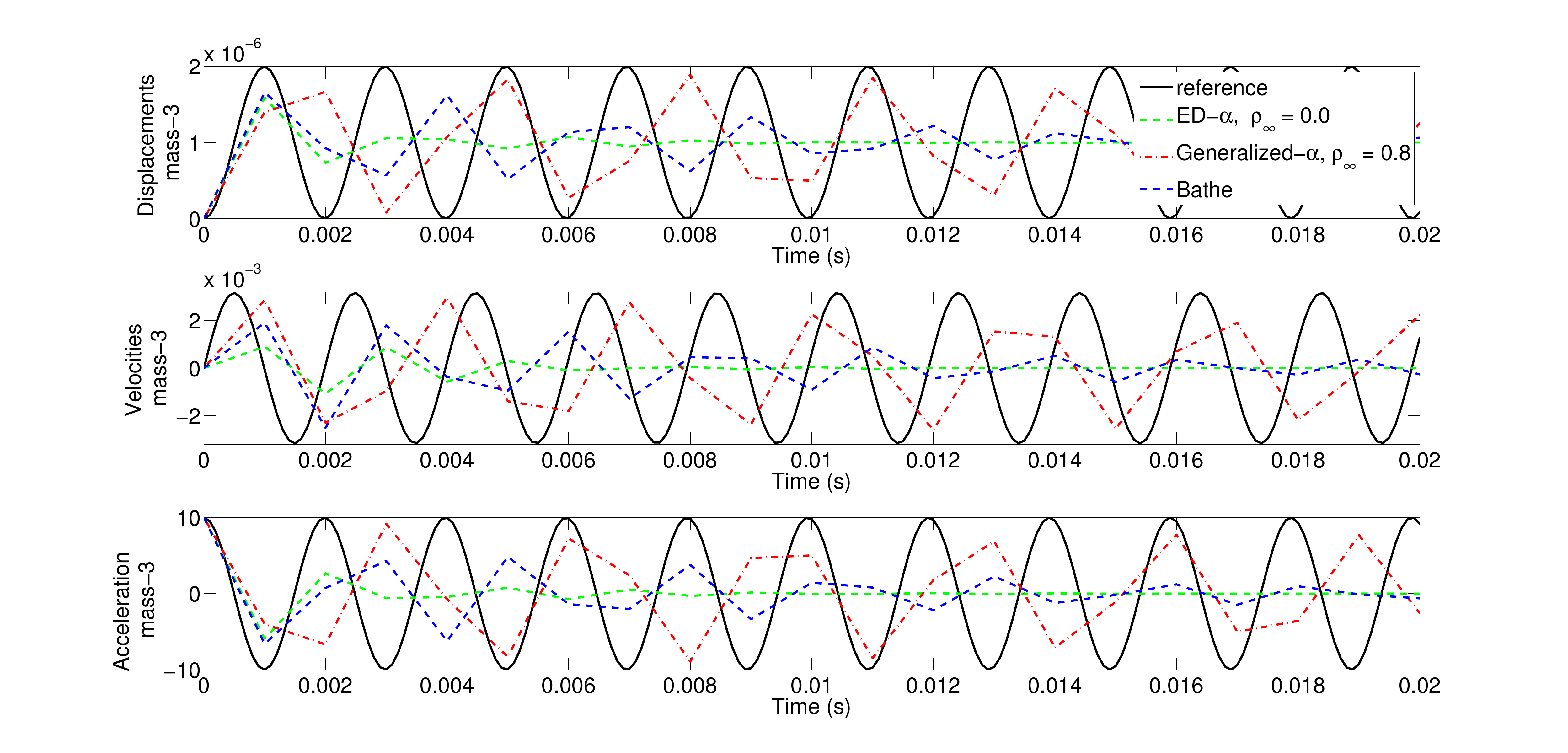}
  \caption{Displacement of mass 3 for various integration methods.}
  \label{fig:balls_a_2_2}
\end{figure}

The highest frequency in case (a) is $f_{\text{max}} = \unit[503.3]{Hz}$, which is due to the stiff spring with stiffness $k_0$. The exact trajectory for mass number 3 is plotted using $\Delta t = \unit[10^{-4}]{s}$ to represent it more properly, and the time-step size for comparison of different base integration schemes is chosen as $\Delta t = \unit[10^{-3}]{s}$. Figure~\ref{fig:balls_a_1} shows the response of the first two degrees of freedom which are connected with the soft springs. As we observe, the response from different base integration schemes follows the reference solution coming from modal analysis. The only important difference is the phase shift due to numerical integration which is minimum for the ED-$\alpha$ method. The response of mass number 3 which is connected to the rigid wall with a stiff spring is shown in Fig.~\ref{fig:balls_a_2_2}. The generalized-$\alpha$ method with $\rho_{\infty}=0.8$ is able to represent the high frequency vibration with some phase shift, however the other methods damp out this high frequency mode in a few time steps. Using $\Delta t = \unit[10^{-3}]{s}$, i.e., a sampling frequency $f_s = \unit[1000]{Hz}$, we have to make sure that all frequencies higher than $\frac{f_s}{2} = \unit[500]{Hz}$ are damped out, in order to avoid chattering in the response.\par
\begin{figure}[ht]
  \centering
  \includegraphics[width=\columnwidth]{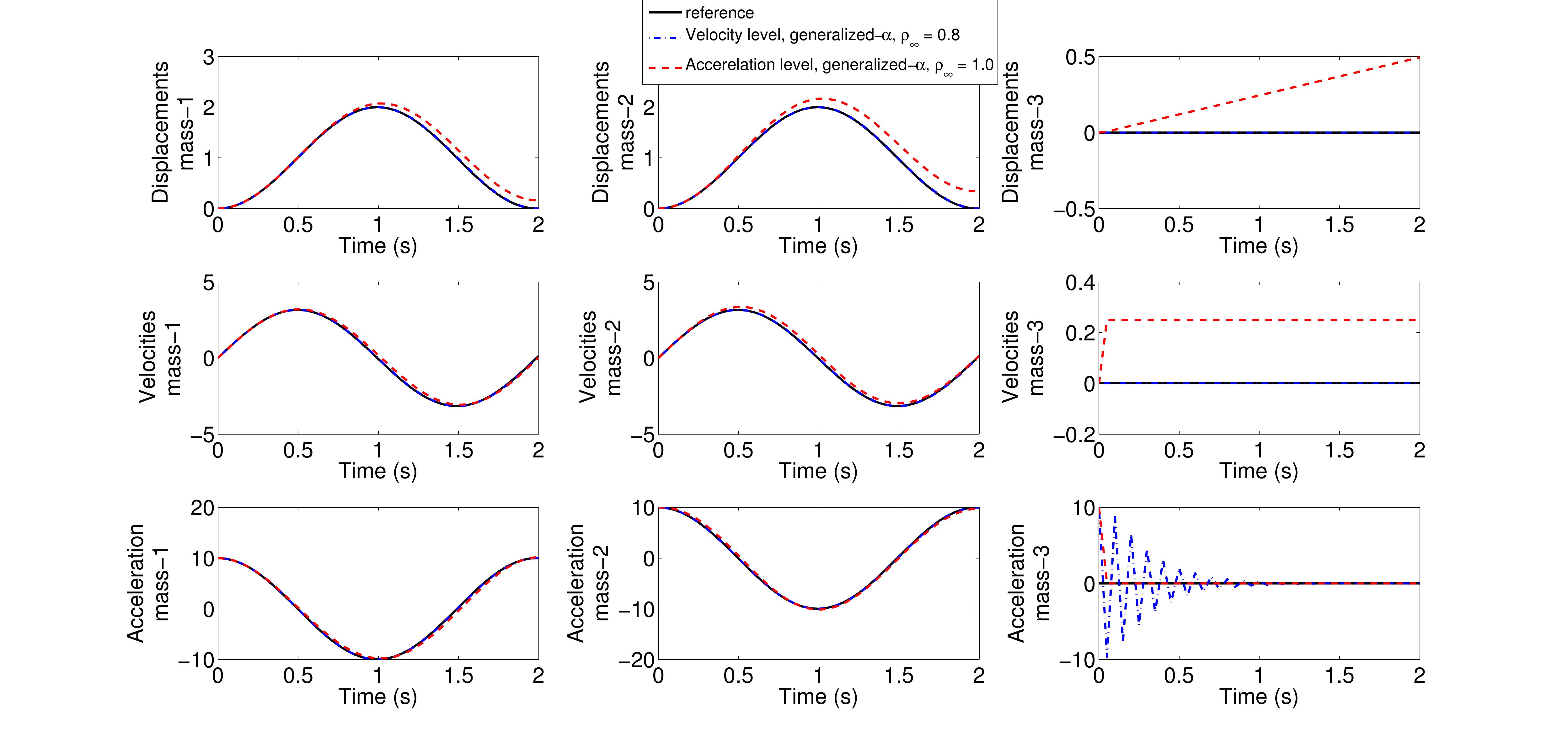}
  \caption{Displacement, velocity and acceleration for mass 1, 2 and 3.}
  \label{fig:balls_b_1}
\end{figure}

\begin{figure}[ht]
  \centering
  \includegraphics[width=\columnwidth]{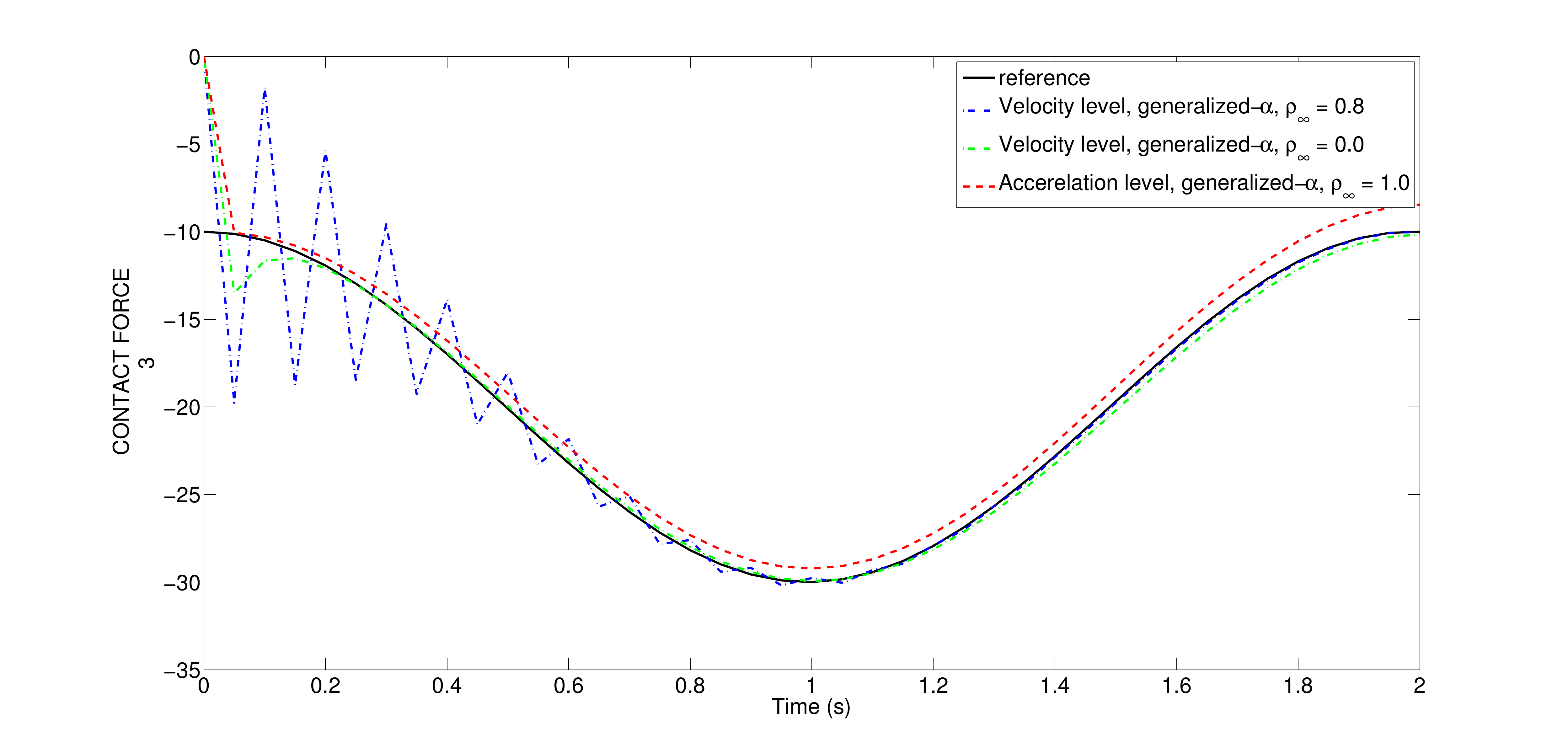}
  \caption{Comparison of acceleration and velocity level for the calculation of contact forces.}
  \label{fig:balls_b_2}
\end{figure}

Figure~\ref{fig:balls_b_1} shows the difference between simulation results if we try to satisfy the constraint on acceleration and velocity level. We conclude that using the acceleration level results in a constant residual velocity for mass 3, and of course a linear drift-off effect and violation of the constraint. The calculation of contact forces is straightforward and follows the reference solution (Fig.~\ref{fig:balls_b_2}). Using the velocity level, the contact force needs to be calculated with an additional integration. This calculation results in a constant position error for this model problem. Using the velocity level results in an oscillatory behavior for the acceleration of mass 3 which means that we expect the same oscillatory behavior for the calculation of contact forces (Fig.~\ref{fig:balls_b_2}). Inserting additional damping for the numerical integration results in better approximations of contact forces and avoids those artificial high frequencies. The results using the Bathe-method or the ED-$\alpha$ method are similar to the results using the generalized-$\alpha$ method with $\rho_{\infty}=0.0$. It is worth to mention that the highest frequency in the system is $f_{\text{max}} = \unit[0.8717]{Hz}$, i.e., $\frac{\Delta t}{T}=0.0436$, which means that the oscillatory behavior of the contact force in case of no damping comes from the structure of the generalized-$\alpha$ method and has nothing to do with the poor representation of high frequencies in the system (model (a)). We conclude that the unsymmetrical structure of the Bathe-method or the ED-$\alpha$ method helps to improve the calculation of contact forces in case of model (b).
A similar observation has already been formulated in~\cite{Sch14a,Sch14b} for the application of half-explicit timestepping schemes on velocity level on examples, which discuss an impacting elastic bar or a rubbing rotor. 

\section{Application to a flexible multibody system}\label{sec:application}
\begin{figure}[ht]
  \centering
  \footnotesize
  \def\svgwidth{0.95\columnwidth}
  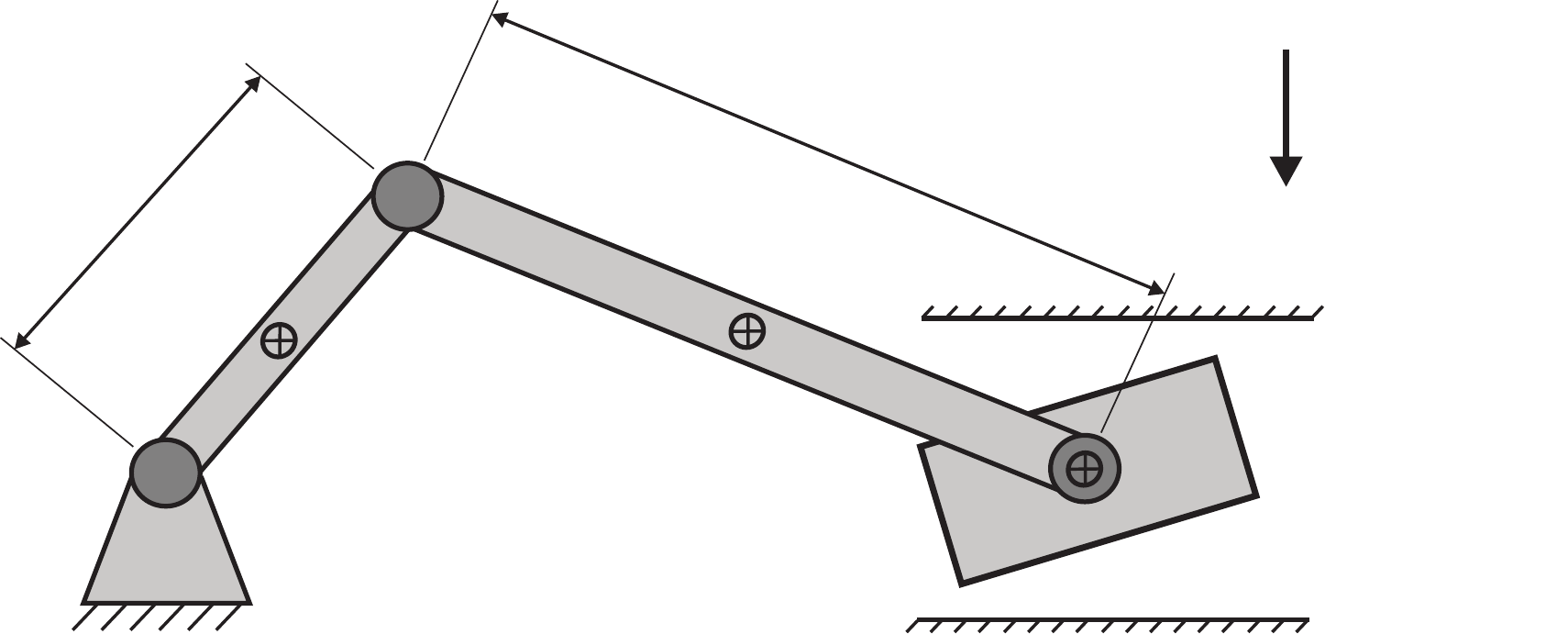
  \caption{Slider-crank mechanism.}
  \label{fig:slider_crank_main}
\end{figure}
We consider the slider-crank mechanism shown in Fig.~\ref{fig:slider_crank_main}, where $l_1$ is the length of the crank and $l_2$ the initial length of the connecting rod for the undeformed state. The inertia of rigid crank, flexible connecting rod and rigid slider consist of translational masses, i.e., $m_1$, $m_2$ and $m_3$, as well as of rotational inertia values, i.e., $J_1$, $J_2$ and $J_3$ for the undeformed state.  For the flexible connecting rod, we consider $\rho$ and $E$ as density and Young's modulus, respectively. The cross-sectional area is given by $A=HD$, i.e., the product of height $H=m_2/(\rho l_2D)$ and thickness $D=\sqrt{12 J_2/m_2-l_2^2}$ of the rod, and the second moment of area is given by $I=\frac{1}{12}HD^3$. The system is subject to gravitation~$\gamma$~\cite{Flo10}.\par 
In order to describe the flexible system, the local coordinate system of the connecting rod is located tangentially in the joint between the crank and the connecting rod. Thus, the basis for a floating frame of reference formulation is accomplished. Such a description is characterized by a separation of the coordinates of an elastic body into reference and elastic coordinates. The reference coordinates delineate the rigid body movement and consist of the translational coordinates describing the absolute position of the local coordinate system and the rotational coordinates describing the orientation by angles. The elastic coordinates capture the flexible movement. Crank and slider are described by minimal coordinates. For the evaluation of the equations of motion, we consider the inertia coupling between the different sets of coordinates. The derivation of the system matrices for the given slider-crank mechanism is studied in~\ref{sec:appD}.

\section{Results}\label{sec:examples}
The application of the generalized-$\alpha$ method, the Bathe-method and the ED-$\alpha$ method as base integration schemes of the overall framework to nonsmooth problems with unilateral contacts with friction is discussed in the following section. First, we show the spatial convergence of the schemes and validate the results comparing to the rigid case. Finally, we discuss the time integration schemes concerning different aspects.
\subsection{Validation for a rigid slider-crank mechanism}
\begin{figure}[ht]
  \centering
  \includegraphics[width=\columnwidth]{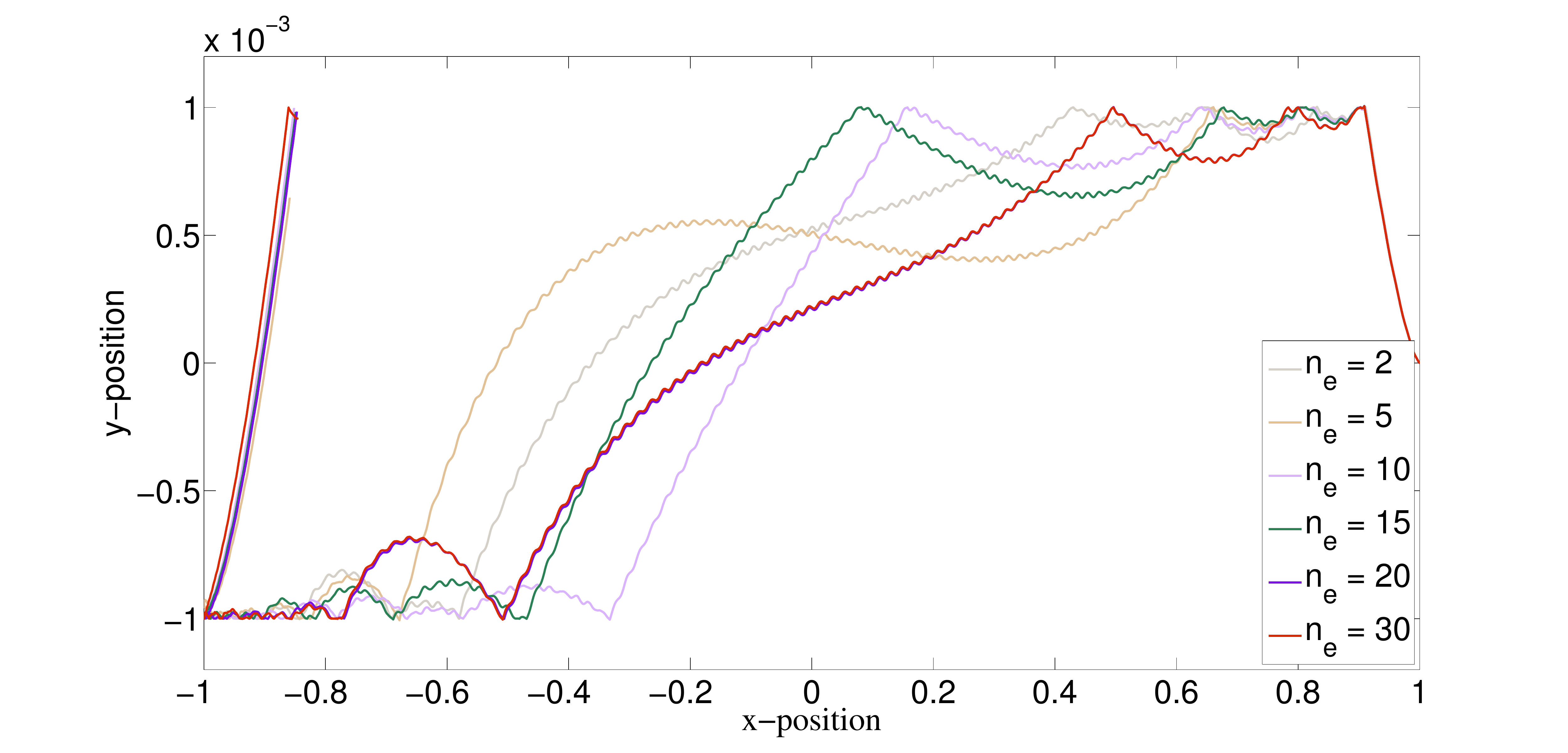}
  \caption{Mass center movement for different numbers of elements.}
  \label{fig:conv_mass_center_el}
\end{figure}

\begin{figure}[ht]
  \centering
  \includegraphics[width=\columnwidth]{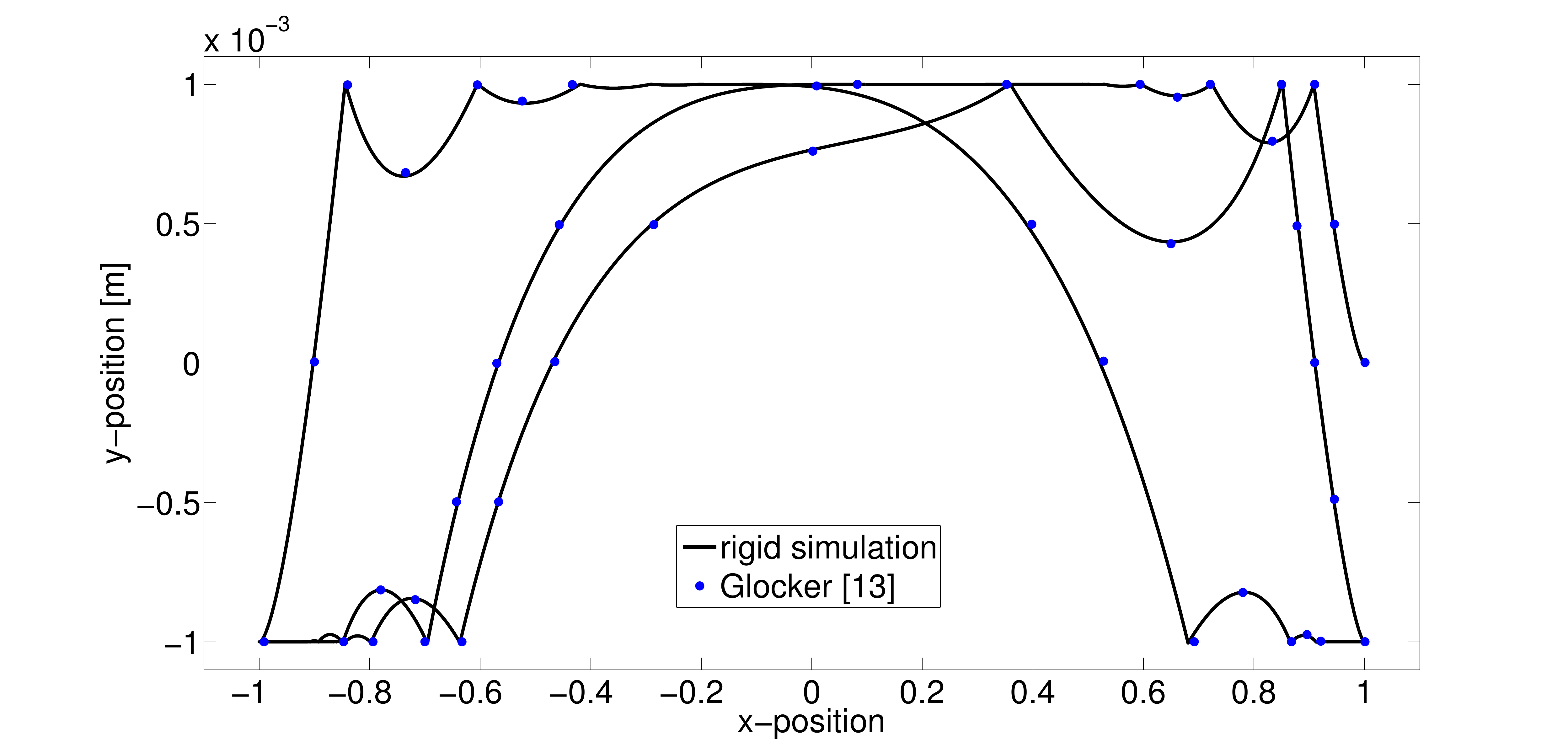}
  \caption{Mass center movement; comparison to a rigid body simulation.}
  \label{fig:val_rigid_uni}
\end{figure}

\begin{figure}[ht]
  \centering
  \includegraphics[width=\columnwidth]{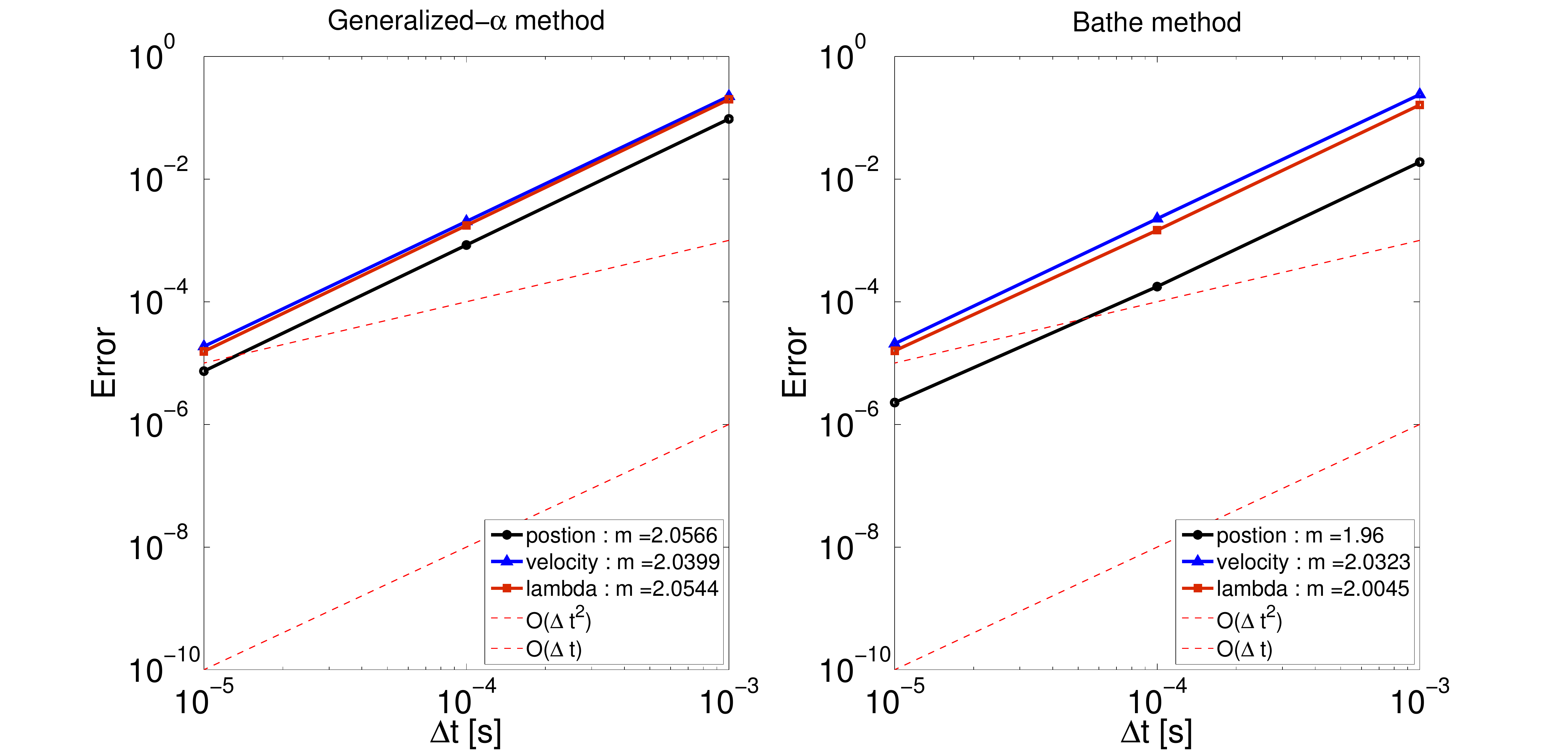}
  \caption{Relative errors of position, velocity and contact force.}
  \label{fig:val_rigid_time}
\end{figure}

Figure~\ref{fig:conv_mass_center_el} shows the convergence for the position of the slider mass center when we increase the number of elements and run the simulation with $\Delta t=\unit[10^{-5}]{s}$ using the Bathe-method. Specific characteristics are given in Table~\ref{tab:slider_crank_characteristics}.
\begin{table}[ht]
  \footnotesize
  \centering
  \begin{tabular}{ll}
    \hline
    Geometrical characteristics & $l_1=\unit[0.1530]{m}$ (length crank)\\
    \ & $l_2=\unit[0.3060]{m}$ (length rod)\\
    \ & $a=\unit[0.0500]{m}$ (half-length slider)\\
    \ & $b=\unit[0.0250]{m}$ (half-height slider)\\
    \ & $c=\unit[0.0010]{m}$ (gap)\\ \hline 
    Inertia properties & $m_1=\unit[0.0380]{kg}$ (mass crank)\\
    \ & $m_2=\unit[0.0380]{kg}$ (mass rod)\\
    \ & $m_3=\unit[0.0760]{kg}$ (mass slider)\\
    \ & $J_1=\unit[7.4\cdot10^{-5}]{kgm^2}$ (inertia crank)\\
    \ & $J_2=\unit[5.9\cdot10^{-4}]{kgm^2}$ (inertia rod)\\
    \ & $J_3=\unit[2.7\cdot10^{-6}]{kgm^2}$ (inertia slider)\\ \hline 
    \ Force elements & $\gamma=\unit[9.81]{m/s^2}$ (gravitation)\\ \hline 
    \ Contact parameters& $\vepsilon_{N_1}=\vepsilon_{N_2}=\vepsilon_{N_3}=\vepsilon_{N_4}=0.4$\\ 
    \ for slider corners & $\vepsilon_{T_1}=\vepsilon_{T_2}=\vepsilon_{T_3}=\vepsilon_{T_4}=0.0$\\ 
    \ & $\mu_1=\mu_2=\mu_3=\mu_4=0.01$\\ \hline 
    \ Initial conditions & $\theta_{1_0}=0.0$\\
    \ & $\theta_{2_0}=0.0$\\
    \ & $\theta_{3_0}=0.0$\\
    \ & $\omega_{1_0}=\unit[150.0]{rad/s}$\\
    \ & $\omega_{2_0}=\unit[-75.0]{rad/s}$\\
    \ & $\omega_{3_0}=\unit[0.0]{rad/s}$\\ \hline
		\ Material properties & $E=\unit[2\cdot10^{11}]{N/m^2}$\\
    \ of flexible rod & $\rho=\unit[7800]{kg/m^3}$\\ \hline 
  \end{tabular}
  \caption{Characteristics of the slider-crank mechanism with unilateral constraints and friction~\cite{Flo10}.}
  \label{tab:slider_crank_characteristics}
\end{table}
 The convergence results are comparable for all base integration schemes.
For the validation of the results, we compare an almost rigid system ($E=\unit[10^{15}]{N/m^2}$) with a rigid multibody system with 3 degrees of freedom. Figure~\ref{fig:val_rigid_uni} shows the results for the position of the slider mass center compared to the same simulation with rigid bodies~\cite{Flo10}.
The generalized-$\alpha$ method and the Bathe-method are both second order accurate. In order to show this behavior in a simulation, we consider the specific case of a bilateral contact by setting the gap $c=\unit[0]{m}$. Figure~\ref{fig:val_rigid_time} shows the convergence when we decrease $\Delta t$. The Slope of the line $m$ shows the accuracy of the calculations in logarithmic scale. The reference solution is calculated with Simpack\footnote{http://www.simpack.com/}.
\subsection{Comparison between generalized-$\alpha$ and Bathe-method}
\begin{figure}[ht]
  \centering
  \footnotesize
  \def\svgwidth{0.6\columnwidth}
  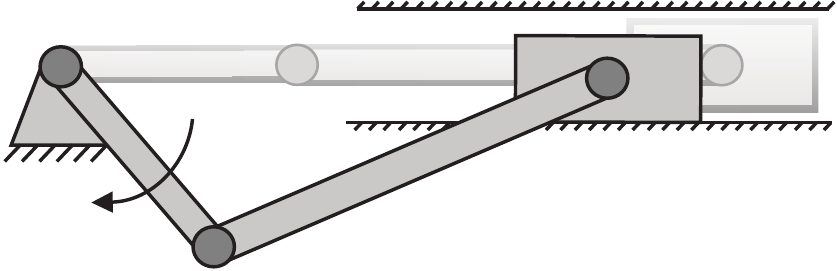
  \caption{Slider crank modified configuration.}
  \label{fig:torque}
\end{figure}

\begin{figure}[ht]
  \centering
  \includegraphics[width=\columnwidth]{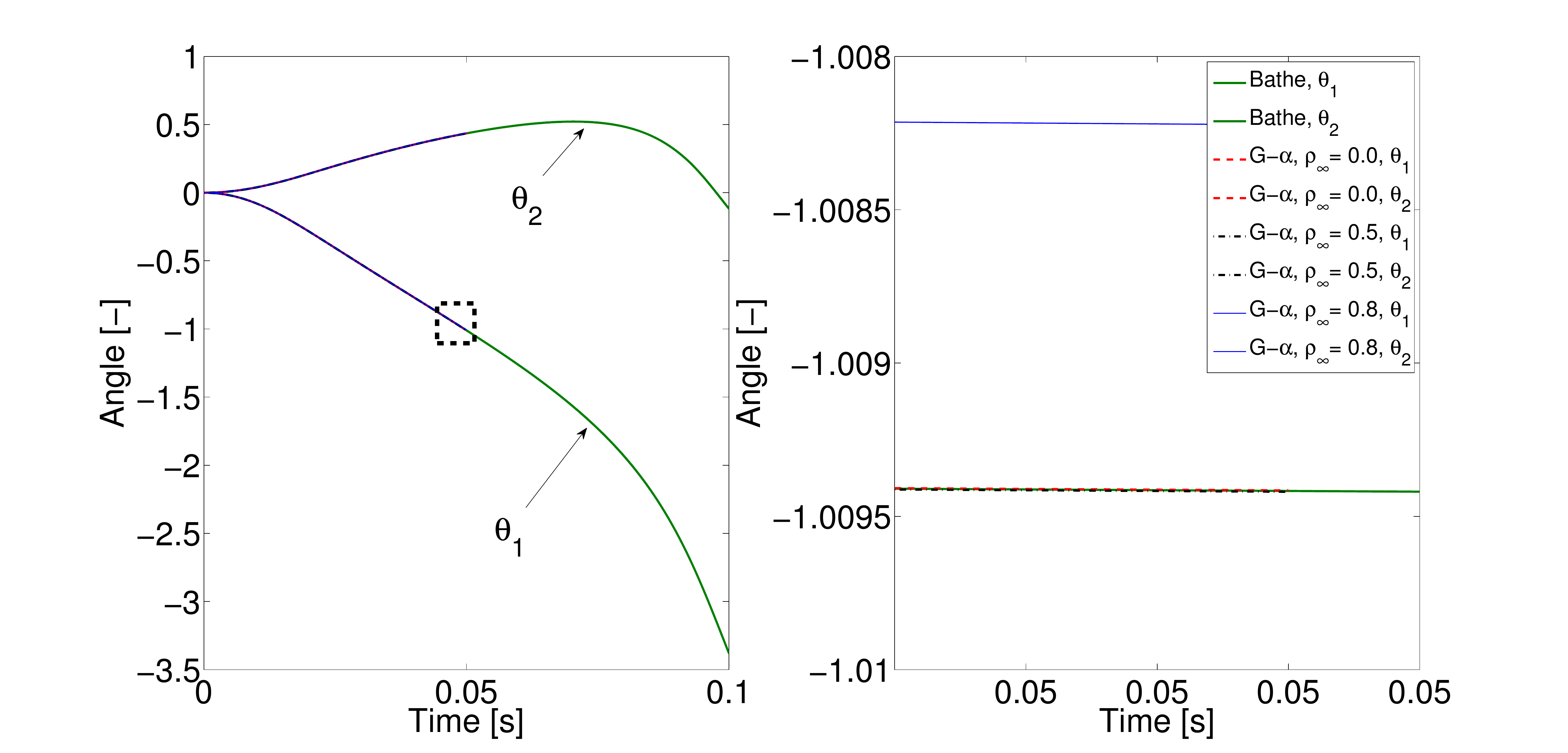}
  \caption{Comparison between generalized-$\alpha$ method and Bathe-method for the angles $\theta_1$ and $\theta_2$.}
  \label{fig:b_vs_g_pos}
\end{figure}

\begin{figure}[ht]
  \centering
  \includegraphics[width=\columnwidth]{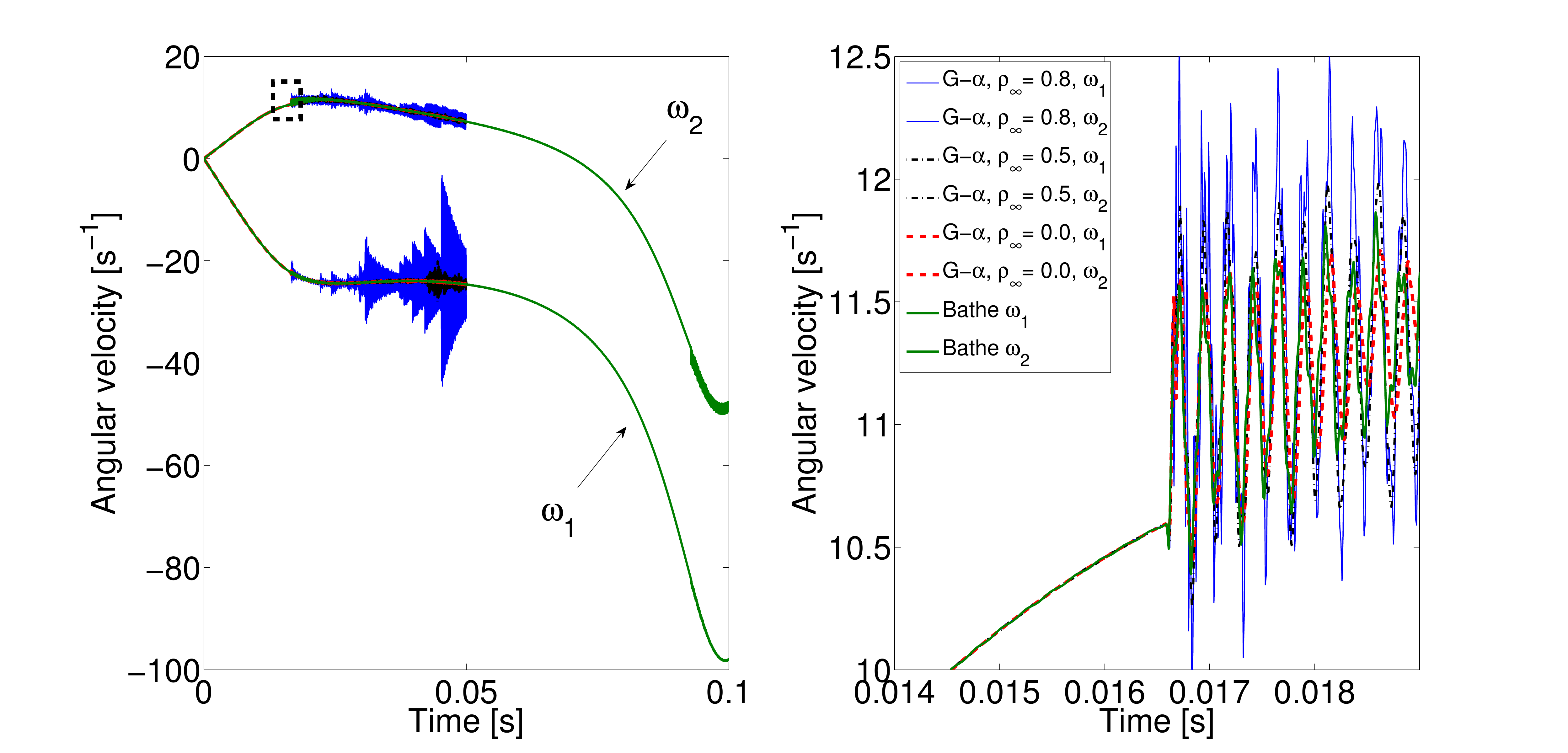}
  \caption{Comparison between generalized-$\alpha$ method and Bathe-method for the angular velocities $\omega_1$ and $\omega_2$.}
  \label{fig:b_vs_g_vel}
\end{figure}

\begin{figure}[ht]
  \centering
  \includegraphics[width=\columnwidth]{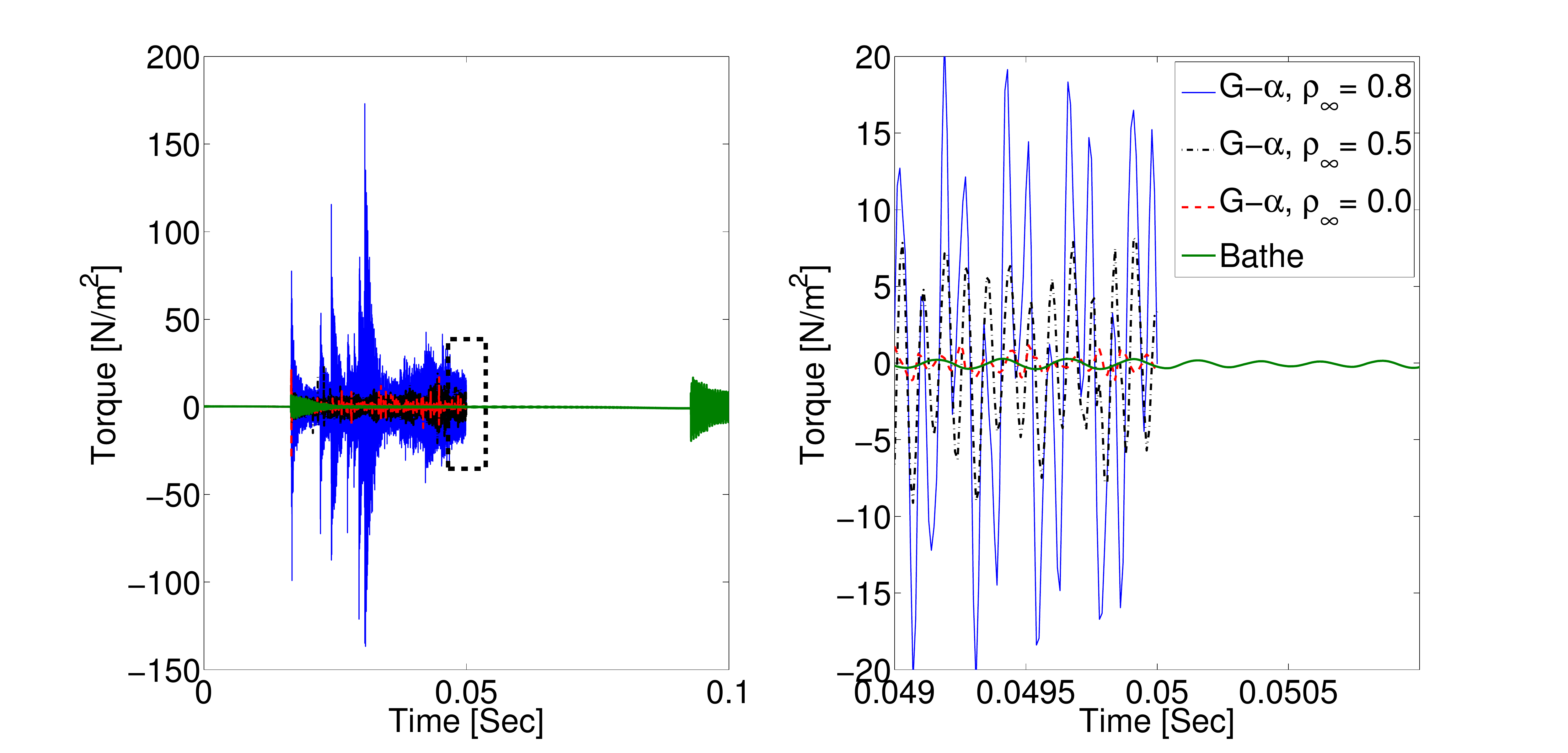}
  \caption{Comparison between generalized-$\alpha$ method and Bathe-method for the torque at the beam root.}
  \label{fig:b_vs_g_mom}
\end{figure}

\begin{figure}[ht]
  \centering
  \includegraphics[width=\columnwidth]{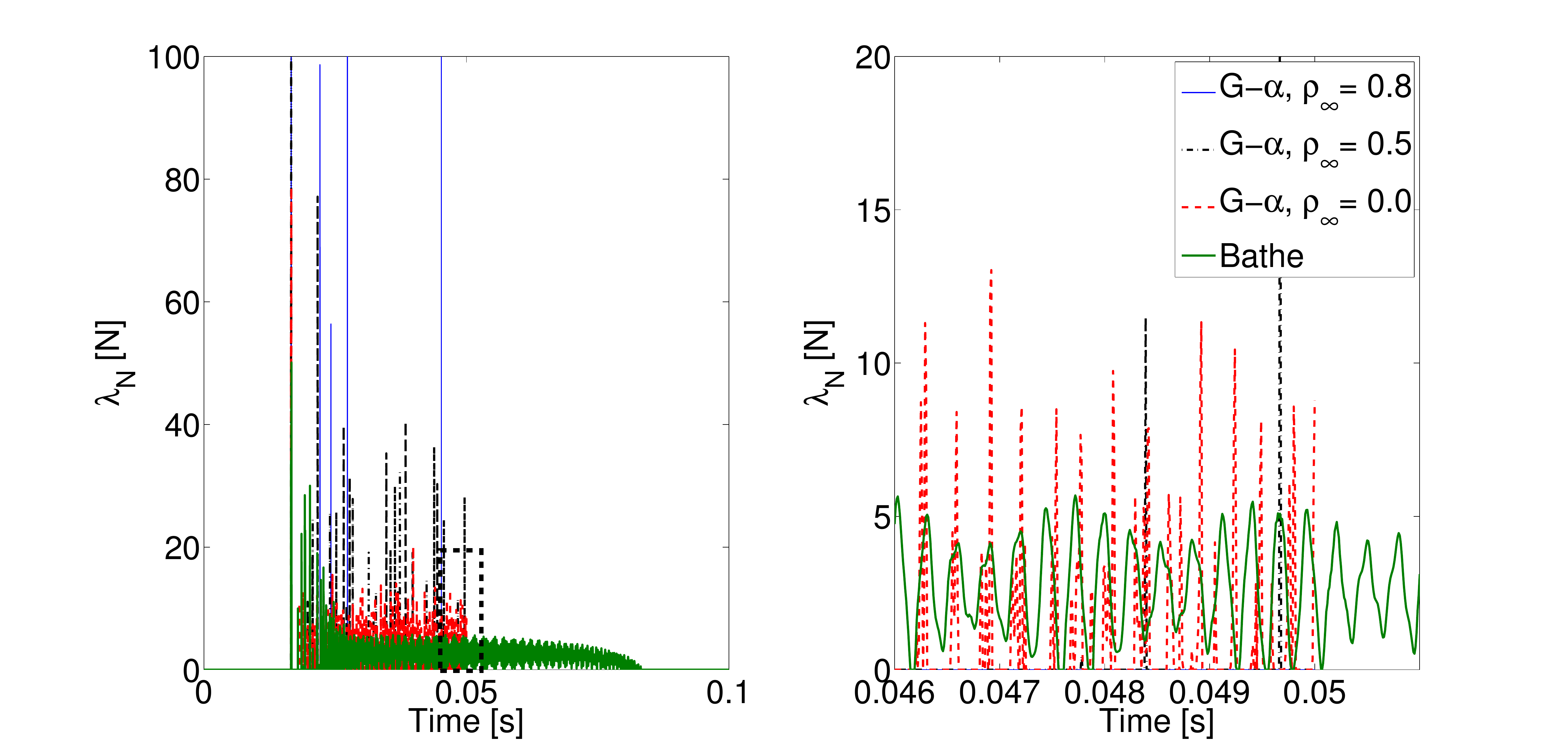}
  \caption{Comparison between generalized-$\alpha$ method and Bathe-method for the active component of the normal contact force.}
  \label{fig:b_vs_g_lambda}
\end{figure}

In order to have a better insight for closed gap situations and the calculation of contact forces, some initial properties of the slider-crank mechanism are changed according to Table~\ref{tab:slider_crank_modified_characteristics}, whereas the other characteristics are set according to Table~\ref{tab:slider_crank_characteristics} (Fig.~\ref{fig:torque}). The time-step size is $\Delta t=\unit[10^{-5}]{s}$. 
\begin{table}[ht]
  \footnotesize
  \centering
  \begin{tabular}{ll}
    \hline
    Geometrical characteristics & $c=\unit[0.0005]{m}$\\ \hline 
    \ Driven torque for crank & $T=\unit[1]{N/m}$\\ \hline 
    \ Contact parameters & $\vepsilon_{N_1}=\vepsilon_{N_2}=\vepsilon_{N_3}=\vepsilon_{N_4}=0.1$\\ 
    \ for slider corners & $\mu_1=\mu_2=\mu_3=\mu_4=0.1$\\ \hline 
    \ Initial conditions & $\omega_{1_0}=\unit[0.0]{rad/s}$\\
    \ & $\omega_{2_0}=\unit[0.0]{rad/s}$\\
    \ & $\omega_{3_0}=\unit[0.0]{rad/s}$\\ \hline 
  \end{tabular}
  \caption{Modified characteristics of the slider-crank mechanism.}
  \label{tab:slider_crank_modified_characteristics}
\end{table}
In Fig.~\ref{fig:b_vs_g_pos}, we see the comparison of the schemes for the angular displacements $\theta_1$ and $\theta_2$, i.e., the inclinations of crank and connecting rod, ($\theta_3=0$). Both algorithms with different damping values behave almost the same, before and after the impact. 
Figures~\ref{fig:b_vs_g_vel} and \ref{fig:b_vs_g_mom} show the comparison for angular velocity and torque at the beam root. 
If we do not consider enough damping for the time integration schemes, high frequency oscillations corrupt the system response, especially in the velocity and stress fields. In this particular case, the results for the generalized-$\alpha$ method with $\rho_{\infty}=0.8, 0.5$ get unstable soon after $t=\unit[0.05]{s}$. This algorithm transfers energy from the higher (artificial) to the lower meaningful modes. Figure~\ref{fig:b_vs_g_lambda} shows the comparison of the contact force using the velocity level approach.
The Bathe-method is able to represent the vibration as the contact is closed, whereas the generalized-$\alpha$ method results in an additional vibration for the velocity which causes the contact condition to go on and off repeatedly in time.
\subsection{Comparison between ED-$\valpha$ and Bathe-method}
\begin{figure}[ht]
  \centering
  \includegraphics[width=\columnwidth]{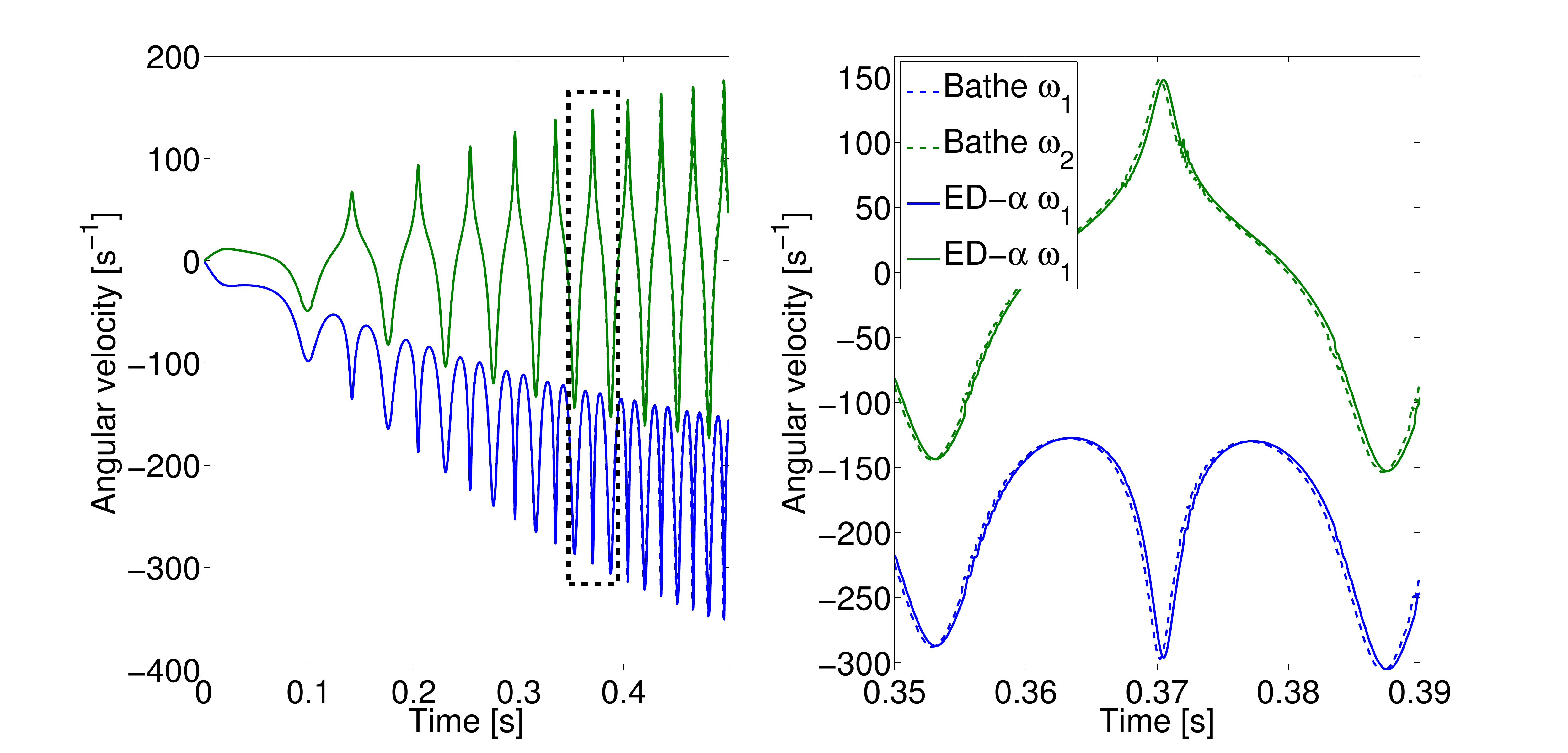}
  \caption{Comparison between ED-$\alpha$ method with $\rho_{\infty}=0$ and Bathe-method for the angular velocities.}
  \label{fig:ed_bathe_compare_vel}
\end{figure}

\begin{figure}[ht]
  \centering
  \includegraphics[width=1\columnwidth]{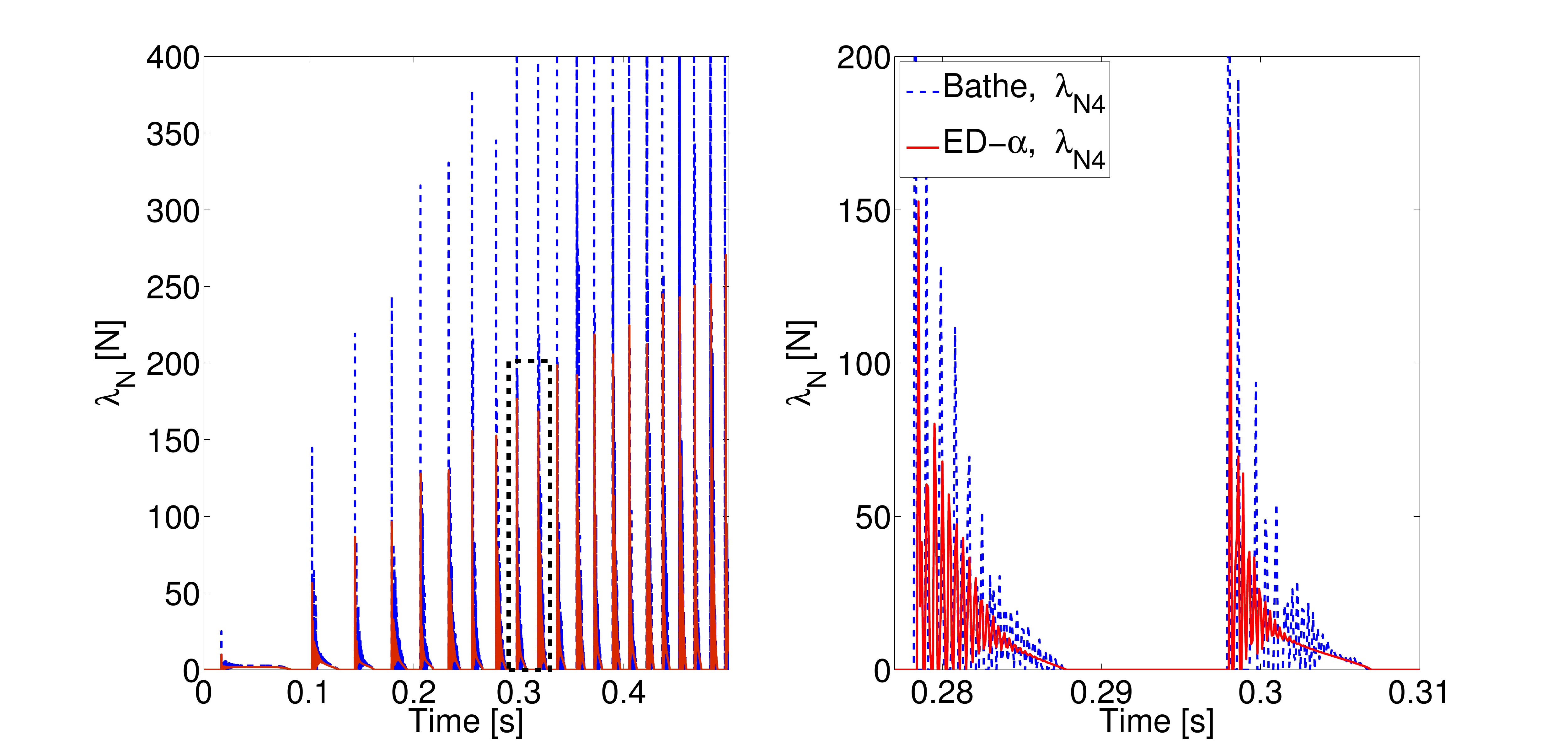}
  \caption{Comparison between ED-$\alpha$ method with $\rho_{\infty}=0$ and Bathe-method for the normal contact force.}
  \label{fig:ed_bathe_compare_lambda}
\end{figure}

Using the data in Table~\ref{tab:slider_crank_modified_characteristics} with time-step size $\Delta t=\unit[10^{-4}]{s}$, we run the same simulations up to $t=\unit[0.5]{s}$ to show the stability and robustness of the ED-$\valpha$ method and the Bathe-method. Figure~\ref{fig:ed_bathe_compare_vel} shows the comparison of angular velocities which match perfectly. 
The phase shift can be explain with the relative period diagram in Fig.~\ref{fig:ed_bathe_compare}. The effect of the phase shift is also noticeable in the calculation of the normal contact force in Fig.~\ref{fig:ed_bathe_compare_lambda}. 
The general behavior and amplitude of the contact force is in good agreement using energy decaying methods.
\clearpage
\subsection{Comparison of computing time}
Based on the setting in Table~\ref{tab:slider_crank_characteristics}, we analyze the relative central processing unit (CPU) time for the computation of a rigid slider-crank mechanism. Thereby, we compare the generalized-$\alpha$ method, the Bathe-method and the ED-$\alpha$ method also with Moreau's midpoint rule, which is a classic timestepping scheme, and the half-explicit timestepping scheme (HETS) proposed in~\cite{Sch14a,Sch14b}.\par
First, we regard the relative CPU time per time-step, exemplary for $\Delta t = \unit[10^{-5}]{s}$ in Table~\ref{tab:time_compare_bi_1} for the bilateral case.
\begin{table}[ht]
  \footnotesize
  \centering
  \begin{tabular}{lccccc}
  \toprule
     & Moreau & HETS & generalized-$\alpha$ & Bathe & ED-$\alpha$ \\
  \midrule
  Rel. CPU time & 1.0 & 1.35 & 4.15 & 4.35 & 5.00\\
  \bottomrule 
  \end{tabular}
  \caption{Relative CPU time for $\Delta t = \unit[10^{-5}]{s}$ (bilateral).}
  \label{tab:time_compare_bi_1}
\end{table}
We compare the necessary time-step sizes $\Delta t$ and their corresponding relative error with respect to the reference Simpack solution in Table~\ref{tab:time_compare_bi_2}. Thereby, we compute the relative error examplary for the connecting rod inclination $\theta_2$ for a set of considered time instances $\left\{t_k\right\}_{k=1}^M$:
\begin{align}
	\text{err} = \norm{\left(\cdots,\frac{\abs{\theta_2(t_k)-\theta_{2_\text{ref}}(t_k)}}{\abs{\theta_{2_\text{ref}}(t_k)}},\cdots\right)}_2 \;.
\end{align}
\begin{table}[ht]
  \footnotesize
  \centering
  \begin{tabular}{lcccc}
  \toprule
	
  $\Delta t$  & $\unit[10^{-4}]{s}$ & $\unit[10^{-5}]{s}$ & $\unit[10^{-6}]{s}$ \\
  \midrule
	Moreau & $1.3\cdot 10^{1}$ & $1.3\cdot 10^{0}$ & $1.3\cdot 10^{-1}$ \\
	HETS & $0.1\cdot 10^{0}$ & $0.1\cdot 10^{-2}$ & $0.1\cdot 10^{-4}$ \\
	generalized-$\alpha$ & $0.1\cdot 10^{0}$ & $0.1\cdot 10^{-2}$ & $0.1\cdot 10^{-4}$ \\
	Bathe & $0.7\cdot 10^{-1}$ & $0.7\cdot 10^{-3}$ & $1.6\cdot 10^{-5}$ \\
	ED-$\alpha$ & $1.3\cdot 10^{-1}$ & $1.3\cdot 10^{-3}$ & $1.6\cdot 10^{-5}$ \\	
  \bottomrule 
  \end{tabular}
  \caption{Comparison of the error $\epsilon$ for different $\Delta t$ ($T=\unit[0.05]{s}$, bilateral).}
  \label{tab:time_compare_bi_2}
\end{table}
For a given time-step size, the computational effort of the new methods is minimal using the generalized-$\alpha$ method as there is no additional midpoint calculation. For the Bathe-method and the ED-$\alpha$ method, we use one additional point in the time integration algorithm which explains the increase in the computing time. In case of the ED-$\alpha$ method, the additional point can be interpreted as a jump at the beginning of the interval and according to \eqref{eq:fixed_point_ed}, we increase the unknowns which have to be solved simultaneously by a factor of 2. The Bathe-method selects the additional point in the middle of the interval. The unknowns are solved independently from the unknowns at $t_{i+1}$. Smaller relative time-steps ($\frac{\Delta t}{2}$) result in smaller changes of the mass matrix, stiffness matrix and force vector and thus in less calculation time for the unknowns compared to the ED-$\alpha$ method as we need less iterations to find the unknowns at the middle of the interval compared to the end of the interval. The old methods are much faster per time-step because of low-order or explicit evaluations. However for a complete comparison, we have to consider also Table~\ref{tab:time_compare_bi_2}. For the same error, e.g. $10^{-1}$, the time-step size for the classic timestepping has to be chosen much smaller than for the remaining four schemes. The benefit of the new schemes in comparison to the half-explicit timestepping, which is also of second order, is the high-frequency damping.\par
For the unilateral case, the solution of the ED-$\alpha$ method with $\rho_{\infty}=0.0$, $\Delta t=\unit[10^{-7}]{s}$ is chosen as the reference solution. Results are given in Tables~\ref{tab:time_compare_uni_1} and \ref{tab:time_compare_uni_2}.
\begin{table}[ht]
  \footnotesize
  \centering
  \begin{tabular}{lccccc}
  \toprule
     & Moreau & HETS & generalized-$\alpha$ & Bathe & ED-$\alpha$ \\
  \midrule
  Rel. CPU time  &  1.0  &  1.15  &  2.72 &  3.15 & 3.40 \\
  \bottomrule 
  \end{tabular}
  \caption{Relative CPU time for $\Delta t = \unit[10^{-5}]{s}$ (unilateral).}
  \label{tab:time_compare_uni_1}
\end{table}
\begin{table}[ht]
  \footnotesize
  \centering
  \begin{tabular}{lcccc}
  \toprule
  $\Delta t$ & $\unit[10^{-4}]{s}$ & $\unit[10^{-5}]{s}$ & $\unit[10^{-6}]{s}$ \\
  \midrule
	Moreau & $6.7\cdot 10^{-2}$ & $6.4\cdot 10^{-3}$ & $6.4\cdot 10^{-4}$ \\
	HETS & $1.1\cdot 10^{-3}$ & $2.8\cdot 10^{-5}$ & $2.3\cdot 10^{-6}$ \\
	generalized-$\alpha$ & $2.0\cdot 10^{-3}$ & $1.2\cdot 10^{-4}$ & $1.4\cdot 10^{-5}$ \\
	Bathe & $1.4\cdot 10^{-3}$ & $7.8\cdot 10^{-5}$ & $8.9\cdot 10^{-6}$ \\
	ED-$\alpha$ & $9.4\cdot 10^{-4}$ & $2.8\cdot 10^{-5}$ & $2.3\cdot 10^{-6}$ \\	
  \bottomrule 
  \end{tabular}
  \caption{Comparison of error $\epsilon$ for different $\Delta t$ ($T=\unit[0.05]{s}$, unilateral).}
  \label{tab:time_compare_uni_2}
\end{table}
They confirm the results of the bilateral case and are even advantageous for the new methods because of their increased stability.
\subsection{Modal approach} \label{sec:modal}
A major advantage of the floating frame of reference formulation is that the finite element nodal coordinates can be easily reduced using modal analysis techniques, based on a reduced set of eigenvectors (\ref{sec:appD}).\par
\begin{figure}[ht]
  \centering
  \includegraphics[width=1.0\columnwidth]{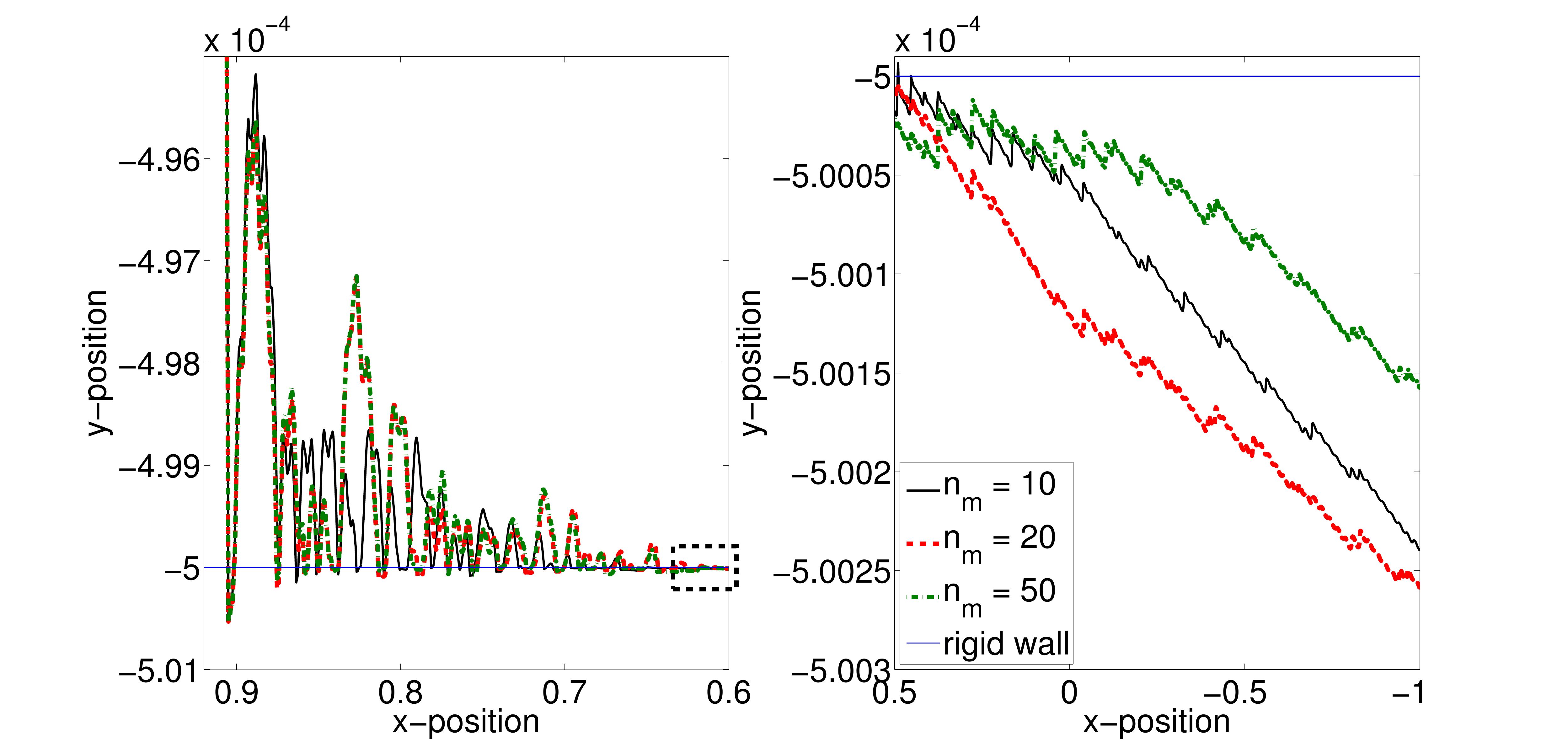}
  \caption{Mass center portrait for modal approach and tangential boundary conditions.}
  \label{fig:clamped_modal_conv_mass_center}
\end{figure}
\begin{figure}[htbp]
  \centering
  \includegraphics[width=1.0\columnwidth]{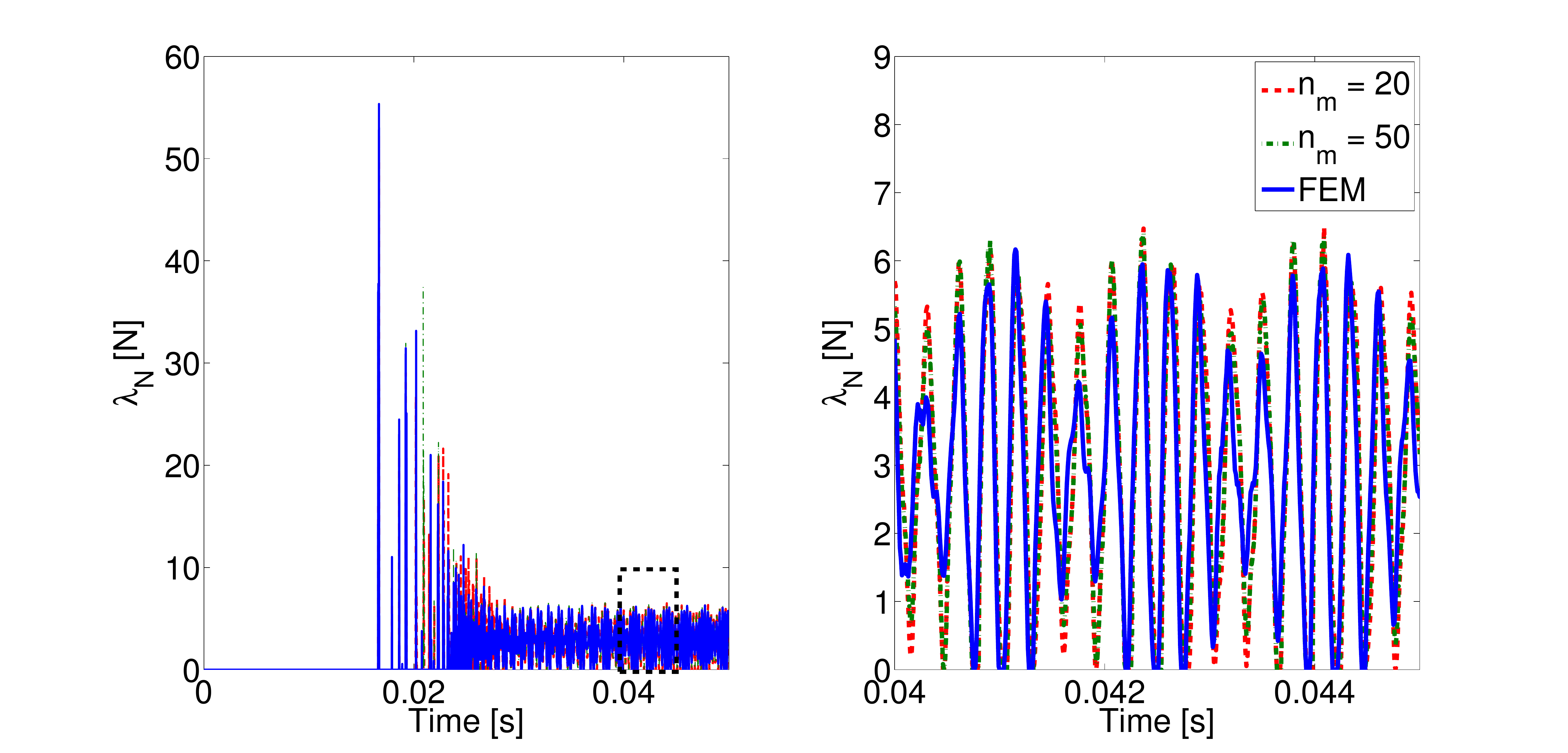}
  \caption{Normal contact force convergence for modal approach.}
  \label{fig:clamped_modal_conv_lambda}
\end{figure}

The spatial convergence is plotted in Fig.~\ref{fig:clamped_modal_conv_mass_center} and Fig.~\ref{fig:clamped_modal_conv_lambda} when we increase the numbers of mode shapes. Thereby, the simulation is based on Table~\ref{tab:slider_crank_modified_characteristics}, the time-step size is $\Delta t=\unit[10^{-5}]{s}$, the boundary condition is considered as clamped-free (tangential), and we use the Bathe-method for time integration. After an impact occurred, depending on the boundary condition and impact points, we need more number of modes to describe the vibration behavior than in the non-impulsive period (before the impact). As far as we deal with nonlinear problems, slightly different conditions before the impact may result in large changes in time after the impact, which explains the differences in the results. Considering spectral analysis for the Bathe-method (Fig.~\ref{fig:bathe_newmark_compare}) after $\frac{f_c}{f_s}=\frac{\Delta t}{T}=10^2$, there will be no effective mode in the response. We simply calculate a sufficient number of modes considering the high frequency dissipation of the Bathe-method. In this case, all frequencies up to $f_c=\frac{10^2}{\unit[10^{-5}]{s}}=\unit[10^7]{1/s}$ (which covers the first 52 modes out of total 63 modes) can be effective in the final solution. Results are in good agreement between modal and full FEM simulation for the calculation of normal contact forces.\par
\begin{figure}[ht]
  \centering
  \includegraphics[width=1.0\columnwidth]{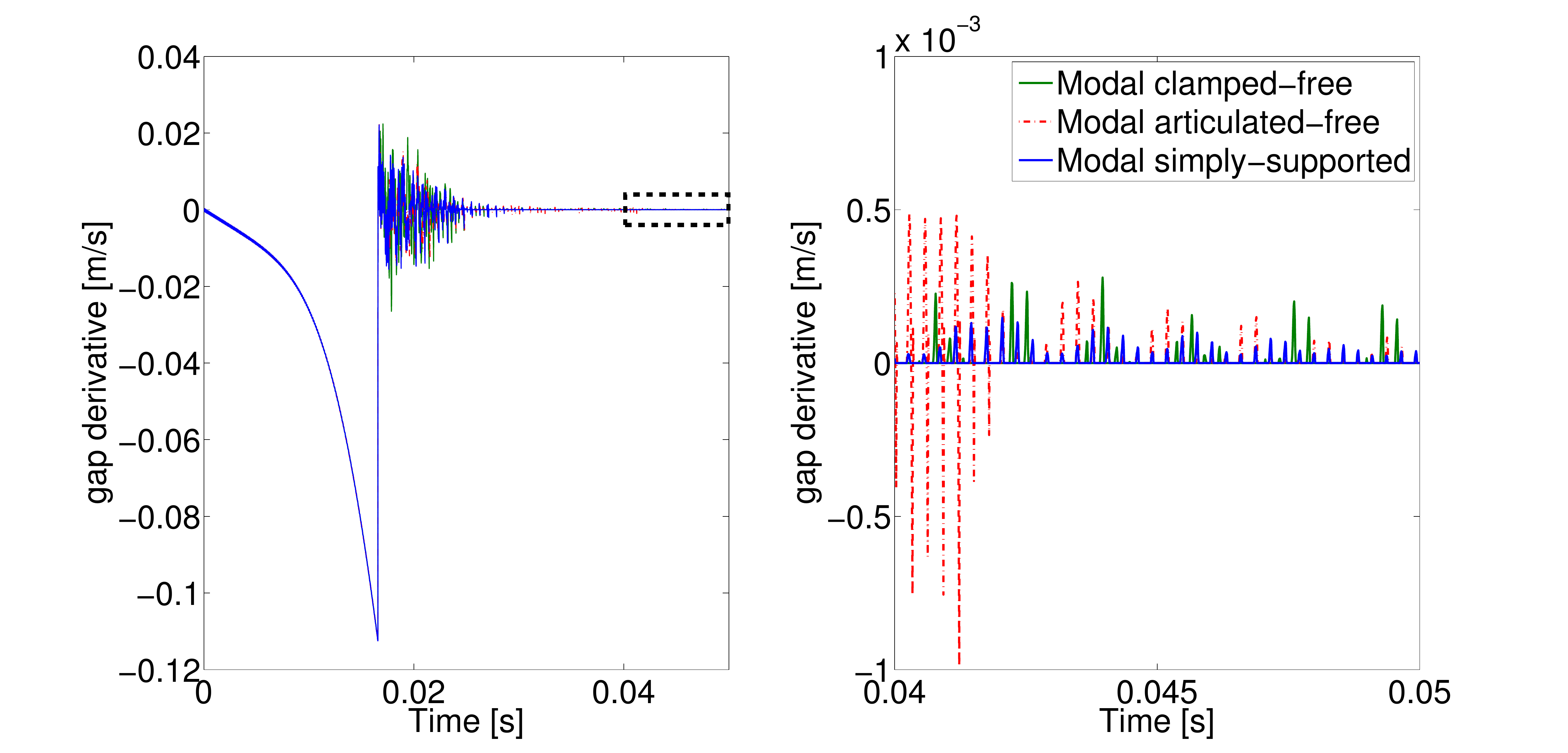}
  \caption{Comparison between FEM and modal solution for different boundary conditions.}
  \label{fig:fem_modal}
\end{figure}

Figure~\ref{fig:fem_modal} shows the comparison between modal solutions for different boundary conditions. Different boundary conditions result in different mode shapes and frequencies. Assuming the same impact point and condition we get different vibration behavior, but the general pattern may be similar.

\section{Summary and Conclusion\label{sec:summary}}
This work deals with the consistent and efficient integration of nonsmooth flexible multibody systems with impacts and dry friction. We present a timestepping scheme, which evolves from the idea of time-discontinuous Galerkin methods~\cite{Sch14a,Sch14b}. However, we abstract this origin and develop a framework which improves a non-impulsive trajectory by impulsive correction after each time-step if necessary. This correction is automatic and is evaluated on the same kinematic level as the piecewise non-impulsive trajectory, i.e., on velocity level. The resulting overall mixed timestepping scheme is consistent for impulses and benefits from higher order in non-impulsive periods and all advantages of the base integration schemes used to calculate the approximation per time-step.\par
We present a nonsmooth adaptation of the generalized-$\alpha$ method, the Bathe method and the ED-$\alpha$ method. It is applied to a slider-crank mechanism with a flexible connecting rod, impacts and dry friction. The elastic behavior of the connecting rod is compared using the different base integration schemes and a modal approach. The results are validated with respect to~\cite{Flo10} and a rigid body simulation in Simpack. It is shown, that introducing enough damping for Newmark-type integrators like the generalized-$\alpha$ method results in more stable solutions for nonsmooth multibody systems especially on velocity and acceleration level in comparison to~\cite{Sch14a,Sch14b}. It seems that schemes like the Bathe-method and the ED-$\alpha$ method behave more robust in the nonlinear regime. These methods are more expensive per time-step but less steps can be used and the methods remain stable, even if the Newmark-type integrator fails for large deformations and long time duration dynamic response calculations.\par
For the future, it is valuable to pursue deeper mathematical analysis of the overall framework with different base integration schemes to prove the characteristics observed numerically especially for nonlinear multibody simulation.

\appendix
\section{Generalized-$\valpha$ method} \label{sec:appA}
For the analysis of the generalized-$\alpha$ method, it is advantageous to reduce the coupled equation of motion to a series of uncoupled single degree of freedom systems using modal analysis and using eigenvector orthogonality. The linear single degree of freedom system with angular frequency $\omega$ is given by
\begin{align}
  \ddot{q} + \omega^2 q  =  0 \;,\label{eq:SDOF}
\end{align}
where the terms related to external damping and forces are set to zero to study accuracy and stability properties of the algorithm. The generalized-$\alpha$ method described in \eqref{eq:mechanical_system} to \eqref{eq:mechanical_system_initial_acceleration} can be written in the compact form
\begin{align}
  \vX_{i+1}= \vA_{g\alpha}\vX_i\;,\quad i\in\left\{ 0,~1,...,~N-1 \right\} \;, \label{eq:X_A}
\end{align}
where $\vX_i=\left(q_i,~\Delta t v_i,~\Delta t^2 a_i \right)^T$ and $\vA_{g\alpha}$ is the amplification matrix for the generalized-$\alpha$ method. With   
\begin{align}
  \begin{aligned}
    D&=1-\alpha_m+\left(1-\alpha_f\right)\beta\Omega^2 \;,\\
    \Omega&=\omega \Delta t  \;,\\  
    \omega&=\sqrt{k/m}\;,
  \end{aligned}
\end{align}
it is
\begin{align}
  \vA_{g\alpha}&= \frac{1}{D}\begin{pmatrix} 
    1-\alpha_m-\alpha_f\beta\Omega^2 & 1-\alpha_m & \left(\frac{1}{2}-\beta\right)\left(1-\alpha_m\right)-\beta\alpha_m \\ 
    -\gamma\Omega^2 & 1-\alpha_m-\left(1-\alpha_f\right)\left(\gamma-\beta\right)\Omega^2 & \left(1-\gamma\right)\left(1-\alpha_m\right)-\gamma\alpha_m-\left(1-\alpha_f\right)\left(\frac{\gamma}{2}-\beta\right)\Omega^2 \\ 
    -\Omega^2 & -\left(1-\alpha_f\right)\Omega^2 & -\left(1-\alpha_f\right)\left(\frac{1}{2}-\beta\right)\Omega^2-\alpha_m \end{pmatrix}\;.\label{eq:A}
\end{align}
The accuracy of an algorithm can be determined using the difference equation in terms of the displacement
\begin{align}
  q_{i+1} - A_1 q_i + A_2 q_{i-1} - A_3 q_{i-2}=0 \;,\label{eq:difference_disp}
\end{align}
where $A_1$ is the trace of $\vA_{g\alpha}$, $A_2$ is the sum of the principal minors of $\vA_{g\alpha}$ and $A_3$ is the determinant of $\vA_{g\alpha}$. It can be shown~\cite{Chu93} that the algorithm is second order accurate for unconstrained mechanical systems if
\begin{align}
  \gamma=\frac{1}{2}-\alpha_m+\alpha_f \;. \label{eq:gamma}
\end{align}
The parameter $\gamma$ is responsible for numerical dissipation. For $\gamma=0.5$, there is no numerical dissipation, whereas for values $\gamma>0.5$ the numerical dissipation increases. The stability and numerical behavior of an algorithm depends on the eigenvalues of the amplification matrix. The spectral radius $\rho$ of an algorithm is defined by
\begin{align}
  \rho = \max{\left(|\lambda_1|,~|\lambda_2|,~|\lambda_3|\right)}\;,
\end{align}
where $\lambda_i$ is the $i$th eigenvalue of $\vA_{g\alpha}$. An algorithm is unconditionally stable for linear problems if $\rho\leq1$ for all $\Omega\in\left[0,\infty\right)$.
The generalized-$\alpha$ method is unconditionally stable provided
\begin{align}
  \alpha_m\leq\alpha_f\leq \frac{1}{2}\;, \quad \beta \geq \frac{1}{4}+\frac{1}{2}\left(\alpha_f-\alpha_m\right) \;.\label{eq:stability_condition}
\end{align}
The spectral radius is a measure for numerical dissipation. A smaller spectral radius corresponds to greater numerical dissipation. Desirable dissipation properties have spectral radius close to unity in the low-frequency domain; the value smoothly decreases as $\Omega$ increases. Typically, in the low-frequency domain, $|\lambda_3|\le|\lambda_{1,2}|$. To preserve the smoothness as $\Omega$ increases, $|\lambda_3|\le|\lambda_{1,2}|$ for all $\Omega\in\left[0,\infty\right)$. Violation of this condition will result in a \emph{cusp} in the spectral radius plot where $\rho$ increases as $\Omega$ increases (point A in Fig.~\ref{fig:alpha_space}). Calculating $\lim_{\Omega\rightarrow\infty}\vA_{g\alpha}$ in \eqref{eq:A}, we get the following eigenvalues of the amplification matrix in the high-frequency domain:
\begin{align}
  \begin{aligned}
    \lambda^{\infty}_{1,2}&=\frac{1}{4\beta}\left(4\beta-\left(2\gamma+1\right)\pm j\sqrt{16\beta-\left(2\gamma+1\right)^2}\right)\;,\\
    \lambda^{\infty}_{3}&=\frac{\alpha_f}{\alpha_f-1} \;,
  \end{aligned}\label{A_eigenvalues}
\end{align}
where $j=\sqrt{-1}$. High-frequency dissipation is maximized if the principal roots ($\lambda^{\infty}_{1,2}$) become real, i.e., $\Im(\lambda^{\infty}_{1,2})=0$. It can be shown from \eqref{A_eigenvalues} that for the generalized-$\alpha$ method, this condition is satisfied if
\begin{align}
  \beta = \frac{1}{4}\left(1-\alpha_m+\alpha_f\right)^2 \;,\label{eq:beta}
\end{align}
which satisfies the second condition in \eqref{eq:stability_condition}. The generalized-$\alpha$ method can be described in terms of the two remaining free parameters $\alpha_m$ and $\alpha_f$. Using \eqref{eq:gamma} and \eqref{eq:beta}, we rewrite $\lambda^{\infty}_{1,2}$ in \eqref{A_eigenvalues} as
\begin{align}
  \lambda^{\infty}_{1,2} = \frac{\alpha_f-\alpha_m-1}{\alpha_f-\alpha_m+1}\;.
\end{align}
Let $\rho_{\infty}$ denote the user-specified value of the spectral radius in the high-frequency limit. Since we require that $\lambda_{3} \leq \lambda_{1,2}$ for all $\Omega$, it is $\rho=|\lambda^{\infty}_{1,2}|$. It has been shown that for a given level of high-frequency dissipation, i.e., for fixed $\rho_{\infty}$, low-frequency dissipation is minimized if $\lambda^{\infty}_{1,2}=\lambda^{\infty}_{3}$, i.e., if $\alpha_f=\left(\alpha_m+1\right)/3$ (dotted line in Fig.~\ref{fig:alpha_space}). It is more convenient to describe this optimal case by defining $\alpha_m$ and $\alpha_f$ in terms of $\rho_{\infty}$:
\begin{align}
  \alpha_m=\frac{2\rho_{\infty}-1}{\rho_{\infty}+1}\;,\quad\alpha_f=\frac{\rho_{\infty}}{\rho_{\infty}+1} \;.\label{eq:alpha_rho} 
\end{align}
It is worth to mention that for the exact solution of \eqref{eq:SDOF}, we get:
\begin{align}
  \rho_{\text{ex}}=1\;, \quad T_{\text{ex}}=\frac{1}{\omega\Delta t}\;.
\end{align}

\section{Bathe-method} \label{sec:appB}
Applying the Bathe-method (Fig.~\ref{fig:bathe_vector}) to the single degree of freedom system introduced in \eqref{eq:SDOF}, we can write the amplification matrix~\cite{Bat12}:
\begin{figure}[ht]
  \centering
  \footnotesize
  \def\svgwidth{0.6\columnwidth}
  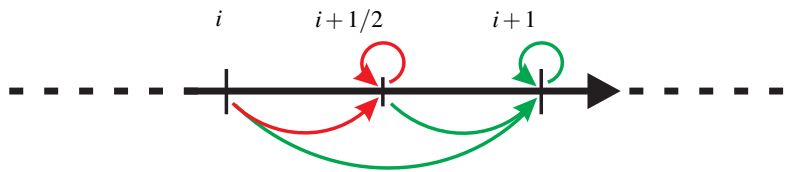
  \caption{Bathe-method in time; red: trapezoidal rule, green: Euler backward rule.}
  \label{fig:bathe_vector}
\end{figure}
\begin{align}
  \vA_{\text{Bathe}}=\frac{1}{D}\begin{pmatrix} 
  -144\Omega^2+19\Omega^4& -144\Omega^2+5\Omega^4 & -28\Omega^2 \\ 
  -96\Omega^2+\Omega^4 & 144-47\Omega^2 & 48-4\Omega^2 \\ 
  144-19\Omega^2 & 144-5\Omega^2 & 28 \end{pmatrix} \;,\label{eq:A_bathe}
\end{align}
where
\begin{align}
  D=\left(16+\Omega^2\right)\left(9+\Omega^2\right) \;.
\end{align}
Analyzing the eigenvalues of $\vA_{\text{Bathe}}$ shows that the method is unconditionally stable and $\rho_{\infty}=0$.

\section{Floating frame of reference description for a slider-crank mechanism}\label{sec:appD}
\begin{figure}[ht]
  \centering
  \footnotesize
  \def\svgwidth{0.95\columnwidth}
  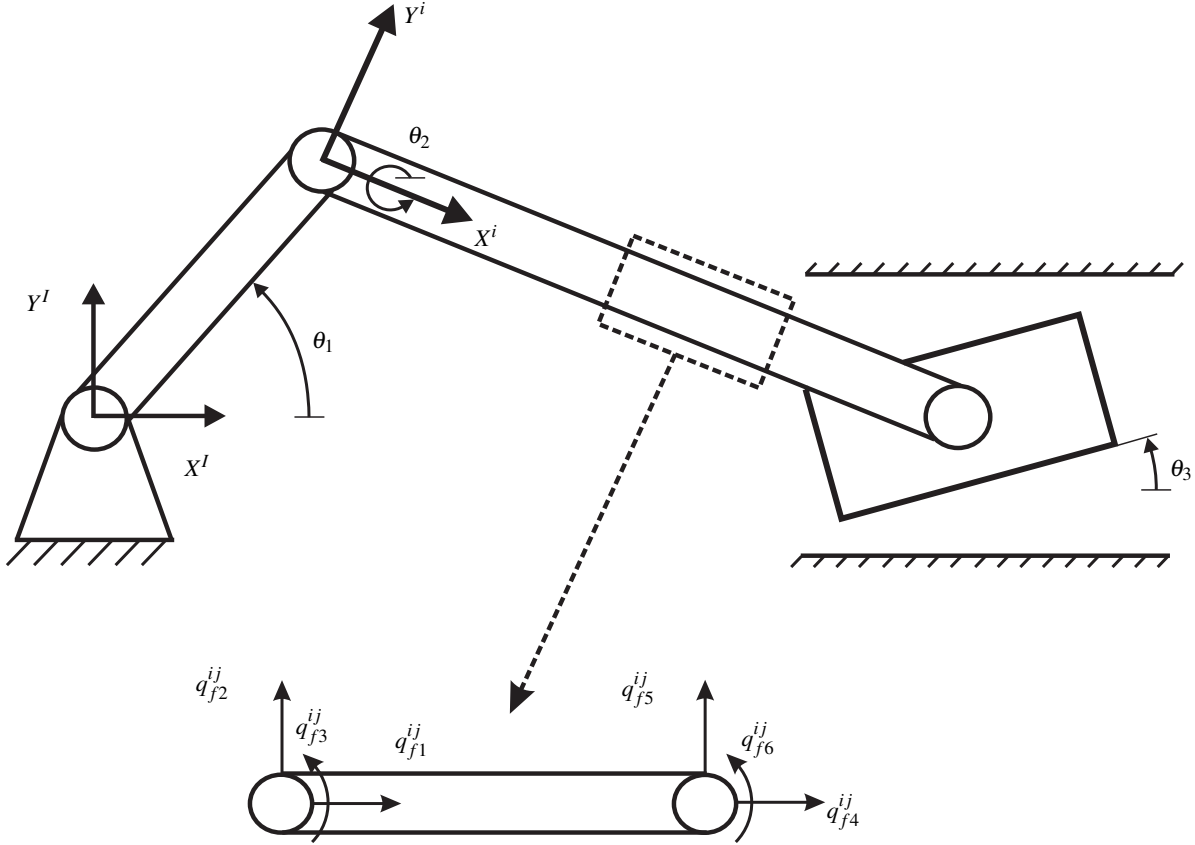
  \caption{The slider-crank mechanism with flexible rod and unilateral constraints, element $j$ of body $i$.}
  \label{fig:slider_floating_frame}
\end{figure}

The following section presents the derivation of the system matrices for the given slider-crank mechanism in Fig.~\ref{fig:slider_floating_frame}. 
\subsection{Finite element formulation}
For the present problem, we consider the two-dimensional beam element shown in Fig.~\ref{fig:slider_floating_frame}. We describe the displacement field within the element by the following polynomials in two directions $X^i$ and $Y^i$:
\begin{align}
  &w_{x^i}=a_0 + a_1 x_1 \;,\\
  &w_{y^i}=a_2 + a_3 x_1 + a_4 x^2_1 + a_5 x^3_1 \;.
\end{align}
Using the above relation for the displacement and considering nodal coordinates (Fig.~\ref{fig:slider_floating_frame}) for each element, we obtain space-independent shape functions of the beam element:
\begin{align}
  \vS^{ij}=\begin{pmatrix} 1-\xi & 0 & 0 & \xi & 0 & 0\\
0 & 1-3\xi^2+2\xi^3 & l\left(\xi-2\xi^2+\xi^3\right) & 0 & 3\xi^2-2\xi^3 & l\left(\xi^3-\xi^2\right) 
    \end{pmatrix} \;,\label{eq:shape_functions}
\end{align}
where $\xi = x/l$.
\subsection{Floating frame of reference}\label{sec:FFR}
\begin{figure}[ht]
  \centering
  \footnotesize
  \def\svgwidth{0.95\columnwidth}
  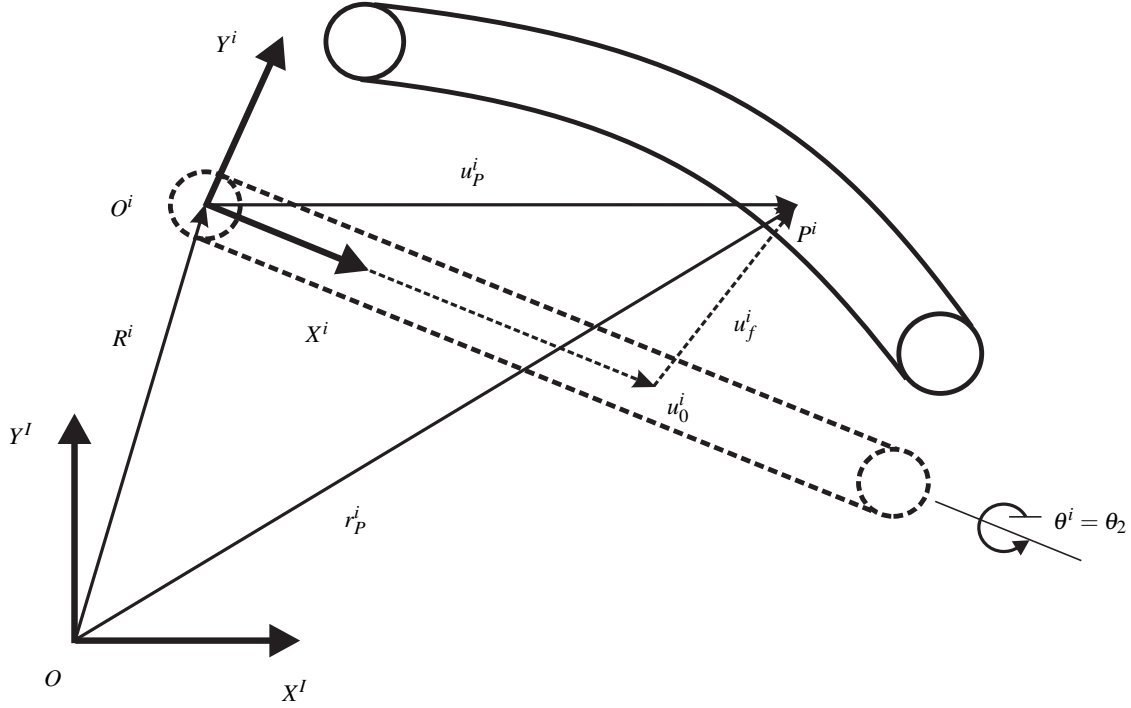
  \caption{Deformable body coordinate.}
  \label{fig:disp_vector}
\end{figure}

In the floating frame of reference formulation presented in this section, the configuration of each deformable body in the multibody system is identified by using two sets of coordinates: reference and elastic coordinates~\cite{Sha05}. Reference coordinates define the location and orientation of a selected body reference ($X^i$, $Y^i$ and $\theta^i$ in Fig.~\ref{fig:disp_vector}). Elastic coordinates describe the body deformation with respect to the body reference ($q^{ij}_1$ to $q^{ij}_6$ in Fig.~\ref{fig:slider_floating_frame}).
The motion of the body is defined as the motion of its reference plus the motion of the material points on the body with respect to its reference (Fig.~\ref{fig:disp_vector}). We write: 
\begin{align}
  \vu^{ij}_f=\vS^{ij}\vq^{ij}_f \;,\label{eq:u_f}
\end{align}
where $\vu^{ij}_f = \begin{pmatrix}u_{f1} & u_{f2} \end{pmatrix}^T$ is the deformation vector of element $j$ of deformable body $i$, $\vS^{ij}$ is the shape matrix of element $j$, $\vq^{ij}_f$ is the vector of elastic coordinates that contains the time dependent nodal values $q^{ij}_1$ to $q^{ij}_6$. For an arbitrary body $i$ of the system, e.g. the flexible rod in this example, we select a body reference $\{X^i,Y^i\}$, the location and orientation of which with respect to the global coordinate system are defined by a set of coordinates called reference coordinates and denoted as $\vq^i_r$. For the planar motion of deformable bodies, which is a special case of three-dimensional motion, the vector $\vq^i_r$ can be written in a partitioned form as
\begin{align}
  \vq^i_r = \begin{pmatrix}\vR^i & \theta^i \end{pmatrix}^T\;,
\end{align}
where $\vR^i$ is a set of Cartesian coordinates that define the location of the origin of the body reference (Fig.~\ref{fig:disp_vector}) and $\theta^i$ is a set of rotational coordinates that describe the orientation of the selected body reference (in the present planar case, it is a scalar value). There is no rigid body motion between the body and its coordinate system. The floating frame of reference formulation does not lead to a separation between the rigid body motion and the elastic deformation.\par
In the case of a rigid body, the global position of an arbitrary point $P$ on the rigid body can be written in planar analysis as:
\begin{align}
  \vr^i_P=\vR^i+\vA^i\vu^i_P\;,
\end{align}
where $\vu^i_P$ is the local position of point $P$ and $\vA^i$ is the transformation matrix defined as
\begin{align}
  \vA^i = \begin{pmatrix} \cos\theta^i & -\sin\theta^i \\ \sin\theta^i & \cos\theta^i\end{pmatrix}\;.
\end{align}
For deformable bodies, the distance between two arbitrary points on the deformable body does not, in general, remain constant because of the relative motion between the particles forming the body. In this case, the vector $\vu^i_P$ can be written as:
\begin{align}
  \vu^{ij}_P=\vu^{ij}_0+\vu^{ij}_f = \vu^{ij}_0+\vS^{ij}\vq^{ij}_f\;.
\end{align}
According to Fig.~\ref{fig:disp_vector} and what we have discussed so far, we describe the new position of an arbitrary point $P^i$ on the flexible body based on reference and elastic coordinates as:
\begin{align}
  \vr^{ij}_P=\vR^i+\vA^i\vu^{ij}_P=\vR^i+\vA^i\left(\vu^i_0+\vS^{ij}\vq^{ij}_f \right)\;.\label{eq:deformable_position}
\end{align}
We summarize all the unknowns which are necessary to calculate the new position of the arbitrary point $P$ in the vector $\vq^{ij}$:
\begin{align}
  \vq^{ij} = \begin{pmatrix} \vR^i \\ \theta^i \\ \vq^{ij}_f\end{pmatrix}\;.
\end{align}
Differentiating \eqref{eq:deformable_position} with respect to time yields
\begin{align}
  \dot{\vr}^{ij}_P=\dot{\vR}^i+\dot{\vA}^i\vu^{ij}_P+\vA^i\dot{\vu}^{ij}_P= \dot{\vR}^i+\dot{\vA}^i\vu^{ij}_P+\vA^i\vS^{ij}\dot{\vq}^{ij}_f \;,
\end{align}
where in case of planar motion, we have $\dot{\vA}=\vA_{\theta}\dot{\theta}$. Then, the velocity vector can be written as
\begin{align}
  \dot{\vr}^{ij}_P = \begin{pmatrix}\vI & \vA^i_{\theta}\vu^{ij} & \vA^i\vS^{ij} \end{pmatrix}\begin{pmatrix} \dot{\vR}^i \\ \dot{\theta}^i \\ \dot{\vq}^{ij}_f\end{pmatrix} \;, \label{eq:deformable_velocity}
\end{align}
where $\vI$ is identity matrix, and $\vA^i_{\theta}$ is the partial derivative of the transformation matrix with respect to the rotational coordinate $\theta^i$:
\begin{align}
  \vA^i_{\theta} = \begin{pmatrix} -\sin\theta^i & -\cos\theta^i \\ \cos\theta^i & -\sin\theta^i\end{pmatrix}\;.
\end{align}
 Equation~\eqref{eq:deformable_velocity} can also be written as
\begin{align}
  \dot{\vr}^{ij}_P = \vL^{ij} \dot{\vq}^{ij}\;,
\end{align}
where $\vL^{ij}= \begin{pmatrix}\vI & \vA^i_{\theta}\vu^{ij} & \vA^i\vS^{ij} \end{pmatrix}$. 
\subsection{Constructing the mass matrix} \label{sec:mass}
In this section, we develop the kinetic energy of deformable bodies and point out the differences between the inertia properties of deformable bodies that undergo finite rotations and the inertia properties of both rigid and structural systems. Then, we explain how to assemble the total mass matrix considering two rigid bodies (crank and slider) and joint constraints.\par
For constructing the mass matrix the following definition of the kinetic energy is used for element $j$ of deformable body $i$:
\begin{align}
  T^{ij} = \frac{1}{2}\int_{V^{ij}}\rho^{ij}\dot{\vr}^{{ij}^T}\dot{\vr}^{ij} dV^{ij}\;,
\end{align}
where $\rho^{ij}$ and $V^{ij}$ are, respectively, the mass density and volume of the element $j$, $\dot{\vr}^{ij}$ is the global velocity vector of an arbitrary point of the element. Using the expression of the velocity vector \eqref{eq:deformable_velocity}, we write the kinetic energy as
\begin{align}
  T^{ij} = \frac{1}{2}\dot{\vq}^{ij^T}  \left[\int_{V^{ij}}\rho^{ij} \vL^{ij^T}\vL^{ij} d V^{ij}\right]  \dot{\vq}^{ij}\;, \label{eq:KE_flexible}
\end{align}
where $\vM^{ij}$ is recognized as the symmetric mass matrix of body $i$. It is defined as
\begin{align}
  \vM^{ij}&=\int_{V^{ij}}\rho^{ij} \vL^{ij^T}\vL^{ij} \mathrm{d} V^{ij} =\int_{V^{ij}}\rho^{ij}  \begin{pmatrix} \vI \\ \left(\vA^i_{\theta}\vu^{ij}\right)^T \\ \left(\vA^i\vS^{ij}\right)^{T} \end{pmatrix} \begin{pmatrix}\vI & \vA^i_{\theta}\vu^{ij} & \vA^i\vS^{ij} \end{pmatrix} \mathrm{d} V^{ij} \nonumber\\
          &=\int_{V^{ij}}\rho^{ij} \begin{pmatrix}\vI & \vA^i_{\theta}\vu^{ij} & \vA^i\vS^{ij} \\ & \vu^{ij^T}\vu^{ij} & \vu^{ij^T}\tilde{\vI}\vS^{ij} \\ \text{symm.} & & \vS^{ij^T}\vS^{ij}     \end{pmatrix} \mathrm{d} V^{ij} \;,	\label{eq:element_mass}
\end{align}
where the orthogonality of the transformation matrix $\vA^{i^T}\vA^i=\vI$ is used in order to simplify the sub-matrix in the lower right-hand corner. Further, $\vA^{i^T}_{\theta}\vA^i=\tilde{\vI}$ with
\begin{align}
 \tilde{\vI} = \begin{pmatrix} 0 & 1 \\ -1 & 0 \end{pmatrix}\;.
\end{align}
The mass matrix~\eqref{eq:element_mass} can also be written as
\begin{align}
  \vM^{ij}=\begin{pmatrix}\vm_{RR} & \vm_{R\theta} & \vm_{R f} \\ & {m}_{\theta\theta} & \vm_{\theta f} \\  \text{symm.} & & \vm_{ff} \end{pmatrix}^{ij}\;,
\end{align}
where
\begin{align}
	\begin{aligned}
    \vm^{ij}_{RR}&=\int_{V^{ij}}\rho^{ij}\vI\mathrm{d}V^{ij}\;, & \vm^{ij}_{R\theta}&=\int_{V^{ij}}\rho^{ij}\vA^i_{\theta}\vu^{ij}  \mathrm{d}V^{ij}\;,\\
    \vm^{ij}_{R f}&=\int_{V^{ij}}\rho^{ij}\vA^i\vS^{ij}\mathrm{d}V^{ij}\;, & {m}^{ij}_{\theta\theta}&=\int_{V^{ij}}\rho^{ij}\vu^{ij^T}\vu^{ij}\mathrm{d}V^{ij}\;,\\
    \vm^{ij}_{\theta f}&=\int_{V^{ij}}\rho^{ij}\vu^{ij^T}\tilde{\vI}\vS^{ij}\mathrm{d}V^{ij}\;, & \vm^{ij}_{ff}&=\int_{V^{ij}}\rho^{ij}\vS^{ij^T}\vS^{ij}\mathrm{d}V^{ij}\;.
	\end{aligned}\label{eq:element_mass_submatrix}
\end{align}
Note that the two sub-matrices $\vm^{ij}_{RR}$ and $\vm^{ij}_{ff}$, which are associated, respectively, with the translational reference and elastic coordinates, are constant. Other matrices, however, depend on the system generalized coordinates.\par
The mass matrix in the case of a rigid body motion can be written as
\begin{align}
  \vM^{i}_{\text{rigid}}=\begin{pmatrix}\vm_{RR} & \vm_{R \theta} \\ \text{symm.} & m_{\theta\theta}  \end{pmatrix}^{i}\;.
\end{align}
In the case of structural systems, the reference coordinates remain constant with respect to time and the mass matrix of the body in this case is the constant matrix $\vm^{ij}_{ff}$. When a deformable body undergoes rigid body motion, the mass matrix is defined by \eqref{eq:element_mass_submatrix} and the sub-matrices $\vm^{ij}_{R f}$ and $\vm^{ij}_{\theta f}$ represent the coupling between the reference motion and the elastic deformation.\par
\begin{figure}[ht]
  \centering
  \footnotesize
  \def\svgwidth{0.75\columnwidth}
  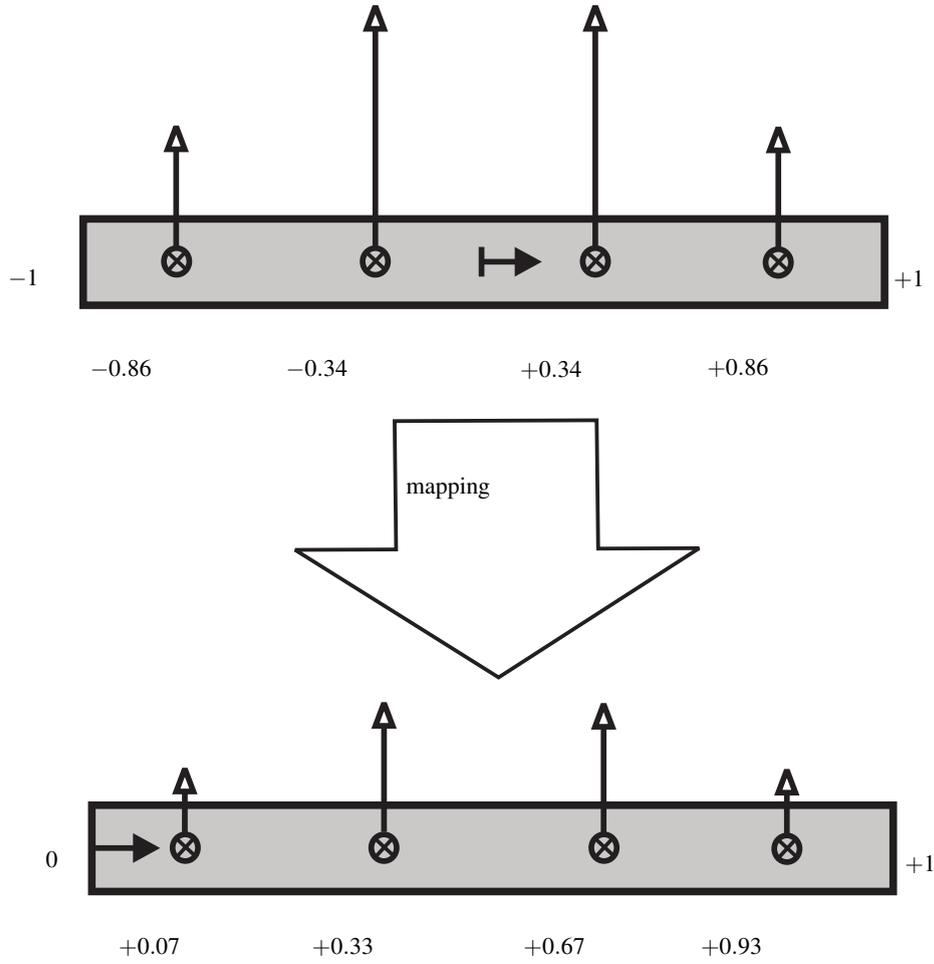
  \caption{Gaussian point distribution for exact integration of polynomial representations in mass matrix evaluations.}
  \label{fig:gp_mass}
\end{figure}

In the following, we detail each sub-matrix to find a simplified version for more efficient calculation. The matrix $\vm^{ij}_{RR}$ can be defined as
\begin{align}
  \vm^{ij}_{RR}=\int_{V^{ij}}\rho^{ij}\vI\mathrm{d}V^{ij}=\begin{pmatrix} m^{ij} & 0 \\  0 & m^{ij}  \end{pmatrix}=\vI^{ij}_0\;,
\end{align}
where $\vI$ is the identity matrix and $m^{ij}$ is the mass of the element $j$ of the deformable body $i$. We write the sub-matrix $\vm^{ij}_{R \theta}$ as
\begin{align}
  \vm^{ij}_{R \theta}&=\int_{V^{ij}}\rho^{ij}\vA^i_{\theta}\vu^{ij}\mathrm{d} V^{ij}=\vA^i_{\theta}\int_{V^{ij}}\rho^{ij}\left[\vu^{ij}_0+\vu^{ij}_f\right] \mathrm{d} V^{ij}=\vA^i_{\theta} \left[\vI^{ij}_1+\bar{\vS}^{ij}\vq^{ij}_f\right]\;,
\end{align}
where the matrices $\vI^{ij}_1$ and $\bar{\vS}^{ij}$ are defined as
\begin{align}
  \vI^{ij}_1= \int_{V^{ij}}\rho^{ij}\vu^{ij}_0 \mathrm{d} V^{ij}\;, \quad \bar{\vS}^{ij}=\int_{V^{ij}}\rho^{ij}\vS^{ij} \mathrm{d} V^{ij}\;. \label{eq:S_bar} 
\end{align}
The vector $\vI^{ij}_1$ is the moment of mass of the body about the axes of the body reference in the undeformed state. It may vanish if the origin of the body reference is initially attached to the body center of mass. The vector $\bar{\vS}^{ij}\vq^{ij}_f$ represents the change in the moment of mass due to the deformation. Using \eqref{eq:element_mass_submatrix}, we verify that
\begin{align}
  \vm^{ij}_{R f}=\vA^i\bar{\vS}^{ij}\;.
\end{align}
The expression for $m^{ij}_{\theta\theta}$ is
\begin{align}
  m^{ij}_{\theta\theta}&= \int_{V^{ij}}\rho^{ij}\left[\vu^{ij}_0+\vu^{ij}_f\right]^T\left[\vu^{ij}_0+\vu^{ij}_f\right] \mathrm{d} V^{ij}= \int_{V^{ij}}\rho^{ij}\left[\vu^{ij^T}_0\vu^{ij}_0+2\vu^{ij^T}_0\vu^{ij}_f+\vu^{ij^T}_f\vu^{ij}_f\right] \mathrm{d} V^{ij}\nonumber\\
                       &=\left( m^{ij}_{\theta\theta}\right)_{rr}+\left( m^{ij}_{\theta\theta}\right)_{rf}+\left( m^{ij}_{\theta\theta}\right)_{ff}\;,
\end{align}
in which the sub-matrix $ m^{ij}_{\theta\theta}$ reduces to a scalar that can be written as the sum of three components. The first component, $\left( m^{ij}_{\theta\theta}\right)_{rr}$, is the mass moment of inertia in the undeformed state:
\begin{align}
  \left( m^{ij}_{\theta\theta}\right)_{rr}&= \int_{V^{ij}}\rho^{ij}\vu^{ij^T}_0\vu^{ij}_0 \mathrm{d} V^{ij}=\int_{V^{ij}}\rho^{ij}\left[\left(x^{ij}\right)^2+\left(y^{ij}\right)^2 \right] \mathrm{d} V^{ij}=I^{ij}_2 \;.
\end{align}
Clearly, this integral has a constant value and does not depend on the body deformation. The last two scalar components, $\left( m^{ij}_{\theta\theta}\right)_{rf}$ and $\left( m^{ij}_{\theta\theta}\right)_{ff}$, represent the change in the mass moment of inertia of the body due to deformation. These two components are evaluated according to
\begin{align}
  \left( m^{ij}_{\theta\theta}\right)_{rf}&=2\int_{V^{ij}}\rho^{ij}\vu^{ij^T}_0\vu^{ij}_f \mathrm{d} V^{ij}=2\left[\int_{V^{ij}}\rho^{ij}\vu^{ij^T}_0\vS^{ij} \mathrm{d} V^{ij} \right]\vq^{ij}_f=2~\vI^{ij}_3~\vq^{ij}_f\;,\\
  \left( m^{ij}_{\theta\theta}\right)_{ff}&=\int_{V^{ij}}\rho^{ij}\vu^{ij^T}_f\vu^{ij}_f \mathrm{d} V^{ij}=\vq^{ij^T}_f \left[\int_{ V^{ij}}\rho^{ij}\vS^{ij^T}\vS^{ij} \mathrm{d} V^{ij}\right]\vq^{ij}_f\;.
\end{align}
If we use definition~\eqref{eq:element_mass_submatrix}, we can write the following:
\begin{align}
  \left( m^{ij}_{\theta\theta}\right)_{ff}=\vq^{ij^T}_f\vm^{ij}_{ff}\vq^{ij}_f\;,
\end{align}
where 
\begin{align}
  \vm^{ij}_{ff}=\int_{V^{ij}}\rho^{ij}\vS^{ij^T}\vS^{ij}\mathrm{d} V^{ij}=\vS^{ij}_{ff} \;.\label{eq:S_ff}
\end{align}
Finally, we introduce
\begin{align}
  \vm^{ij}_{\theta f}=\int_{V^{ij}}\rho^{ij}\left[\vu^{ij}_0+\vu^{ij}_f\right]^T\tilde{\vI}~\vS^{ij}\mathrm{d} V^{ij} = \vI^{ij}_4+\vq^{{ij}^T}_f\tilde{\vS}^{ij}\;,
\end{align}
where the constant skew symmetric matrix $\tilde{\vS}^{ij}$ is defined as
\begin{align}
  \vI^{ij}_4=\int_{V^{ij}}\rho^{ij}\vu^{ij^T}_0\tilde{\vI}~\vS^{ij}\mathrm{d} V^{ij}\;,\quad \tilde{\vS}^{ij}=\int_{V^{ij}}\rho^{ij}\vS^{ij^T}\tilde{\vI}~\vS^{ij}\mathrm{d} V^{ij} \;.\label{eq:S_tilde}         
\end{align}
We conclude that, to completely describe the inertia properties of the deformable body in plane motion, a set of inertia shape integrals is required. These integrals, which depend on the assumed displacement field, can be obtained using the Gaussian quadrature method. As we are dealing with polynomials for describing the displacement field, for finding the optimum number of Gaussian points, we have to know the highest degree of the polynomials which appear in the calculation of the mass matrix. We expect polynomials with degree 6 at most, which need 4 Gaussian points to yield an exact integration. Figure~\ref{fig:gp_mass} shows Gaussian points which are necessary to evaluate exact values on a sample element. To evaluate the integrals for $\xi \in \left[0,1\right]$, we have to map Gaussian points and weights.\par
\begin{figure}[ht]
  \centering
  \footnotesize
  \def\svgwidth{\columnwidth}
  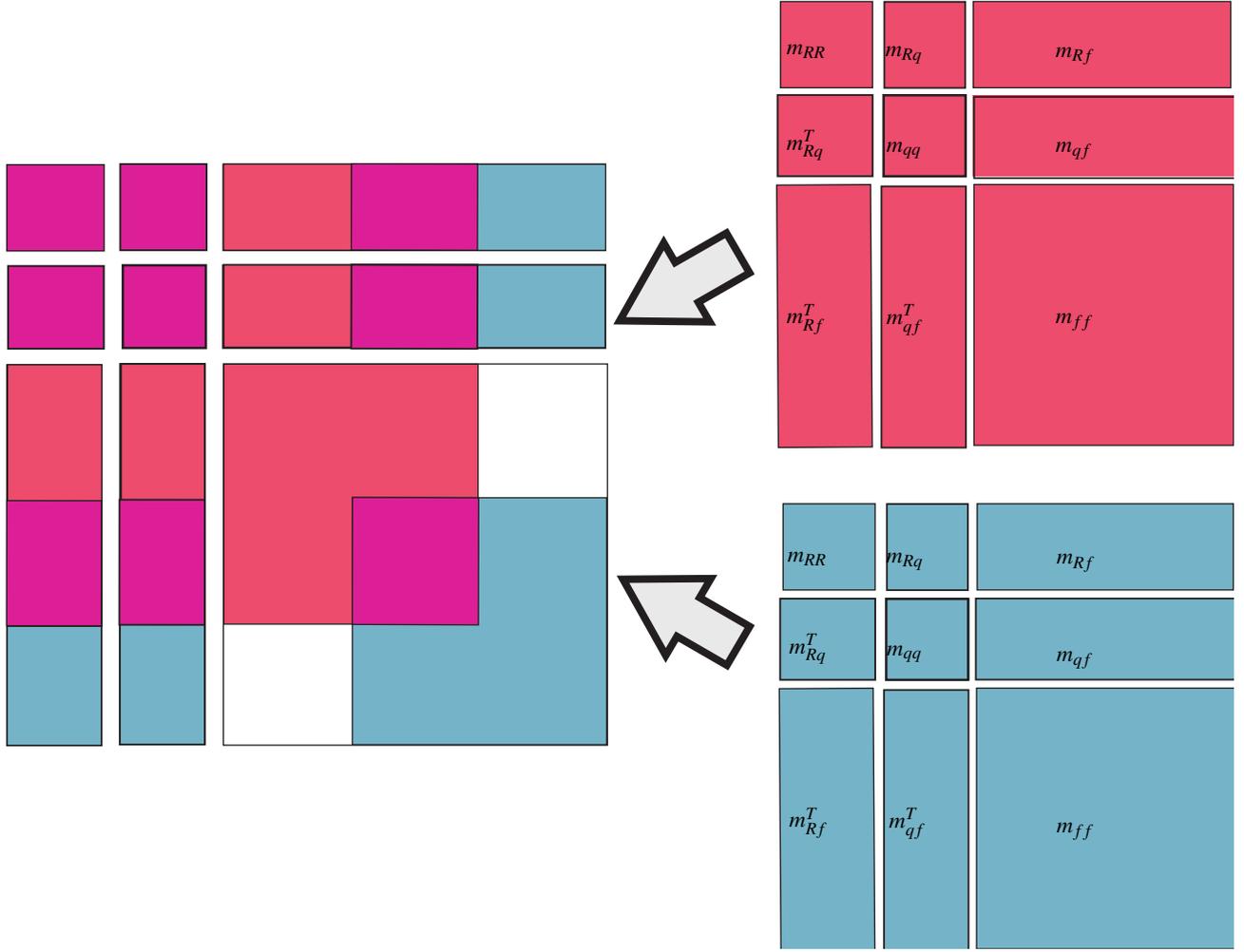
  \caption{Mass matrix results from assembling of element mass matrices $j$ and $j+1$; purple: summation of mutual DOFs, red: element mass matrix $j$, blue: element mass matrix $j+1$.}
  \label{fig:mass_assemble}
\end{figure}

Once we have calculated the mass matrix on element level, we have to assemble the total mass matrix considering mutual degrees of freedom. One should note that entities in the mass matrix of different elements share the same reference coordinates, but perhaps different elastic coordinates. Figure~\ref{fig:mass_assemble} shows the regions in which we have the overlap (purple color) between element $j$ (light red) and element $j+1$ (light blue) in the assembly process.
As the slider and the rod are attached to each other at one revolute joint, the end point of the flexible rod and the mass center of the slider have the same translational velocity. Therefore, we add the effect of the rigid slider to the mass matrix of the last element of the flexible rod, also by adding one additional degree of freedom in the total mass matrix for $\theta_3$ (Fig.~\ref{fig:slider_floating_frame}):
\begin{align}
  T^3 = \frac{1}{2}m_3~\dot{\vq}^{il^T}\left[\vL^{il^T}\vL^{il}\right]\dot{\vq}^{il}+\frac{1}{2}J_3~\dot{\theta}^2_3 \;,\label{eq:slider_mass}
\end{align}
where index 3 is related to the slider (third body), and $l$ is related to the last element of the mesh which is connected to the slider by means of a revolute joint. We assemble the first term of \eqref{eq:slider_mass} in the total mass matrix regarding mutual degrees of freedom and we add one additional row and column for the new degree of freedom $\theta_3$ which contain zero everywhere except one diagonal term which is $J_3$.\par
To add the effect of the rigid crank mass into the total mass matrix, first we define the constraint for the connecting joint between crank and flexible rod:
\begin{align}
  \vC^1 = \vR^i-\begin{pmatrix} l_1\cos\theta_1 \\ l_1\sin\theta_1 \end{pmatrix}=\vnull \;,\label{eq:first_con}
\end{align}
where $\vR^i = \begin{pmatrix} x^i & y^i \end{pmatrix}^T$ is the translational coordinate of the rod reference frame. Equation~\eqref{eq:first_con} shows that the rod reference coordinate and therefore its derivative with respect to time can be expressed in terms of $\theta_1$. Therefore and for rewriting the total mass matrix in terms of $\theta_1$, we need the time derivative of the constraint $\vC^1$:
\begin{align}
  \dot{\vC}^1 = \dot{\vR}^i - \begin{pmatrix} -l_1\sin\theta_1 \\ l_1\cos\theta_1 \end{pmatrix}\dot{\theta}_1=
  \dot{\vR}^i - \vC^1_{\theta_1}\dot{\theta}_1 =\vnull \;.\label{eq:first_con_dot}
\end{align}
Using \eqref{eq:first_con_dot} and keeping in mind the kinetic energy formulation \eqref{eq:KE_flexible}, we modify the total mass matrix such that it only depends on the coordinate $\theta_1$ according to the constraint $\vC^1$:
\begin{align}
	\begin{aligned}
    \bar{m}_{RR}&=\vC^{1^T}_{\theta_1}\vm_{RR}\vC^1_{\theta_1}+I_1\;,\\
    \bar{m}_{R\theta}&=\vC^{1^T}_{\theta_1}\vm_{R\theta}\;, \quad \bar{m}_{\theta R} = \bar{m}^{T}_{R\theta}\;, \\
   \bar{\vm}_{Rf}&=\vC^{1^T}_{\theta_1}\vm_{R f}\;, \quad \bar{\vm}_{fR} = \bar{\vm}^{T}_{Rf} \;,
	\end{aligned} \label{eq:m_bar}
\end{align}
where $I_1=J_1+\frac{1}{2} m_1 l^2_1$ is the rigid crank mass moment of inertia with respect to the joint. We consider the contribution of the rigid crank in the total mass matrix according to the kinetic energy of the crank in terms of its only degree of freedom $\theta_1$:
\begin{align}
  T^1 = \frac{1}{2}I_1\dot{\theta}^2_1 \;,\label{eq:crank_mass}
\end{align}
where index 1 denotes the body number of the crank.\par
\begin{figure}[ht]
  \centering
  \footnotesize
  \def\svgwidth{0.7\columnwidth}
  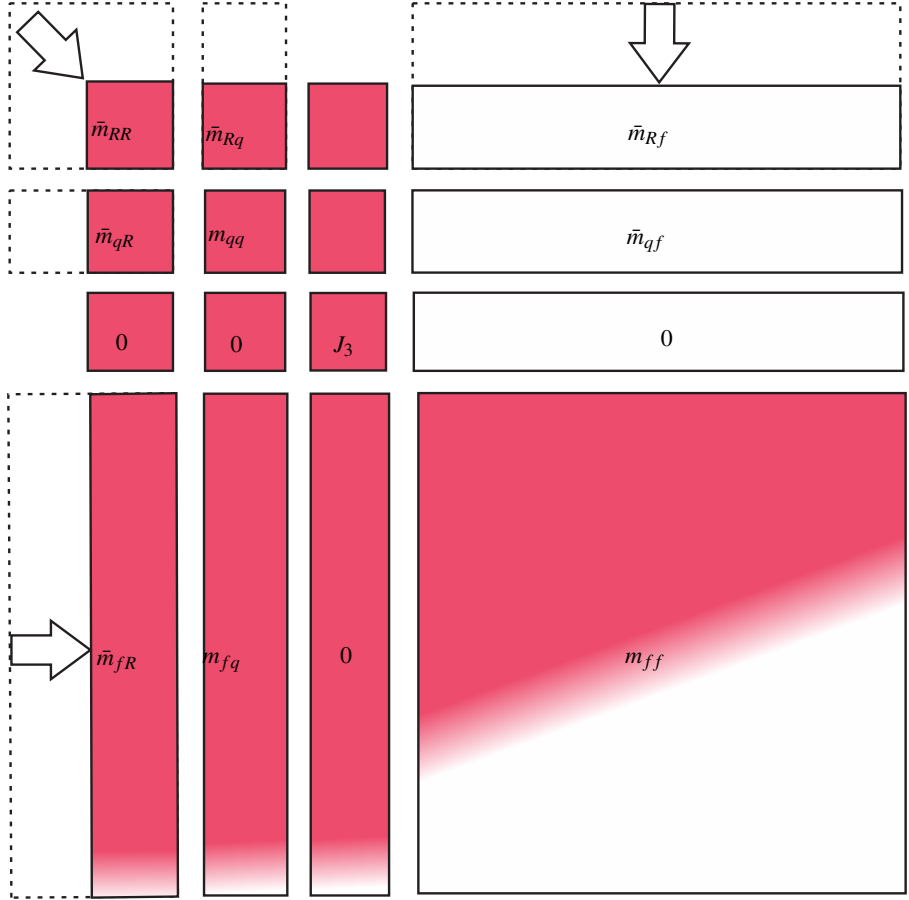
  \caption{Total mass matrix for the slider-crank mechanism considering joint constraints.}
  \label{fig:total_mass}
\end{figure}

Figure \ref{fig:total_mass} shows the final mass matrix configuration after considering all the effects in the slider-crank mechanism.
\subsection{Quadratic velocities}
For the quadratic velocity vector, we need to calculate the derivative of the mass matrix with respect to time and the derivative of the kinetic energy with respect to the degrees of freedom. For the time-derivative of the mass matrix, we write
\begin{align}
	\begin{aligned}
    \dot{\vM}^{ij}&=\int_{V^{ij}}\rho^{ij} \left( \dot{\bar{\vL}}^{ij^T}\bar{\vL}^{ij}+\bar{\vL}^{ij^T}\dot{\bar{\vL}}^{ij} \right) \mathrm{d} V^{ij}=\begin{pmatrix}\dot{\bar{m}}_{RR} & \dot{\bar{m}}_{R\theta} & \dot{\bar{\vm}}_{R f} \\ & \dot{m}_{\theta\theta} & \dot{\vm}_{\theta f} \\ \text{symm.} & & \dot{\vm}_{ff} \end{pmatrix}^{ij}\;,\\
	\end{aligned}\label{eq:element_Dmass_submatrix}
\end{align}
where $\bar{\vL}^{ij}$ and $\dot{\bar{\vL}}^{ij}$ are the modified versions of $\vL^{ij}$ and $\dot{\vL}^{ij}$ according to the hinge constraint defined in \eqref{eq:first_con}:
\begin{align}
  \bar{\vL}^{ij}&=\begin{pmatrix}\vC^{1}_{\theta_1} & \vA^i_{\theta}\vu^{ij} & \vA^i\vS^{ij} \end{pmatrix}, \\
  \dot{\bar{\vL}}^{ij}&=\begin{pmatrix} \vC^{1}_{\theta_1\theta_1}\dot{\theta_1} & \vA^i_{\theta}\vS^{ij}\dot{\vq}^{ij}_f-\vA^i\left(\vu^{ij}_0+\vS^{ij}\vq^{ij}_f\right)\dot{\theta_2} & \vA^i_{\theta}\vS^{ij}\dot{\theta_2} \end{pmatrix}\;.
\end{align}
Thereby, $\vC^1_{\theta_1\theta_1}$ is the second derivative of the joint constraint:
\begin{align}
 \dot{\vC}^1_{\theta_1} = \begin{pmatrix} l_1\cos\theta_1 \\ l_1\sin\theta_1 \end{pmatrix}\dot{\theta}_1 = \vC^1_{\theta_1\theta_1}\dot{\theta}_1\;. \label{eq:first_con_ddot}
\end{align}
After assembling the element contributions according to the pattern described in Fig.~\ref{fig:mass_assemble}, we add the effect of the rigid slider to the matrix of the last element of the flexible rod, noting that $\dot{J}_3$  is zero:
\begin{align}
  \dot{\vM}^3 = m_3\left( \dot{\bar{\vL}}^{il^T}\bar{\vL}^{il}+\bar{\vL}^{il^T}\dot{\bar{\vL}}^{il} \right)\;.\label{eq:slider_Dmass}
\end{align}
Note that the effect of the rigid crank is already taken into account by modifying the $L$ matrix.\par
To calculate the second term of the quadratic velocity, first we write the expression for the total kinetic energy as
\begin{align}
  T&=\frac{1}{2}\bar{m}_{RR}\dot{\theta}^2_1 + \dot{\theta}_1\bar{m}_{R\theta}\dot{\theta}_2 + \frac{1}{2}m_{\theta\theta}\dot{\theta}^2_2 + \dot{\theta}_1\bar{\vm}_{Rf}\dot{\vq}_f + \dot{\theta}_2\vm_{\theta f}\dot{\vq}_f+\frac{1}{2}\vm_{f f}\dot{\vq}^2_f + \frac{1}{2}J_3\dot{\theta}^2_3\;.
\end{align}
Keeping in mind definition \eqref{eq:m_bar}, the derivative of the kinetic energy with respect to the nodal coordinates is
\begin{align}
  \frac{\partial T}{\partial\vq} = \begin{pmatrix} \dot{\vC}^{1^T}_{\theta_1}\left(m_2+m_3\right)\vC^1_{\theta_1}\dot{\theta_1}+\dot{\vC}^{1^T}_{\theta_1}\vA^i_{\theta}\left(\vI_1+\bar{\vS}\vq_f\right)\dot{\theta_2}+\dot{\vC}^{1^T}_{\theta_1}\vA\bar{\vS}\dot{\vq}_f \\
  -\vC^{1^T}_{\theta_1}\vA\left(\vI_1+\bar{\vS}\vq_f\right)\dot{\theta_1}\dot{\theta_2} + \vC^{1^T}_{\theta_1}\vA^i_{\theta}\bar{\vS}\dot{\theta_1}\dot{\vq}_f\\
  0\\
  \vC^{1^T}_{\theta_1}\vA^i_{\theta}\bar{\vS}\dot{\theta_1}\dot{\theta_2} + \left(\vI_3+\vS_{ff}\right)\vq_f\dot{\theta}^2_2 + \tilde{\vS}\dot{\theta}_2\dot{\vq}_f
\end{pmatrix}\;,
\end{align}
where $\bar{\vS}$, $\vS_{ff}$, $\tilde{\vS}$, $\vI_1$, and $\vI_3$ are the assembled versions of $\bar{\vS}^{ij}$, $\vS^{ij}_{ff}$, $\tilde{\vS}^{ij}$, $\vI^{ij}_1$ and $\vI^{ij}_3$ defined in \ref{sec:mass}. Finally, we have
\begin{align}
  \vQ_v = -\begin{pmatrix}\dot{\bar{m}}_{RR} & \dot{\bar{m}}_{R\theta} & 0 & \dot{\bar{\vm}}_{R f} \\  & \dot{m}_{\theta\theta} & 0 & \dot{\vm}_{\theta f} \\ & & 0 & 0 \\ \text{symm.} & & & \dot{\vm}_{ff}\end{pmatrix} \begin{pmatrix}\dot{\theta_1} \\ \dot{\theta_2} \\ \dot{\theta_3} \\ \dot{\vq}_f\end{pmatrix} + \frac{\partial T}{\partial\vq}\;.
\end{align}
\subsection{Stiffness matrix}
\begin{figure}[ht]
  \centering
  \includegraphics[width=0.8\columnwidth]{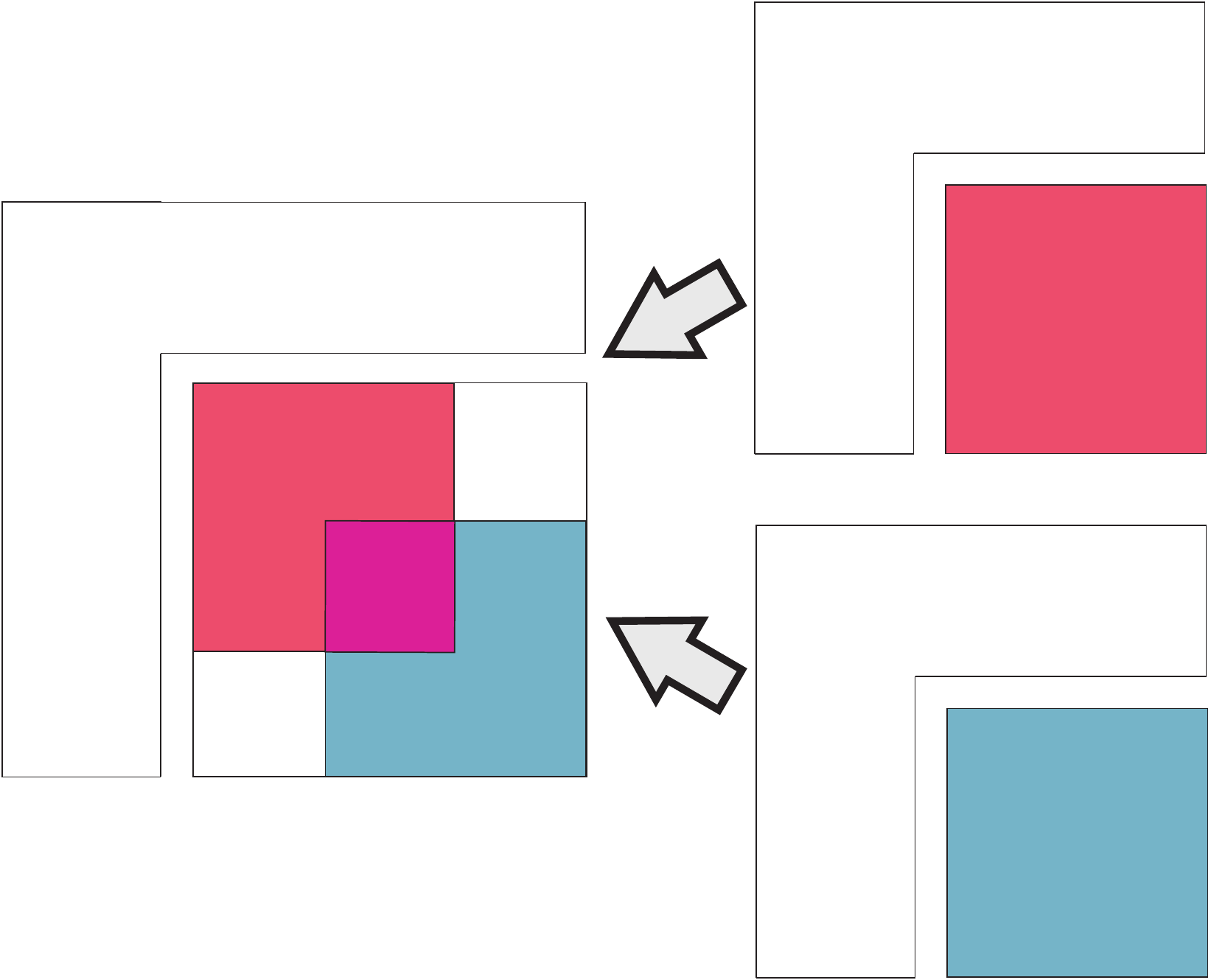}
  \caption{Stiffness matrix assembling of two element $j$ and $j+1$ stiffness matrices; purple: summation of mutual DOFs, blank: zero entries.}
  \label{fig:stiffness_assemble}
\end{figure}

Considering a linear isotropic material, the virtual work due to the elastic forces for element $j$ of body $i$ can be written as
\begin{align}
  \delta W^{ij}_s=-\int_{V^{ij}}\vsigma^{{ij}^T}\delta\vepsilon^{ij}\mathrm{d} V^{ij}\;,
\end{align}
where $\vsigma^{ij}$ and $\vepsilon^{ij}$ are, respectively, the stress and strain tensors. Since the rigid body motion corresponds to the case of constant strains and since we defined the deformation of flexible bodies with respect to the body reference, the strain displacement relations can be written in the following form:
\begin{align}
  \vepsilon^{ij}=\vD^{ij}\vu^{ij}_f\;,
\end{align}
where $\vD^{ij}$ is a differential operator. We write
\begin{align}
  \vepsilon^{ij}=\vD^{ij}\vS^{ij}\vq^{ij}_f\;.
\end{align}
For a linear isotropic material, the constitutive equations can be written as
\begin{align}
  \vsigma^{ij}=\vC^{ij}\vepsilon^{ij}\;,
\end{align}
where $\vC^{ij}$ is the symmetric matrix of elastic coefficients. We conclude
\begin{align}
  \vsigma^{ij}=\vC^{ij}\vD^{ij}\vS^{ij}\vq^{ij}_f\;,
\end{align}
where the stress tensor is written in terms of the elastic generalized coordinates of body $i$. Finally, we have
\begin{align}
  \delta W^{ij}_s=-\vq^{{ij}^T}_f \left[\int_{V^{ij}} \left(\vD^{ij}\vS^{ij}\right)^T\vC^{ij}\vD^{ij}\vS^{ij} \mathrm{d} V^{ij} \right] \delta\vq^{ij}_f = -\vq^{{ij}^T}_f\vK^{ij}_{ff} \delta\vq^{ij}_f \;, \label{eq:elastic_virtual_work}
\end{align}
where
\begin{align}
  \vK^{ij}_{ff} =\int_{V^{ij}} \left(\vD^{ij}\vS^{ij}\right)^T\vC^{ij}\vD^{ij}\vS^{ij} \mathrm{d} V^{ij}\;. \label{eq:stiffness_matrix}
\end{align}
Neglecting the shear deformation and using the assumptions of Euler Bernoulli beam theory, the strain energy for the element $j$ of the elastic rod can be written as
\begin{align}
U^{ij} = \frac{1}{2} \int^{l^{ij}}_0 \begin{pmatrix} u'_{f1} & u''_{f2} \end{pmatrix}^{ij} \begin{pmatrix} E^{ij}A^{ij} & 0 \\ 0 & E^{ij}I^{ij} \end{pmatrix} \begin{pmatrix}{c} u'_{f1} \\ u''_{f2} \end{pmatrix} ^{ij}\mathrm{d} x\;, \label{eq:strain_energy}
\end{align}
where $l^{ij}$ is the length of element $j$ of the beam, $E^{ij}$ is the modulus of elasticity of the beam, $A^{ij}$ is the cross-sectional area, $I^{ij}$ is the second moment of area, $u^{ij}_{f1}$ and $u^{ij}_{f2}$ are the axial and transverse displacements of element $j$ respectively and $\left('\right)$ denotes differentiation with respect to the spatial coordinate. Taking the derivative of \eqref{eq:strain_energy} with respect to $\vq_f$ and comparing with \eqref{eq:elastic_virtual_work}, we write
\begin{align}
  \vC^{ij} = \begin{pmatrix} E^{ij}A^{ij} & 0 \\ 0 & E^{ij}I^{ij} \end{pmatrix}\;.
\end{align}
Using the shape functions introduced in \eqref{eq:shape_functions}, we conclude:
\begin{align}
\vD^{ij}\vS^{ij} = \vS'^{ij} &= \begin{pmatrix} \frac{1}{l^{ij}} & 0 \\ 0 & \frac{1}{l^{ij^2}} \end{pmatrix} \begin{pmatrix} -1 & 0 & 0 & 1 & 0 & 0 \\ 0 & -6+12 \xi & l^{ij}\left(-4+6\xi\right) & 0 & 6-12\xi & l^{ij}\left(6\xi-2\right) \end{pmatrix} \;,
\end{align}
where the first matrix is the inverse of the Jacobian for transforming to local coordinates with $\xi=x/l$, and the second matrix is the derivative of $\vS^{ij}$ with respect to $\xi$. We calculate the stiffness matrix for each element using \eqref{eq:stiffness_matrix}. It is worth to mention that as we are dealing with derivatives of the shape functions, the highest degree of polynomials which appears in the integral of \eqref{eq:stiffness_matrix} is $2$. Therefore, we only need $2$ Gaussian points for the computation of the stiffness matrix on element level. Figure~\ref{fig:stiffness_assemble} shows the stiffness matrices for element $j$ and $j+1$ and the assembling procedure. Blank areas, which correspond to the reference degrees of freedom, i.e., $\theta_1$, $\theta_2$ and $\theta_3$, are filled with zeros.
\subsection{External forces}
The virtual work of all external forces $\vF^i$ acting on element $j$ of body $i$ in the multibody system can be written as
\begin{align}
  \delta W^{ij}_e=\int_{V^{ij}}  \vF^{ij^T} \delta\vr^{ij} \mathrm{d} V^{ij}\;,
\end{align}
where
\begin{align}
  \delta\vr^i = \begin{pmatrix}\vI & \vA^i_{\theta}\vu^{ij} & \vA^i\vS^{ij} \end{pmatrix}  \begin{pmatrix} \delta\vR^i \\ \delta\theta^i \\ \delta\vq^{ij}_f\end{pmatrix} = \vL^{ij}\delta\vq^{ij}\;.
\end{align}
Combining the above two equations, we get the following expression for the virtual work of external forces:
\begin{align}
  \delta W^{ij}_e=\begin{pmatrix}\vQ^{ij^T}_R & Q^{ij^T}_{\theta} & \vQ^{ij^T}_f \end{pmatrix}  \begin{pmatrix} \delta\vR^i \\ \delta\theta^i \\ \delta\vq^{ij}_f\end{pmatrix}=\vQ^{ij^T}_e \delta\vq^{ij}\;,
\end{align}
where
\begin{align}
  \begin{aligned}
    \vQ^{ij^T}_R&= \int_{V^{ij}} \vF^{ij^T} \mathrm{d} V^{ij}\;,\\  
    Q^{ij^T}_{\theta}&=\int_{V^{ij}}\vF^{ij^T}\vA^i_{\theta}\vu^{ij} \mathrm{d} V^{ij}\;,\\
    \vQ^{ij^T}_f&=\int_{V^{ij}}\vF^{ij^T}\vA^i\vS^{ij} \mathrm{d} V^{ij}\;.
  \end{aligned} \label{eq:element_force_submatrix}
\end{align}
Assembling of the external force vector can be done with the same procedure as for the mass matrix. To add the virtual work of the rigid slider, we use the revolute joint constraint between slider and rod
\begin{align}
  \vQ^{3^T}_e = \vF^T \vL^{il} \;,\label{eq:slider_force}
\end{align}
where index 3 is related to the slider (third body), and index $l$ is related to the last element of the mesh 
which is connected to the slider by means of a revolute joint. We assemble \eqref{eq:slider_force} in the total external force vector regarding mutual degrees of freedom and adding one additional row and column for the degree of freedom $\theta_3$ which is zero as there is no moment in this direction. To add the effect of the rigid crank mass into the total external force vector, we use the constraint for the connecting joint between crank and flexible rod \eqref{eq:first_con}. We modify the sub-vector introduced in \eqref{eq:element_force_submatrix}, such that it only depends on coordinate $\theta_1$ according to the constraint $\vC^1$:
\begin{align}
  \bar{Q}^{ij^T}_R=\vQ^{ij^T}_R \vC^1_{\theta_1}-m g\frac{l}{2}\sin\theta_1\;,
\end{align}
where the second term comes from the work of the gravitational force with the rigid crank. It is worth to mention that according to the polynomial order in the external force integral, 2 Gaussian points are sufficient for each element.
\subsection{Unilateral constraints}
\begin{figure}[ht]
  \centering
  \footnotesize
  \def\svgwidth{0.9\columnwidth}
  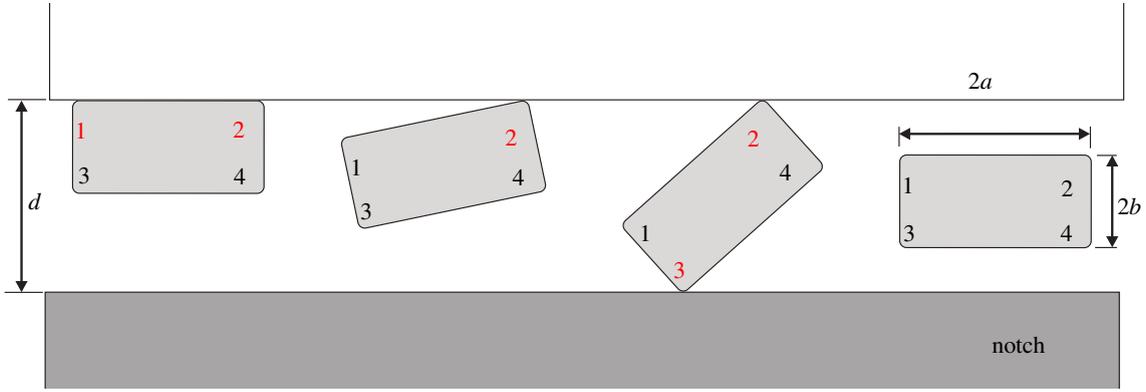
  \caption{Different scenarios for the slider and notch interaction. Red points indicate active contact points.}
  \label{fig:slider_scenarios}
\end{figure}

\begin{figure}[ht]
  \centering
  \footnotesize
  \def\svgwidth{0.8\columnwidth}
  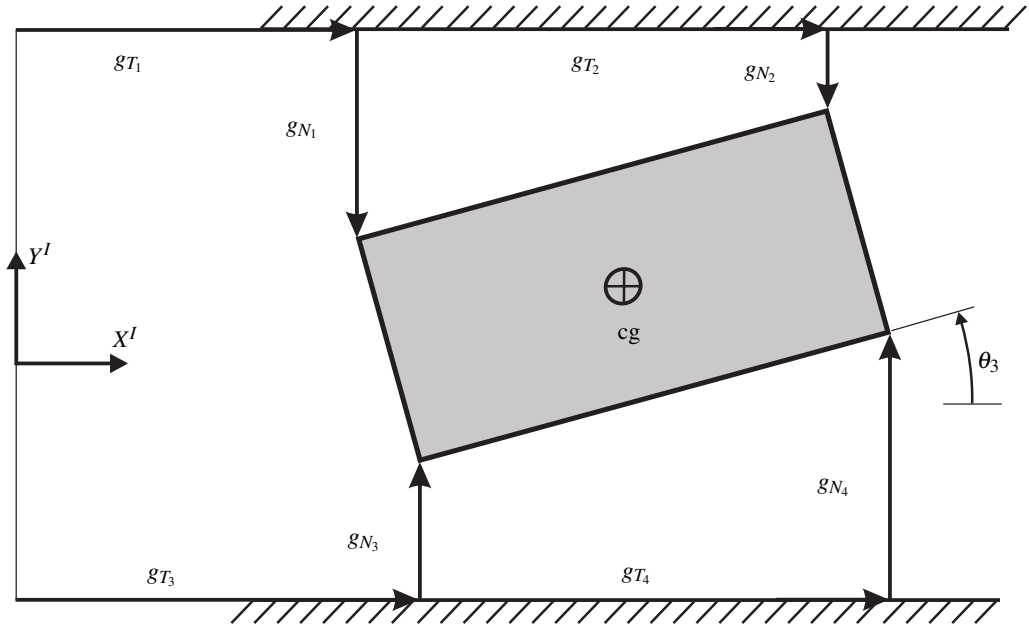
  \caption{Definition of the gap functions for the slider-crank mechanism with unilateral constraints.}
  \label{fig:slider_gap}
\end{figure}

Figure~\ref{fig:slider_scenarios} shows the geometric characteristics of the translational clearance joint. It is $2a$ the length and $2b$ the height of the slider. The height of the notch is given by $d$. The existence of a clearance in a translational joint introduces two extra degrees of freedom. Hence, the slider can move freely inside the guiding limits until it reaches the surfaces.
Considering the geometry of the slider according to Fig.~\ref{fig:slider_gap} and keeping in mind the revolute joint between slider and end point of the rod, we write the position of the slider center according to \eqref{eq:deformable_position} and \eqref{eq:first_con} as
\begin{align}
  \vr_{cg}=\begin{pmatrix} r^x_{cg} \\ r^y_{cg} \end{pmatrix}=\vR^i+\vA^i\vu_{cg}=\begin{pmatrix} l_1\cos\theta_1 \\ l_1\sin\theta_1 \end{pmatrix}+\vA^i\left(\theta_2\right) \left( \begin{pmatrix} l_2 \\ 0 \end{pmatrix}+\vS^{i l}\vq^{i l}_f \right) \;,\label{eq:cg_position}
\end{align}
where index $l$ is related to the last element of the mesh of the flexible rod with initial length $l_2$. Therefore, the nonlinear normal and tangential gap functions split up for each corner:
\begin{align}
  g_{N_1}\left(q\right)&=\frac{d}{2}-r^y_{cg}+a\sin\theta_3-b\cos\theta_3\;,\\
  g_{N_2}\left(q\right)&=\frac{d}{2}-r^y_{cg}-a\sin\theta_3-b\cos\theta_3\;,\\
  g_{N_3}\left(q\right)&=\frac{d}{2}+r^y_{cg}-a\sin\theta_3-b\cos\theta_3\;,\\
  g_{N_4}\left(q\right)&=\frac{d}{2}+r^y_{cg}+a\sin\theta_3-b\cos\theta_3\;,\\
  g_{T_1}\left(q\right)&=r^x_{cg}-a\cos\theta_3-b\sin\theta_3\;,\\
  g_{T_2}\left(q\right)&=r^x_{cg}+a\cos\theta_3-b\sin\theta_3\;,\\
  g_{T_3}\left(q\right)&=r^x_{cg}-a\cos\theta_3+b\sin\theta_3\;,\\
  g_{T_4}\left(q\right)&=r^x_{cg}+a\cos\theta_3+b\sin\theta_3\;.
\end{align}
The matrices of generalized force directions are the derivatives of $\vg_N$ and $\vg_T$ with respect to $\vq$:
\begin{align}
  \vW_N^T\left(\vq\right)&=\begin{pmatrix}
  -l_1\cos\theta_1 & -\left[\vA^i_{\theta}\vu_{cg}\right]_y & a\cos\theta_3+b\sin\theta_3 & -\left[\vA^i~\check{\vS}\right]_y \\
  -l_1\cos\theta_1 & -\left[\vA^i_{\theta}\vu_{cg}\right]_y & -a\cos\theta_3+b\sin\theta_3 & -\left[\vA^i~\check{\vS}\right]_y \\
  l_1\cos\theta_1 & \left[\vA^i_{\theta}\vu_{cg}\right]_y & -a\cos\theta_3+b\sin\theta_3 & \left[\vA^i~\check{\vS}\right]_y \\
  l_1\cos\theta_1 & \left[\vA^i_{\theta}\vu_{cg}\right]_y & a\cos\theta_3+b\sin\theta_3 & \left[\vA^i~\check{\vS}\right]_y \end{pmatrix}\;,\label{eq:generalized_force_direction_normal} \\
  \vW_T^T\left(\vq\right)&=\begin{pmatrix}
    -l_1\sin\theta_1 & \left[\vA^i_{\theta}\vu_{cg}\right]_x & a\sin\theta_3-b\cos\theta_3 & \left[\vA^i~\check{\vS}\right]_x \\
    -l_1\sin\theta_1 & \left[\vA^i_{\theta}\vu_{cg}\right]_x & -a\sin\theta_3-b\cos\theta_3 & \left[\vA^i~\check{\vS}\right]_x \\
    -l_1\sin\theta_1 & \left[\vA^i_{\theta}\vu_{cg}\right]_x & a\sin\theta_3+b\cos\theta_3 & \left[\vA^i~\check{\vS}\right]_x \\
  -l_1\sin\theta_1 & \left[\vA^i_{\theta}\vu_{cg}\right]_x & -a\sin\theta_3+b\cos\theta_3 & \left[\vA^i~\check{\vS}\right]_x \end{pmatrix} \label{eq:generalized_force_direction_tangential}
\end{align}
where $\left[*\right]_x$ and $\left[*\right]_y$ mean the component of the vector respectively in $x-$ and $y-$direction. The matrix $\check{\vS}$ is defined as follows:
\begin{align}
\check{\vS}=\left(\begin{array}{r|cccccc} 0 \cdots 0 & 0 & 0 & 0 & 1 & 0 & 0\\0 \cdots 0 & 0 & 0 & 0 & 0 & 1 & 0 \end{array}\right)_{2\times N}\;,
\end{align}
where $N=3\times\left(n_{ele}+1\right)$ is the number of elastic degrees of freedom. The last six columns are filled with $\vS\left(\xi=1\right)$ related to the degrees of freedom of the last element of the mesh.
\subsection{Boundary conditions}  
\begin{figure}[ht]
  \centering
  \footnotesize
  \def\svgwidth{\columnwidth}
  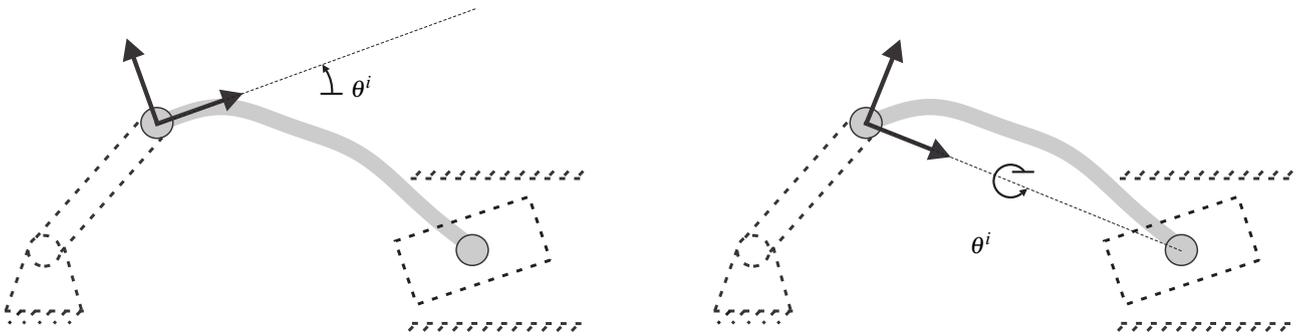
  \caption{Tangential and pinned reference system.}
  \label{fig:bc_floating_frame}
\end{figure}

When solving a finite element problem with the floating frame of reference formulation, the system becomes singular. This is because rigid body motion is added to the equations of motion at the same time as the deformation field also contains rigid body motion. A modal analysis would show that the first three eigenvalues are equal to zero which correspond to the degrees of freedom of the rigid body in planar motion. In this model, boundary conditions for the body reference system are introduced to avoid the singularity in the system of equations. There are different ways to define the boundary conditions. Shabana~\cite{Sha96} shows that two sets of modes associated with two sets of boundary conditions can be used to obtain the same solution if the coordinate system is properly selected. Therefore, the physical deformation is unique in the inertial frame. For a beam, it is most common to use a clamped-free or a simply supported reference system (Fig.~\ref{fig:bc_floating_frame}), where three conditions are given in both cases.
Clampled-free, i.e., tangential, means that the reference system is tangential to the beam deflection at the root of the beam, i.e., both displacements and rotations are equal to zero:
\begin{align}
  q^{i1}_{f1}=q^{i1}_{f2}=q^{i1}_{f3}=0\;.
\end{align}
For a simply supported, i.e., pinned, reference system, the root of the beam is locked and the end of the beam is moving but only along the local $x$ direction:
\begin{align}
  q^{i1}_{f1}=q^{i1}_{f2}=q^{il}_{f5}=0\;,
\end{align}
where index $l$ stands for the last beam element. In the present simulation, the tangential reference system is chosen. We compare different boundary conditions using modal coordinates in Sect.~\ref{sec:modal} as introduced in~\ref{sec:appE}.

\section{Modal analysis} \label{sec:appE}
\begin{figure}[ht]
  \centering
  \includegraphics[width=0.8\columnwidth]{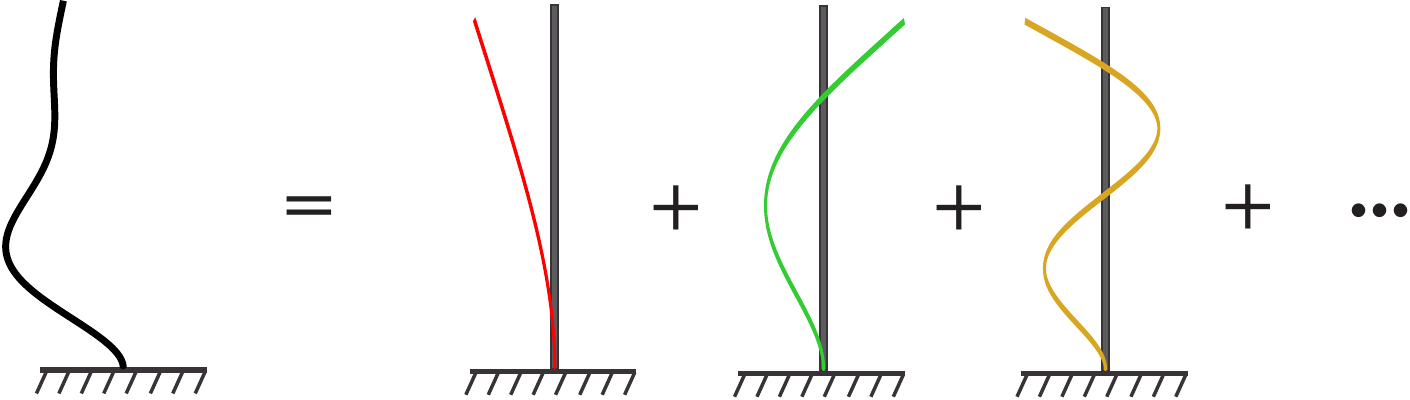}
  \caption{Superposition of the modes in modal analysis.}
  \label{fig:modal}
\end{figure}

Modal analysis is a process where the nodal displacement vector is approximated by a linear combination of dominant eigenvectors (also called mode shapes) as it is shown in Fig.~\ref{fig:modal}. Elimination of high-frequency mode shapes decreases the number of numerical operations per time~step because the size of the matrices in the equations of motion is much less than in the non-reduced case.\par
\begin{figure}[ht]
  \centering
  \footnotesize
  \def\svgwidth{0.8\columnwidth}
  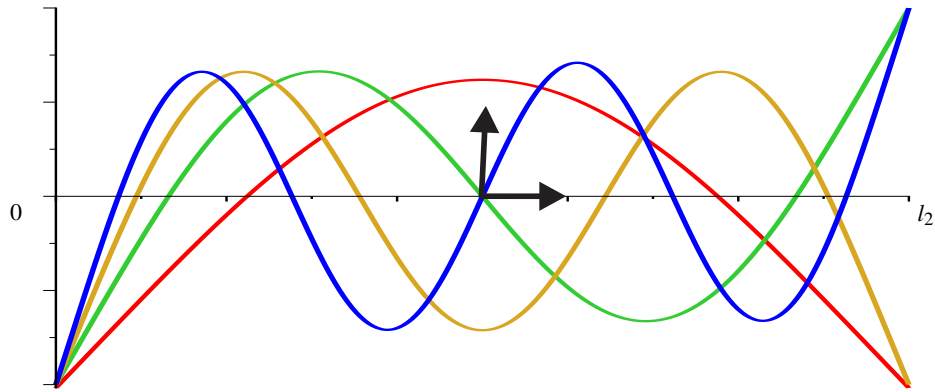
  \caption{Free-free boundary condition and new reference frame at the mass center for decoupling the reference and the elastic coordinates.}
  \label{fig:free_free}
\end{figure}

\begin{figure}[ht]
  \centering
  \footnotesize
  \def\svgwidth{0.8\columnwidth}
  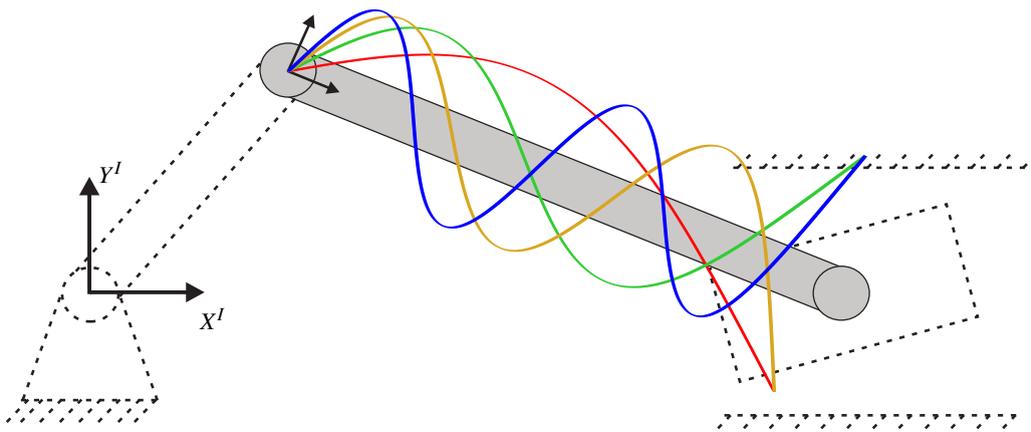
  \caption{Articulated-free boundary condition for partially decoupling the reference and the elastic coordinates.}
  \label{fig:modal_floating_frame}
\end{figure}

For the flexible beam with mass matrix $\vm_{ff}$ and stiffness matrix $\vK_{ff}$, we obtain the angular frequency $\omega_k$ and the relative mode shape $\vphi_k$ using the free vibration equations of motion. The free vibration equations of motion can be derived from \eqref{eq:mechanical_system_position} in case of no external forces and damping:
\begin{align}
  \left(\vK_{ff}-\omega^2_k~\vm_{ff}\right)\vphi_k=0\;.
\end{align}
The deformations are expressed as
\begin{align}
  q_{f_j}=\sum\limits_{k=1}^{n_m} \vPhi_{jk}q_{m_k}\;,
\end{align}
where $q_{f_j}$ is the deformation of the degree of freedom number $j$. Considering that the transformation matrix $\vPhi$ consists of the eigenvectors $\vphi_k$, the values $\vPhi_{jk}$ and $q_{m_k}$ are components of the eigenvectors and modal coordinates (the new unknowns). The number of reduced coordinates $n_m$ is chosen depending on the accuracy. We rewrite the above equation:
\begin{align}
  \vq_{f}=\vPhi~\vq_{m}\;.
\end{align}
The deformation field \eqref{eq:u_f} is described with
\begin{align}
  \vu^{ij}_f=\vS^{ij}\vPhi^{ij}\vq_{m}=\vS^{ij}_m\vq_{m} \;.\label{eq:modal_S}
\end{align}
Using characteristics of normalized orthogonal eigenvectors ($\phi_k$), we decouple the sub-matrices $\vm_{ff}$ and $\vK_{ff}$ in \eqref{eq:stiffness_matrix} and \eqref{eq:S_ff} respectively:
\begin{align}
  \vPhi^T\vm_{ff}\vPhi=\vI_{n_m \times n_m}\;,\quad \vPhi^T\vK_{ff}\vPhi=\omega^2_k~\vdelta_{kl}\;.
\end{align}
In order to complete the decoupling in the mass matrix, we deal with the sub-matrices $\vm_{R\theta}$, $\vm_{Rf}$ and $\vm_{\theta f}$ on element level using the modified shape functions introduced in \eqref{eq:modal_S}.  Moving the body coordinate system to the mass center, we show that the integrals in \eqref{eq:S_bar} and \eqref{eq:S_tilde} vanish if
\begin{align}
  \int_{V}{\phi}_k \left(\vr\right) \mathrm{d} V =0 \;,\quad \int_{V}\vr~ {\phi}_k \left(\vr\right) \mathrm{d} V=\vnull \;.\label{eq:free_cond}
\end{align}
It can be shown that the free-free modes satisfy the condition above and also the mean-axis conditions which is obtained by minimizing the kinetic energy of the elastic motion with respect to an observer sitting on the flexible body~\cite{Agr86}. This ideal coordinate system and the free-free eigenfunctions are shown in Fig.~\ref{fig:free_free}. Removing the first three modal coordinates (corresponding to the rigid body motion) results in decoupled versions of the total mass matrix $\vM$.
It is clear from Fig.~\ref{fig:free_free} that the deformation at the two ends of the beam do not vanish as defined in this coordinate system. Shabana~\cite{Sha96} modifies the free-free shape functions to satisfy the boundary condition which results in a simply supported, i.e., pinned, mode. However, it cannot be used to decouple the mass matrix anymore. Using articulated-free modes as shown in Fig.~\ref{fig:modal_floating_frame}~\cite{Esc02}, i.e., fixing the position of the left end and leaving the right end free, we can fix the boundary condition at the joint in addition to satisfying the second condition in \eqref{eq:free_cond}, which results in vanishing terms regarding $\vm_{\theta f}$.

\bibliographystyle{plain}
\bibliography{Literatur}

\end{document}